\documentclass[10pt,reqno]{article}
\usepackage[a4paper]{anysize}\marginsize{3cm}{3cm}{2cm}{2cm}
\pdfpagewidth=\paperwidth \pdfpageheight=\paperheight
\usepackage[T1]{fontenc}
\usepackage[english]{babel}
\usepackage{microtype}
\usepackage{amssymb,amsmath,amsthm,braket}
\usepackage{mathtools}
\usepackage{verbatim}
\usepackage{relsize}
\usepackage[colorlinks=true, urlcolor=blue, linkcolor=blue, citecolor=blue]{hyperref}
\usepackage[percent]{overpic}
\usepackage{comment}
\usepackage{algorithm}
\usepackage{algpseudocode}
\usepackage{xcolor}
\definecolor{OliveGreen}{rgb}{0,0.6,0}
\usepackage{caption}
\usepackage{subcaption}
\usepackage{enumitem}
\usepackage{wrapfig}
\newlist{CI}{enumerate}{1}
\setlist[CI]{label=(C\arabic*)}
\usepackage{pgfplots}
\pgfplotsset{compat=1.15}
\usetikzlibrary{arrows}
\usepackage{authblk}
\usepackage{diagbox}
\usepackage{array}
\usepackage{cite}

\numberwithin{equation}{section}
\theoremstyle{plain}
\newtheorem*{theorem*}{Theorem}
\newtheorem{theorem}{Theorem}
\numberwithin{theorem}{section}
\newtheorem{proposition}[theorem]{Proposition}
\newtheorem{lemma}[theorem]{Lemma}

\newtheorem{conjecture}[theorem]{Conjecture}

\theoremstyle{definition}
\newtheorem{definition}[theorem]{Definition}
\newtheorem{remark}[theorem]{Remark}
\newtheorem{example}[theorem]{Example}

\newcommand{\sS}{\mathcal{S}}

\newcommand{\sP}{\mathcal{P}}

\DeclareMathOperator*{\sym}{Sym}

\graphicspath{ {figures/} }

\makeatletter
\newcommand{\ALOOP}[1]{\ALC@it\algorithmicloop\ #1%
  \begin{ALC@loop}}
\newcommand{\ENDALOOP}{\end{ALC@loop}\ALC@it\algorithmicendloop}

\algnewcommand\TRUE{\textbf{true}\space}
\algnewcommand\FALSE{\textbf{false}\space}
\makeatother
\usepackage{float}

\makeatletter
\algnewcommand{\LineComment}[1]{\Statex \hskip\ALG@thistlm \(\triangleright\) #1}
\makeatother

\makeatletter
\newcommand\fs@nocaptionruled{
  \let\@fs@capt\floatc@ruled
  \def\@fs@pre{}
  \def\@fs@post{\kern2pt\hrule\relax}
  \def\@fs@mid{}
  \let\@fs@iftopcapt\iftrue}
\makeatother

\makeatletter
\newcommand\fs@nobottomruled{
\def\@fs@cfont{\bfseries}
\captionsetup[algorithm]{labelfont=bf,justification=raggedright,singlelinecheck=false,labelsep=space}
 \def\@fs@pre{\hrule height.8pt depth0pt \kern2pt}%
 \def\@fs@post{}
 \def\@fs@mid{\kern2pt\hrule\kern2pt}%
\let\@fs@iftopcapt\iftrue
  }
\makeatother

\makeatletter
\newcommand\fs@nobottomandcaptionruled{\def\@fs@cfont{\bfseries}\let\@fs@capt\floatc@ruled
 \let\@fs@capt\relax
 \def\@fs@pre{}
 \def\@fs@post{}
 \def\@fs@mid{}%
 \let\@fs@iftopcapt\iftrue}
\makeatother

\def\NoNumber#1{{\def\alglinenumber##1{}\State #1}\addtocounter{ALG@line}{-1}}

\DeclareMathOperator{\an}{an}

\DeclareMathOperator{\cycle}{\mathcal{SC}}
\DeclareMathOperator{\score}{GS}
\DeclareMathOperator{\moral}{Moral}
\DeclareMathOperator{\Alg}{Alg}
\DeclareCaptionLabelFormat{andtable}{#1~#2  \&  \tablename~\thetable}
\newcommand{\gtp}{\ensuremath{\delta^+_{\tilde{A}_C \setminus Y}}}
\newcommand{\gtm}{\ensuremath{\delta^-_{\tilde{A}_C \setminus Y}}}
\newcommand{\gtjp}{\ensuremath{\delta^+_{\tilde{A}^{Y,J}_C}}}

\newcommand{\eyj}{\ensuremath{\tilde{E}^{Y,J}_C}}
\newcommand{\eyjx}{\ensuremath{\tilde{E}^{Y,J(x)}_C}}
\newcommand{\ayj}{\ensuremath{\tilde{A}^{Y,J}_C}}
\newcommand{\by}{\ensuremath{B^{Y}_C}}

\newtheorem*{claim*}{Claim}
\newcommand{\padj}{\mathrel{\stackrel{p}{\sim}}}

\newcommand\independent{\protect\mathpalette{\protect\independenT}{\perp}}
\def\independenT#1#2{\mathrel{\rlap{$#1#2$}\mkern2mu{#1#2}}}

\title{Causal Structure Learning in Directed, Possibly Cyclic, Graphical Models}
\author{Pardis Semnani and Elina Robeva}
\affil{University of British Columbia}
\date{}

\begin{document}

\maketitle

\begin{abstract}
We consider the problem of learning a directed graph $G^\star$ from observational data. We assume that the distribution which gives rise to the samples is Markov and faithful to the graph $G^\star$ and that there are no unobserved variables. We do not rely on any further assumptions regarding the graph or the distribution of the variables. Particularly, we allow for directed cycles in $G^\star$ and work in the fully non-parametric setting. Given the set of conditional independence statements satisfied by the distribution, we aim to find a directed graph which satisfies the same $d$-separation statements as $G^\star$. We propose a hybrid approach consisting of two steps. We first find a {\em partially ordered partition} of the vertices of $G^\star$ by optimizing a certain score in a greedy fashion. We prove that any optimal partition uniquely characterizes the Markov equivalence class of $G^\star$. Given an optimal partition, we propose an algorithm for constructing a graph in the Markov equivalence class of $G^\star$ whose strongly connected components correspond to the elements of the partition, and which are partially ordered according to the partial order of the partition. Our algorithm comes in two versions -- one which is provably correct and another one which performs fast in practice.

\vspace{5mm}
\noindent{\em Keywords:} causal discovery, directed cyclic graphs, Markov equivalence, faithfulness, $d$-separations, conditional independence

\vspace{5mm}
\noindent{\em 2020 Mathematics Subject Classification:} 62D20 (Primary), 62H22, 05C20 (Secondary)
\end{abstract}

\section{Introduction} \label{casuality: introduction}

\begin{figure}[t]
    \centering
    \includegraphics[width=0.5\linewidth]{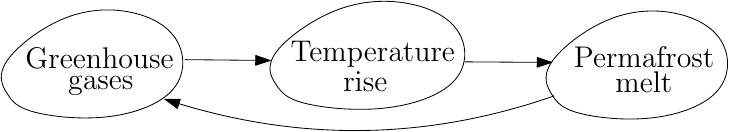}
    \caption{An example of a causal cycle.}
    \label{fig:cycle_example}
\end{figure}
Determining the causal structure between several random variables from observational data is a central task in many  disciplines including computational biology, epidemiology, sociology, and economics~\cite{Friedman2000,pearlcausalitymodels,Robins2000, SpirtesGlymourClark1993}. The causal structure is often modeled by a directed graph, where the vertices correspond to the variables of interest, and the directed edges represent the direct causal effects those variables have on one another. 
One of the most common simplifying assumptions is that the directed graph is acyclic, i.e., there are no directed cycles (or feedback loops). However, in many realistic settings, directed cycles do exist (e.g., gene regulatory networks, climate science (see Figure~\ref{fig:cycle_example}), social sciences, feedback systems in electrical engineering, and economic processes \cite{Mason1953, Mason1956, Haavelmo1943, Goldberger1972}).

In this work, we propose a new method for causal discovery in case the causal graph $G^\star$ may contain directed cycles. While most approaches for causal discovery from observational data are either score-based (e.g., \cite{Chickering2002OptimalSI, meek1997graphical}) or constraint-based (e.g., \cite{Verma1990EquivalenceAS, SpirtesGlymourClark1993}), we here extend the line of work based on a hybrid approach~\cite{Teyssier2005OrderingBasedSA, Raskutti2018, Solus2021,Wang2017, mohammadi2017generalized, Bernstein2020} from the acyclic to the cyclic setting. To the best of our knowledge, our algorithm is the first  which works in the fully non-parametric setting and does not assume the graph is acyclic. \cite{richardson2013discovery} and \cite{SAT2013} both propose causal discovery algorithms that work for cyclic graphs with no parametric assumptions, however, they do not produce a graph Markov equivalent to $G^\star$. Rather, they output some of the characteristics that are shared by all of the members of the Markov equivalence class of $G^\star$. Our algorithm not only specifies the Markov equivalence class of $G^\star$, but also outputs a graph Markov equivalent to $G^\star$.

In recent work~\cite{claassen23a}, Claassen and Mooij introduce an efficient algorithm for determining Markov equivalence between directed, possibly cyclic, graphs. Given a directed graph $G$, the algorithm generates a \textit{partial ancestral graph} that represents the invariant ancestral relationships unique to the Markov equivalence class of $G$. Unlike causal discovery algorithms, this algorithm directly uses the graph structure itself -- not the conditional independence statements satisfied by a distribution Markov to the graph, for instance. However, the authors suggest that their characterization of Markov equivalence could inform future causal discovery methods.

While there does not necessarily exist a topological ordering of the variables in a possibly cyclic directed graph, the set of {\em strongly connected components} forms a partially ordered partition of the vertices of the graph.

\begin{definition}\label{defn:strong_component}
Let $G=(V,E)$ be a directed graph. A subset $C\subseteq V$ is a \textit{strongly connected component} of $G$ if $C$ is a maximal set of vertices such that every vertex in $C$ can be reached from all other vertices in $C$ via a directed path. 

We denote the set of strongly connected components of $G$ by $\cycle(G)$, and the unique strongly connected component containing the vertex $v$ in graph $G$ by $C_{v,G}$ for all $v\in V$. 
\end{definition}

\begin{figure}[t]
    \centering
    \includegraphics[scale=0.65]{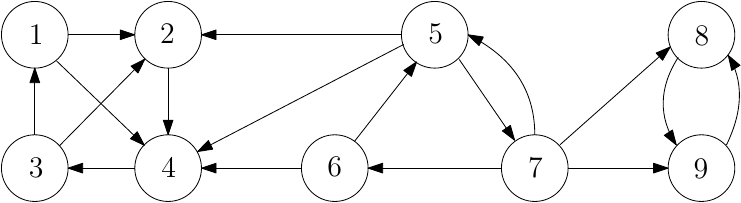}
    \caption{An example of a cyclic graph with strongly connected components $\{1,2,3,4\}, \{5, 6, 7\}, \{8,9\}$ obeying the partial order: $\{5,6,7\}\le_G \{1,2,3,4\}; \{5,6,7\}\le_G \{8,9\}$.
    }
    \label{fig:cyclic_graph_example}
\end{figure}

For example, in Figure~\ref{fig:cyclic_graph_example}, the strongly connected components are $\{1,2,3,4\}, \{5, 6, 7\}, \{8,9\}$.

\begin{definition}
For a directed graph $G=(V,E)$, 
consider the following relation on $\cycle(G)$: For all $C_1,C_2\in \cycle(G)$,
\begin{align*}
    C_1 \leq_G C_2 \Longleftrightarrow &\text{ There is a directed path from any vertex in $C_1$ to any vertex in $C_2$ in $G$,}\\
    &\text{ or $C_1=C_2$.}
\end{align*}
One can verify that $\leq_G$ is a partial order on $\cycle(G)$. 
The pair $(\cycle(G),\leq_G)$ is called \textit{the partially ordered partition associated with $G$}.
\end{definition}

For example, the partially ordered partition associated with the graph $G$ in Figure~\ref{fig:cyclic_graph_example} is the partition $\{1,2,3,4\}, \{5, 6, 7\}, \{8,9\}$ with partial order $\{5,6,7\}\leq_{G} \{1,2,3,4\}$ and $\{5,6,7\}\leq_{G}\{8,9\}$.

Given a set of conditional independence statements satisfied by our unknown distribution, we aim to recover the Markov equivalence class of the causal graph of this distribution as well as a representative of this class. In the first step of our algorithm we greedily find a partially ordered partition of the vertices by optimizing a certain sparsity-based score. This partition, along with the optimized score, completely characterize the Markov equivalence class of the causal graph  (see Section~\ref{poset discovery}). Using this partially ordered partition, in the second step of the algorithm we then build a directed graph which is Markov equivalent to the true causal graph (see Section~\ref{graph discovery}).

In the present work, we assume the global Markov property and faithfulness.
A distribution satisfies the global Markov property with respect to $G^\star$ if every $d$-separation statement arising from $G^\star$ is a conditional independence statement satisfied by the distribution. The distribution is faithful to $G^\star$ when every conditional independence statement satisfied by the distribution is also a $d$-separation statement arising from $G^\star$. Therefore, in this work we assume the conditional independence statements satisfied by the observed distribution are precisely those arising from the $d$-separations of the graph.

  In our experiments (Section~\ref{simulations}) we observe that our causal discovery algorithm performs well on graphs of up to 10 vertices. 
The algorithm takes as input the set of conditional independence statements satisfied by the unknown distribution. In our simulations, this set is generated by deriving the $d$-separation statements from a fixed graph. However, producing this set becomes prohibitively slow for graphs with more than 10 vertices, which is the factor that has limited our testing to smaller graphs. We wish to remark that in practice we only need access to a subset of the conditional independence statements. Therefore, we will not need to test for all of them when working with real data (see Remark~\ref{remark: number of tests}).

Note that in a structural equation model (SEM) corresponding to an acyclic graph, $d$-separation in the graph implies conditional independence in the observational distribution of the SEM (i.e., the distribution of the SEM is Markov to the acyclic graph). However, for SEMs corresponding to cyclic graphs, this is only known to be true when the SEM is linear or discrete~\cite{Spirtes1995,Bongers_2021}. A generalization of the notion of $d$-separation, called $\sigma$-separation, has been introduced in~\cite{forré2017markov}, and is shown to imply conditional independence under certain unique solvability conditions for the SEM~\cite[Theorem 6.3]{Bongers_2021}. Furthermore,~\cite{Bongers_2021} proves that under certain ancestral unique solvability conditions the causal interpretation of an SEM is consistent with its corresponding graph~\cite[Propositions 7.1 and 7.2]{Bongers_2021}.

\textbf{Organization.} 
The remainder of this paper is organized as follows. Section~\ref{characterization} gives an overview of the characterization of Markov equivalence for directed cyclic graphs proposed in \cite{richardson1997char}. In Section~\ref{poset discovery} we present our greedy algorithm for discovering the Markov equivalence class of the unknown graph. In Section~\ref{graph discovery}, we propose an algorithm which finds a graph in this Markov equivalence class. In Section~\ref{simulations} we show simulations illustrating the efficiency of our algorithm. We conclude with a discussion in Section~\ref{discussion}.

\section{Characterization of Markov Equivalence in Directed Cyclic Graphs}\label{characterization}
We assume that the input to our causal discovery algorithm is the set of  conditional independence statements arising from the observed distribution, and that these correspond precisely to the $d$-separation statements that the unknown graph $G^\star$ satisfies. 
We aim to recover a graph satisfying the same $d$-separation statements as the true graph $G^\star$, i.e., a graph {\em Markov equivalent} to $G$.

A known characterization of Markov equivalence for DAGs, proved in \cite{Frydenberg1990, Verma1990EquivalenceAS}, is that two DAGs are Markov equivalent if and only if they have the same skeleton (i.e., the same underlying undirected graph) and the same set of immoralities (i.e., the same set of triples $i\to j\leftarrow k$ where $i$ and $k$ are not adjacent). However, this statement does not generally hold for all directed graphs.

\begin{example}\label{example: Markov equivalence doesn't work for DAG}
The graphs shown in Figure~\ref{fig3: non Markov equivalent graphs} have the same skeleton and $a \rightarrow b \leftarrow c$ is the only immorality in both of the graphs. However, $a$ and $c$ are $d$-connected in the left graph via the path $a \rightarrow b \rightarrow d \rightarrow c$ while $a$ and $c$ are $d$-separated in the right graph.

On the other hand, in Figure~\ref{fig4: Markov equivalent graphs}, no $d$-separation statement arises from either of the graphs, which means that they are Markov equivalent. Nevertheless, they do not have the same skeleton or the same set of immoralities.
\begin{figure}[h]
    \centering
    \begin{subfigure}{1\textwidth}
        \centering
        \includegraphics[scale=0.6]{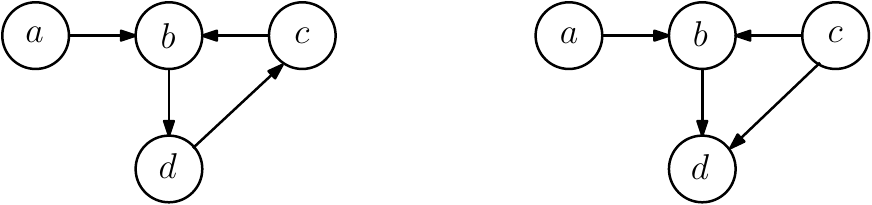}
        \caption{Two graphs with the same skeleton and the same set of immoralities that are not Markov equivalent.}
        \label{fig3: non Markov equivalent graphs}
    \end{subfigure}
    \hfill
    \begin{subfigure}{1\textwidth}
        \centering
        \includegraphics[scale=0.65]{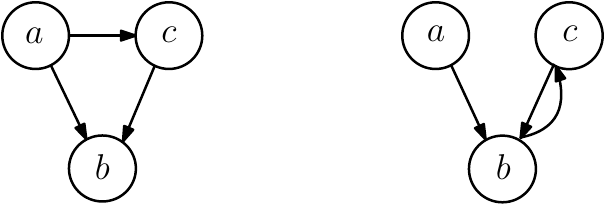}
        \caption{Two Markov equivalent graphs with different skeletons and immoralities.}
        \label{fig4: Markov equivalent graphs}
    \end{subfigure}
\caption{The local characterization of Markov equivalence for DAGs does not hold in general for all directed graphs.}
\end{figure}
\end{example}

\begin{figure}
    \centering
    \includegraphics[width=0.5\textwidth]{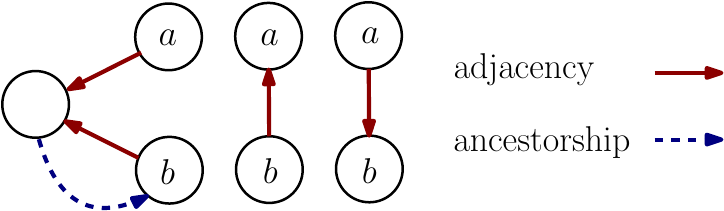}
    \caption{A schematic representation of the three conditions any of which can make two vertices $a$ and $b$ $p$-adjacent in a directed graph.}
    \label{fig:p-adjacency}
\end{figure}

Having a local characterization of Markov equivalence for all directed graphs is essential to designing an algorithm for causal discovery. Such a characterization is given in \cite{richardson1997char}. We state it in Theorem~\ref{th: Markov equivalence characterization} after providing the necessary background in the rest of this section.
\begin{definition}
Let $G=(V,E)$ be a directed graph. Two vertices $a,b\in V$ are said to be \textit{pseudo-adjacent} or for short, \textit{$p$-adjacent} if either there is an edge between $a$ and $b$ in $G$, or $a$ and $b$ have a common child in $G$ which is also an ancestor of $a$ or $b$ in $G$. See Figure~\ref{fig:p-adjacency}.
\end{definition}
Note that  the notion of $p$-adjacency is a generalization of the notion of adjacency in DAGs. More precisely, while in a DAG, two vertices $a$ and $b$ are adjacent if and only if $a$ is $d$-connected to $b$ given $S$ for all $S\subseteq V\setminus \{a,b\}$, \cite{richardson1997char} proves that 
two vertices $a,b\in V$ are $p$-adjacent if and only if for all sets $S\subseteq V\setminus \{a,b\}$, $a$ is $d$-connected to $b$ given $S$.

\begin{figure}
    \centering
    \begin{subfigure}[t]{0.3\textwidth}
    \centering
        \includegraphics[scale=0.6]{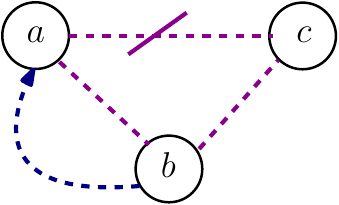}
        \caption{Unshielded conductor.}
    \end{subfigure}
    \hfill
    \begin{subfigure}[t]{0.3\textwidth}
    \centering
        \includegraphics[scale=0.6]{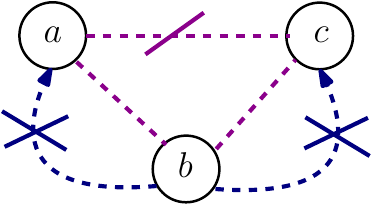}
        \caption{Unshielded non-conductor.}
    \end{subfigure}
    \hfill
    \begin{subfigure}[t]{0.2\textwidth}
    \centering
        \includegraphics[scale=0.6]{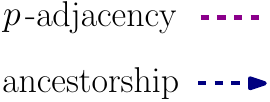}
    \end{subfigure}
    \caption{A schematic representation of a triple $(a,b,c)$ as an unshielded conductor and an unshielded non-conductor in a directed graph.}
    \label{fig: unshielded conductor}
\end{figure}

\begin{definition}
Let $G=(V,E)$ be a directed graph. For a triple $(a,b,c)\in V^3$, suppose $a$ and $c$ are not $p$-adjacent in $G$, but both the pairs $a,b$ and $c,b$ are $p$-adjacent in $G$. The triple $(a,b,c)$ is said to be an \textit{unshielded conductor} if $b$ is an ancestor of $a$ or $c$ in $G$, and is said to be an \textit{unshielded non-conductor} otherwise. See Figure~\ref{fig: unshielded conductor}.
\end{definition}
The notion of unshielded non-conductors is a generalization of the notion of immoralities in DAGs. 

\begin{figure}
    \centering
    \begin{subfigure}[t]{0.3\textwidth}
    \centering
        \includegraphics[scale=0.6]{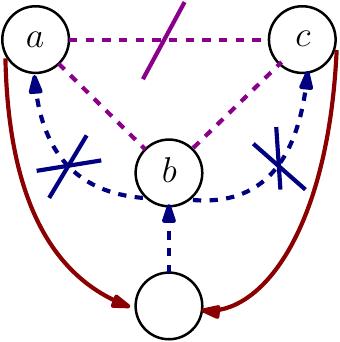}
        \caption{Unshielded perfect non-conductor.}
    \end{subfigure}
    \hfill
    \begin{subfigure}[t]{0.3\textwidth}
    \centering
        \includegraphics[scale=0.6]{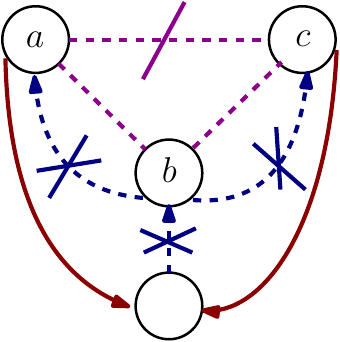}
        \caption{Unshielded imperfect non-conductor.}
    \end{subfigure}
    \hfill
    \begin{subfigure}[t]{0.2\textwidth}
    \centering
        \includegraphics[scale=0.6]{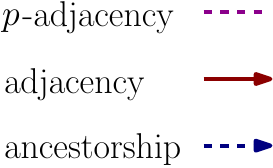}
         \vspace{10mm}
    \end{subfigure}
    \hfill
    \caption{A schematic representation of a triple $(a,b,c)$ as an unshielded perfect non-conductor and an unshielded imperfect non-conductor in a directed graph.}
    \label{fig: unshielded perfect conductor}
\end{figure}

\begin{definition}
Let $G=(V,E)$ be a directed graph and $(a,b,c)\in V^3$ be an unshielded non-conductor in $G$. Then  $(a,b,c)$ is said to be an \textit{unshielded perfect non-conductor} if $b$ is a descendant of a common child of $a$ and $c$ in $G$, and is said to be an \textit{unshielded imperfect non-conductor} otherwise. See Figure~\ref{fig: unshielded perfect conductor}.
\end{definition}
In a DAG, since all unshielded non-conductors are immoralities, they are all perfect.

\begin{figure}
    \centering
    \includegraphics[scale=0.6]{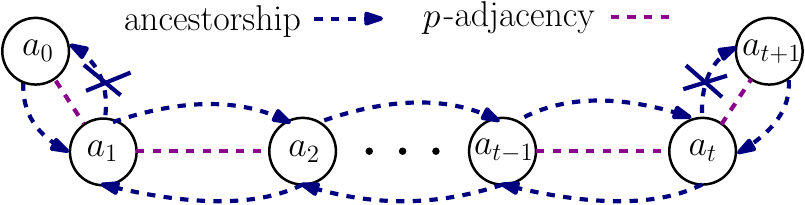}
    \caption{A schematic representation of triples $(a_0,a_1,a_2)$ and $(a_{t-1},a_t,a_{t+1})$ which are mutually exclusive with respect to the uncovered itinerary $(a_0,a_1,\ldots,a_t,a_{t+1})$ in a directed graph.}
    \label{fig: mutually exclusive}
\end{figure}

\begin{figure}
    \centering
    \begin{subfigure}[t]{0.48\textwidth}
    \centering
        \includegraphics[scale=0.6]{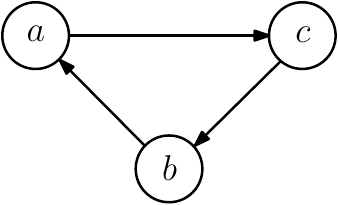}
        \caption{A cyclic graph with no $d$-separation statements.}
    \end{subfigure}
    \hfill
    \begin{subfigure}[t]{0.45\textwidth}
    \centering
        \includegraphics[scale=0.6]{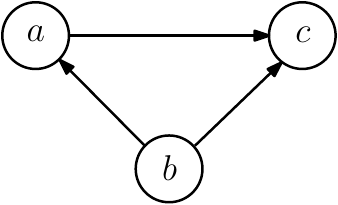}
        \caption{A DAG with no $d$-separation statements.}
    \end{subfigure}
    \hfill
    \caption{A DAG and a directed cyclic graph that are Markov equivalent.}
    \label{fig: MEC with DAG and cyclic}
\end{figure}

\begin{figure}
    \centering
    \includegraphics[scale=0.7]{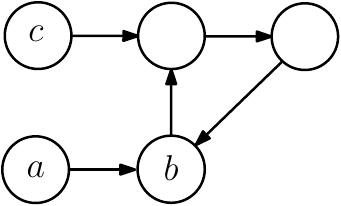}
    \caption{A directed cyclic graph that is not Markov equivalent to any DAGs.}
    \label{fig: MEC without DAG}
\end{figure}

\begin{definition}\label{def: mutually exclusive}
        Let $G=(V,E)$ be a directed graph. For $t\in \mathbb{N}$ and $a_0,a_1,\ldots,a_{t},a_{t+1}\in V$, it is said that the triples $(a_0,a_1,a_2),(a_{t-1},a_t,a_{t+1})$ are \textit{mutually exclusive with respect to the uncovered itinerary $(a_0,a_1,\ldots,a_t,a_{t+1})$} if
        \begin{itemize}
                \item 
                $a_{i}$ and $a_{i-1}$ are $p$-adjacent for all $i\in [t+1]$,
                \item 
                $a_{i}$ and $a_{j}$ are not $p$-adjacent for all $i\in \{2,\ldots,t+1\}$ and $j\in \{0,\ldots,i-2\}$,
                \item 
                $a_i$ is an ancestor of $a_{i+1}$ for all $i\in \{0,\ldots,t-1\}$,
                \item 
                $a_i$ is an ancestor of $a_{i-1}$ for all $i\in \{2,\ldots,t+1\}$, and
                \item 
                $a_1$ is not an ancestor of $a_0$ or $a_{t+1}$.
        \end{itemize}
See Figure~\ref{fig: mutually exclusive}.

\end{definition}
Definition~\ref{def: mutually exclusive} is stated slightly differently in \cite{richardson1997char} in the sense that $t$ is considered to be greater than 1 in \cite{richardson1997char}. However, one can confirm this difference does not affect the statement of Theorem~\ref{th: Markov equivalence characterization} since when $t=1$ in Definition~\ref{def: mutually exclusive}, $(a_0,a_1,a_2)$ is simply an unshielded non-conductor in $G$. We are now prepared to state this theorem, which is the main result in \cite{richardson1997char}.
\begin{theorem}\label{th: Markov equivalence characterization}
Assume that $G_1=(V,E_1)$ and $G_2=(V,E_2)$ are two directed graphs. Then $G_1$ and $G_2$ are Markov equivalent if and only if the following conditions hold:
\begin{enumerate}
    \item $G_1$ and $G_2$ have the same $p$-adjacencies. \label{thm: markov con1}
    \item $G_1$ and $G_2$ have the same set of unshielded conductors.\label{thm: markov con2}
    \item $G_1$ and $G_2$ have the same set of unshielded perfect non-conductors.\label{thm: markov con3}
    \item If $(a,b_1,c)$ and $(a,b_2,c)$ are unshielded imperfect non-conductors (in $G_1$ and $G_2)$, then $b_1$ is an ancestor of $b_2$ in $G_1$ if and only if $b_1$ is an ancestor of $b_2$ in $G_2$.\label{thm: markov con4}
    \item For any $t\in \mathbb{N}$, triples $(a_0,a_1,a_2)$ and $(a_{t-1},a_t,a_{t+1})$ are mutually exclusive with respect to the uncovered itinerary $P=(a_0,a_1,\ldots,a_t,a_{t+1})$ in $G_1$ if and only if $(a_0,a_1,a_2)$ and $(a_{t-1},a_t,a_{t+1})$ are mutually exclusive with respect to $P$ in $G_2$.
    \label{thm: markov con5}
    \item If $(a_0,a_1,a_2)$ and $(a_{t-1},a_t,a_{t+1})$ are mutually exclusive with respect to the uncovered itinerary $(a_0,a_1,\ldots,a_{t+1})$ for some $t\in \mathbb{N}$ and $(a_0,b,a_{t+1})$ is an unshielded imperfect non-conductor (in $G_1$ and $G_2$), then $a_1$ is an ancestor of $b$ in $G_1$ if and only if $a_1$ is an ancestor of $b$ in $G_2$. \label{thm: markov con6}
\end{enumerate}
\end{theorem}
\begin{remark}
    Some Markov equivalence classes contain both DAGs and cyclic graphs. For example, the two graphs in Figure~\ref{fig: MEC with DAG and cyclic} are Markov equivalent since no $d$-separation statement is induced by these graphs. On the other hand, there are Markov equivalence classes which only contain directed cyclic graphs. For instance, the graph in Figure~\ref{fig: MEC without DAG} cannot be Markov equivalent to any DAG since $(a,b,c)$ is an unshielded imperfect non-conductor in this graph, and according to Theorem~\ref{th: Markov equivalence characterization}, this triple is an unshielded imperfect non-conductor in any Markov equivalent graph. However, all unshielded non-conductors are perfect in DAGs.
\end{remark}

\section{Discovering the Markov equivalence class}\label{poset discovery}
We assume that we have observed i.i.d. samples from a distribution $\mathbb P$ which is both Markov and faithful to a directed graph $G^\star=(V^\star,E^\star)$ and have inferred all the conditional independence statements satisfied by this distribution. In this section we try to recover the Markov equivalence class of $G^\star$. We assume throughout that there are no latent variables, and that $V^\star=\{1,\ldots,n\}$, where $n$ is the number of observed variables.

To recover the Markov equivalence class of $G^\star$, we use a greedy procedure that moves through the set of partially ordered partitions of the vertices of the graph and optimizes a certain score. We sweep through the set of partially ordered partitions in a depth-first search manner in order to greedily optimize the score. Our score is sparsity-based and is given by a vector with $n+2$ entries, where $n$ is the number of variables. We compare scores lexicographically. Given a partially ordered partition and the input of conditional independence statements, the score first counts the number of $p$-adjacencies. Once we have found a partially ordered partition with a minimal number of $p$-adjacencies, we then minimize a different quantity based on the characterization of Markov equivalence, etc. The entries of the score correspond to the different characteristics that Markov equivalent graphs share as described in~\ref{th: Markov equivalence characterization}. We show that if a partially ordered partition is a global minimizer of our score, then it characterizes the Markov equivalence class of $G^\star$ (Theorem~\ref{thm:minimal_score_graph}).

\subsection {Defining a score} \label{subsection: defining the score}
In order to recover the Markov equivalence class of $G^\star$, we assign a certain score, described in Definition~\ref{def: score}, to each partially ordered partition on $V^\star$. We prove in Theorem~\ref{thm:minimal_score_graph} that the partially ordered partitions whose score is optimal characterize the Markov equivalence class of $G^\star$. 
\begin{definition}
Define
\begin{align*}
\sS \coloneqq \Set{(\sP,\pi) | \text{$\sP\subseteq 2^{V^\star}$ is a partition of $V^\star$ and $\pi\subseteq \sP\times \sP$ is a partial order on $\sP$}}.
\end{align*}
For all $i\in V^\star$, we denote the unique section of the partition $\sP$ which contains $i$ by $C_{i,\sP}$. For $(\sP,\pi)\in \mathcal{S}$ and for any $C_1,C_2\in \sP$, we write $C_1 \leq_\pi C_2$ whenever $(C_1,C_2)\in \pi$.
\end{definition}

\begin{definition} \label{def: score}
For a given set of conditional independence statements and for any partially ordered partition $(\sP,\pi)\in \mathcal{S}$, we define the \textit{graphical score} of $(\sP,\pi)$, denoted by $\score(\sP,\pi)$, to be
\begin{align*}
    \score (\sP,\pi) \coloneqq \left(\ \left|E^{(1)}_{(\sP,\pi)}\right|,\left|E^{(2)}_{(\sP,\pi)}\right|,\left|E^{(3)}_{(\sP,\pi)}\right|,\left|E^{(4)}_{(\sP,\pi)}\right|,-\left|D^{(2)}_{(\sP,\pi)}\right|,\cdots,-\left|D^{(n-2)}_{(\sP,\pi)}\right|,\left|E^{(6)}_{(\sP,\pi)}\right| \ \right),
\end{align*}
where the sets $E^{(1)}_{(\sP,\pi)}, E^{(2)}_{(\sP,\pi)}, E^{(3)}_{(\sP,\pi)}, E^{(4)}_{(\sP,\pi)},  D^{(t)}_{(\sP,\pi)}$ for $t\in \{2,\ldots,n-2\}$, and $E^{(6)}_{(\sP,\pi)}$ are defined in Propositions~\ref{prop: Rich1},~\ref{prop: Rich2},~\ref{prop: Rich3},~\ref{prop: Rich4},~\ref{prop: Rich5}, and~\ref{prop: Rich6} respectively.

\end{definition}
The graphical score of any partially ordered partition $(\sP,\pi)\in \mathcal{S}$ is a vector in $\mathbb{Z}^{n+2}$. Equipping $\mathbb{Z}^{n+2}$ with the \textit{lexicographical order} allows us to compare the graphical scores of different partially ordered partitions. Our main result in this subsection is the following.
\begin{theorem}\label{thm:minimal_score_graph}
Each minimizer of the graphical score over $\mathcal{S}$,  uniquely determines the Markov equivalence class of $G^\star$.
\end{theorem}
\begin{figure}
    \centering
    \begin{subfigure}[b]
    {0.5\textwidth}
    \centering
    \includegraphics[width= 0.45\textwidth]{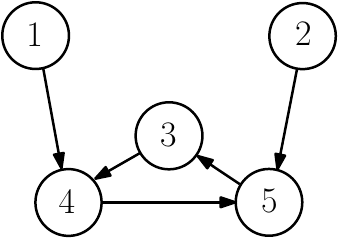}
    \caption{Graph $G^\star$.}
    \label{fig: running example 3.1 graph}
    \end{subfigure}
    ~
    \begin{subfigure}[b]{0.42\textwidth}
    \centering
    \begin{align*}
        &1 \independent 2 \mid \{\} && 1 \independent 5 \mid \{3,4\} \\
        & 1 \independent 2 \mid \{3,4\} && 1 \independent 5 \mid \{2,3,4\} \\
        & 1 \independent 2 \mid \{4,5\} && 2 \independent 3 \mid \{4,5\} \\
        & 1 \independent 2 \mid \{3,4,5\} && 2 \independent 3 \mid \{1,4,5\}
    \end{align*}
    \caption{The set of $d$-separation statements of $G^\star$.}
    \label{fig: running example section 3.1 statements}
    \end{subfigure}
    \caption{A directed graph and the set of $d$-separation statements arising from it.}
    \label{fig: running example section 3.1}
\end{figure}

\begin{example}\label{ex: graphical score}
    Consider the graph $G^\star$ depicted in Figure~\ref{fig: running example 3.1 graph}. The set of $d$-separations arising from $G^\star$ are shown in Figure~\ref{fig: running example section 3.1 statements}. If the observed distribution is Markov and faithful to $G^\star$, then these statements represent the complete set of conditional independence statements satisfied by the distribution. For the purpose of a running example in this subsection, we examine the following partially ordered partitions, where only the non-trivial pairs in each partial order are listed below:
    \begin{align*}
        & \mathcal{P}_1  = \{\{1\}, \{2\}, \{3,4,5\} \}, &&\{1\}\leq_{\pi_1} \{3,4,5\}, \ \{2\} \leq_{\pi_1} \{3,4,5\},\\
        & \mathcal{P}_2 = \{\{1,2\},\{3,4,5\}\}, && \{1,2\} \leq_{\pi_2} \{3,4,5\},\\
        & \mathcal{P}_3  = \{\{1,2\},\{3,4\},\{5\} \}, && \{1,2\} \leq_{\pi_3} \{3,4\},\\
        & \mathcal{P}_4  = \{\{1,2\},\{3,4,5\}\}, && \{3,4,5\} \leq_{\pi_4} \{1,2\}.
        \end{align*}
        Using Definition~\ref{def: score}, the graphical scores of these partially ordered partitions are as follows:
        \begin{align*}
            &\score(\mathcal{P}_1, \pi_1) = (7,4,0,1,-1,0,2),
            &&\score(\mathcal{P}_2, \pi_2) = (7,4,0,1,-1,0,2),\\
            &\score(\mathcal{P}_3, \pi_3) = (9,0,1,4,0,0,4),
            &&\score(\mathcal{P}_4, \pi_4) = (7,5,0,0,0,0,0).
        \end{align*}
        Computation of the graphical scores corresponding to all the possible partially ordered partitions on $V^\star = [5]$ reveals that $\score(\mathcal{P}_1, \pi_1) = \score(\mathcal{P}_2, \pi_2) = \min_{(\sP,\pi)\in \mathcal{S}} \score(\sP, \pi)$, which according to Theorem~\ref{thm:minimal_score_graph}, implies that $(\mathcal{P}_1,\pi_1)$ and $(\mathcal{P}_2, \pi_2)$ determine the Markov equivalence class of $G^\star$. More precisely, for $i\in [2]$, the sets 
        \begin{align} \label{eq: sets in the graphical score}
            E^{(1)}_{(\sP_i,\pi_i)}, E^{(2)}_{(\sP_i,\pi_i)}, E^{(3)}_{(\sP_i,\pi_i)}, E^{(4)}_{(\sP_i,\pi_i)},
        \bigcup_{j=1}^{n-2}D^{(j)}_{(\sP_i,\pi_i)}, \text{ and }E^{(6)}_{(\sP_i,\pi_i)}
        \end{align}
        determine the occurrences of the structures described in conditions (1) to (6) of Theorem~\ref{th: Markov equivalence characterization} in $G^\star$, thus characterizing the Markov equivalence class of $G^\star$.
        For instance for $i\in [2]$, the set $E^{(1)}_{(\sP_i,\pi_i)}$ is equal to the set of $p$-adjacencies of $G^\star$. Throughout this subsection, after we define each of the sets in \eqref{eq: sets in the graphical score}, we demonstrate that set for all $i\in [4]$.
\end{example}
The proof of this theorem can be found at the end of this subsection. We now proceed to carefully define each of the entries of our graphical score.

\begin{definition}
    Let $\pi$ be a partial order on the set $\sP$. For $a,a_1,\ldots,a_t\in \sP$, we write $a\leq_\pi \max\{a_1,\ldots,a_t\}$ if and only if $a\leq_\pi a_i$ for some $i\in [t]$.
\end{definition}

\begin{definition}
    For $a,b\in V^\star$, we write $a\padj b$ to show that $a$ and $b$ are $p$-adjacent in the graph $G^\star$.
\end{definition}

Since $\mathbb{P}$ is both Markov and faithful to $G^\star$, we use $A\independent B \mid C$ and $X_A \independent X_B \mid X_C$ interchangeably for triples $(A,B,C)$ of pairwise disjoint subsets of $V^\star$ from now on.

\begin{proposition} \label{prop: Rich1}
For each $(\sP,\pi)\in \sS$, let
\begin{align*}
    &E^{(1)}_{(\sP,\pi)}\coloneqq \Set{(a,b)\in [n]^2 | a\neq b, \ a \not \independent b \mid \bigcup \Set{C\in \sP |  C\leq_\pi \max\{C_{a,\sP}, C_{b,\sP}\}}\setminus \{a,b\}}.
\end{align*}
Also, define  $S_1 \coloneqq \arg\min_{(\sP,\pi)\in \sS} \left| E^{(1)}_{(\sP,\pi)} \right|$. Then 
\begin{itemize}
    \item the partially ordered partition associated with $G^\star$ is in $S_1$, and
    \item for every $(\sP,\pi)\in S_1$, the set $E^{(1)}_{(\sP,\pi)}$ is equal to the set of $p$-adjacencies in $G^\star$.
\end{itemize}
\end{proposition}
The proof can be found in Appendix~\ref{app:proofs_sec3}.
\begin{example}
For the conditional independence statements and the partially ordered partitions given in Example~\ref{ex: graphical score}, we have
\begin{align*}
    &E^{(1)}_{(\sP_1,\pi_1)} = E^{(1)}_{(\sP_2,\pi_2)} = E^{(1)}_{(\sP_4,\pi_4)} = \{(1,3),(1,4),(2,4),(2,5),(3,4),(3,5),(4,5)\},\\
    &E^{(1)}_{(\sP_3,\pi_3)} =\{(1,3),(1,4),(1,5),(2,3),(2,4),(2,5),(3,4),(3,5),(4,5)\}.
\end{align*}
\end{example}

\begin{proposition} \label{prop: Rich2}
For each $(\sP,\pi)\in \sS$, let
\begin{align*}
    E^{(2)}_{(\sP,\pi)}\coloneqq \big\{(a,b,c)\in [n]^3 \mid &\text{ $a,b,c$ are distinct}, \ (a,b),(c,b)\in E^{(1)}_{(\sP,\pi)},\ (a,c)\not\in E^{(1)}_{(\sP,\pi)},\\
    &C_{b,\sP} \leq_\pi \max\{C_{a,\sP},C_{c,\sP} \}\big\}.
\end{align*}
Also, define  $S_2 \coloneqq \arg\min_{(\sP,\pi)\in S_1} \left| E^{(2)}_{(\sP,\pi)} \right|$. Then 
\begin{itemize}
    \item the partially ordered partition associated with $G^\star$ is in $S_2$, and 
    \item for every $(\sP,\pi)\in S_2$, the set $E^{(2)}_{(\sP,\pi)}$ is equal to the set of unshielded conductors in $G^\star$.
\end{itemize}
\end{proposition}
The proof of this proposition can be found in Appendix~\ref{app:proofs_sec3}.

\begin{example}
For the conditional independence statements and the partially ordered partitions given in Example~\ref{ex: graphical score}, we have
\begin{align*}
    &E^{(2)}_{(\sP_1,\pi_1)}  = E^{(2)}_{(\sP_2,\pi_2)} = \{(1,3,5),(1,4,5),(2,4,3),(2,5,3)\},\\
     &E^{(2)}_{(\sP_3,\pi_3)}  = \{\},\ 
    E^{(2)}_{(\sP_4,\pi_4)} =\{(1,4,2),(1,3,5),(1,4,5),(2,4,3),(2,5,3)\}.
\end{align*}
\end{example}

\begin{proposition} \label{prop: Rich3}
For each $(\sP,\pi)\in \sS$, let
\begin{align*}
    E^{(3)}_{(\sP,\pi)}\coloneqq \big\{(a,b,c)\in [n]^3 \mid &\text{ $a,b,c$ are distinct}, \ (a,b),(c,b)\in E^{(1)}_{(\sP,\pi)},\ (a,c)\not\in E^{(1)}_{(\sP,\pi)},\\ 
    &a \not \independent c \mid \bigcup \Set{C\in \sP |  C\leq_\pi \max\{ C_{a,\sP},C_{b,\sP} , C_{c,\sP}}\} \setminus \{a,c\} \big\}.
\end{align*}
Also, define  $S_3 \coloneqq \arg\min_{(\sP,\pi)\in S_2} \left| E^{(3)}_{(\sP,\pi)} \right|$. Then 
\begin{itemize}
    \item the partially ordered partition associated with $G^\star$ is in $S_3$, and 
    \item for every $(\sP,\pi)\in S_3$, the set $E^{(3)}_{(\sP,\pi)}$ is equal to the set of unshielded perfect non-conductors in $G^\star$.
\end{itemize}
\end{proposition}
The proof of this proposition can be found in Appendix~\ref{app:proofs_sec3}.

\begin{example}
For the conditional independence statements and the partially ordered partitions given in Example~\ref{ex: graphical score}, we have
\begin{align*}
    &E^{(3)}_{(\sP_1,\pi_1)}  = E^{(3)}_{(\sP_2,\pi_2)} = E^{(3)}_{(\sP_4,\pi_4)} =\{\}, \ 
    E^{(3)}_{(\sP_3,\pi_3)}  = \{(1,5,2)\}.
\end{align*}
\end{example}

\begin{proposition} \label{prop: Rich4}
For every $(\sP,\pi)\in \sS$, let
\begin{align*}
    E^{(4)}_{(\sP,\pi)} \coloneqq \big\{\left((a,b_1,c),(a,b_2,c)\right)\in [n]^3\times [n]^3 \mid &\text{ $a,b_1,c$ are distinct}, \ \text{$a,b_2,c$ are distinct}, \\ 
    &(a,b_1),(c,b_1),(a,b_2),(c,b_2) \in E^{(1)}_{(\sP,\pi)}, \ (a,c)\not\in E^{(1)}_{(\sP,\pi)} \\
    &(a,b_1,c),(a,b_2,c) \not \in E^{(2)}_{(\sP,\pi)} , (a,b_1,c),(a,b_2,c) \not \in E^{(3)}_{(\sP,\pi)}\\
    &C_{b_1,\sP} \leq_{\pi} C_{b_2,\sP}
    \big\}.
\end{align*}
Also, define  $S_4 \coloneqq \arg\min_{(\sP,\pi)\in S_3} \left| E^{(4)}_{(\sP,\pi)} \right|$. Then 
\begin{itemize}
    \item the partially ordered partition associated with $G^\star$ is in $S_4$, and 
    \item for every $(\sP,\pi)\in S_4$, $\left((a,b_1,c),(a,b_2,c)\right)\in E^{(4)}_{(\sP,\pi)}$ if and only if $(a,b_1,c),(a,b_2,c)$ are unshielded imperfect non-conductors and $b_1$ is an ancestor of $b_2$ in $G^\star$.
\end{itemize}
\end{proposition}

The proof of this proposition can be found in Appendix~\ref{app:proofs_sec3}.

\begin{example}
For the conditional independence statements and the partially ordered partitions given in Example~\ref{ex: graphical score}, we have
\begin{align*}
    &E^{(4)}_{(\sP_1,\pi_1)}  = E^{(4)}_{(\sP_2,\pi_2)} = \{((1,4,2),(1,4,2))\}, \\
    &E^{(4)}_{(\sP_3,\pi_3)}  = \{((1,3,2),(1,3,2)),((1,3,2),(1,4,2)), ((1,4,2),(1,4,2)),((1,4,2),(1,3,2))\}, \ E^{(4)}_{(\sP_4,\pi_4)} =\{\}.
\end{align*}
\end{example}

\begin{lemma}\label{lemma: mutually exclusive}
Let $(a_0,a_1,\ldots,a_{t+1})\in [n]^{t+2}$ such that for all $i\in \{0,\ldots,t\}$, $a_i \padj a_{i+1}$, and for all $i\in \{2,\ldots,t+1\}$ and for all $j\in \{0,\ldots,i-2\}$, $a_i \not\padj a_j$ . Then one of the following happens:
\begin{enumerate}
    \item Each of the vertices $a_1,\ldots,a_t$ is the ancestor of $a_0$ or $a_{t+1}$ in $G^\star$, or 
    \item there exists $j,\ell\in[t]$ with $j\leq \ell$ such that $(a_{j-1},a_j,a_{j+1})$ and $(a_{\ell-1},a_\ell,a_{\ell+1})$ are mutually exclusive with respect to the uncovered itinerary $(a_{j-1},\ldots,a_{\ell+1})$.
\end{enumerate}
\end{lemma}
The proof of this lemma can be found in Appendix~\ref{app:proofs_sec3}.

\begin{proposition} \label{prop: Rich5}
For all $t\in [n-2]$ and $(\sP,\pi)\in \sS$, let
\begin{align*}
    D^{(t)}_{(\sP,\pi)} \coloneqq \big\{(a_0,a_1,\ldots,a_{t},a_{t+1})\in [n]^{t+2} \mid 
    &\text{ $a_0,a_1,\ldots,a_{t+1}$ are distinct},\\
    &(a_i,a_{i+1}) \in E^{(1)}_{(\sP,\pi)} \text{ for all } i \in \{0,\ldots,t\},\\
    &(a_i,a_j) \not \in E^{(1)}_{(\sP,\pi)} \text{ for all } i \in \{2,\ldots,t+1\}, j \in \{0,\ldots,i-2\},\\
    &C_{a_1,\sP} = C_{a_2,\sP} = \cdots = C_{a_t,\sP},\\
    &C_{a_1,\sP} \not \leq_{\pi} \max \{C_{a_0,\sP}, C_{a_{t+1},\sP} \}
    \big\}.
\end{align*}
Also, define $S^{(1)}_5\coloneqq S_4$ and for all $t\in \{2,\ldots,n-2\}$, $S^{(t)}_5 \coloneqq \arg\max_{(\sP,\pi)\in S^{(t-1)}_5} \left| D^{(t)}_{(\sP,\pi)} \right|$. Then 
\begin{itemize}
    \item the partially ordered partition associated with $G^\star$ is in $S^{(t)}_5$ for all $t\in [n-2]$, and 
    \item for all $t\in [n-2]$ and $(\sP,\pi)\in S^{(t)}_5$, $(a_0,a_1,\cdots,a_{t},a_{t+1})\in D^{(t)}_{(\sP,\pi)}$ if and only if $(a_0,a_1,a_2)$ and $(a_{t-1},a_t,a_{t+1})$ are mutually exclusive with respect to the uncovered itinerary $(a_0,a_1,\ldots,a_{t+1})$ in $G^\star$.
\end{itemize}
\end{proposition}
The proof of this proposition can be found in Appendix~\ref{app:proofs_sec3}.

\begin{example}
For the conditional independence statements and the partially ordered partitions given in Example~\ref{ex: graphical score}, we have
\begin{align*}
    &D^{(1)}_{(\sP_1,\pi_1)}  = D^{(1)}_{(\sP_2,\pi_2)} = \{(1,4,2)\},  &&D^{(2)}_{(\sP_1,\pi_1)}  = D^{(2)}_{(\sP_2,\pi_2)} = \{(1,3,5,2)\}, &&D^{(3)}_{(\sP_1,\pi_1)}  = D^{(3)}_{(\sP_2,\pi_2)} = \{\},\\
    &D^{(1)}_{(\sP_3,\pi_3)} = \{(1,3,2),(1,4,2),(1,5,2)\}, &&D^{(2)}_{(\sP_3,\pi_3)}  = \{\}, &&D^{(3)}_{(\sP_3,\pi_3)}  = \{\},\\
    &D^{(1)}_{(\sP_4,\pi_4)} = \{\},  &&D^{(2)}_{(\sP_4,\pi_4)}  = \{\}, &&D^{(3)}_{(\sP_4,\pi_4)}  = \{\}.
\end{align*}
\end{example}

\begin{proposition} \label{prop: Rich6}
For every $(\sP,\pi)\in \sS$, let
\begin{align*}
    E^{(6)}_{(\sP,\pi)} \coloneqq \bigcup_{t=1}^{n-2} \big\{\left((a_0,a_1,\ldots,a_{t},a_{t+1}),(a_0,b,a_{t+1})\right)\in [n]^{t+2}\times [n]^3 \mid 
    &\text{ $a_0,b,a_{t+1}$ are distinct}, \\ 
    &(a_0,b),(a_{t+1},b) \in E^{(1)}_{(\sP,\pi)}, \\
    &(a_0,b,a_{t+1}) \not \in E^{(2)}_{(\sP,\pi)} , \\
    &(a_0,b,a_{t+1})\not \in E^{(3)}_{(\sP,\pi)},\\
    & (a_0,a_1,\ldots,a_{t},a_{t+1})\in D^{(t)}_{(\sP,\pi)},\\
    &C_{a_1,\sP} \leq_{\pi} C_{b,\sP}
    \big\}.
\end{align*}
Also, define  $S_6 \coloneqq \arg\min_{(\sP,\pi)\in S^{(n-2)}_5} \left| E^{(6)}_{(\sP,\pi)} \right|$. Then 
\begin{itemize}
    \item the partially ordered partition associated with $G^\star$ is in $S_6$, and 
    \item for every $(\sP,\pi)\in S_6$, $\left((a_0,a_1,\cdots,a_{t},a_{t+1}),(a_0,b,a_{t+1})\right)\in E^{(6)}_{(\sP,\pi)}$ if and only if $(a_0,a_1,a_2)$ and $(a_{t-1},a_t,a_{t+1})$ are mutually exclusive with respect to the uncovered itinerary $(a_0,a_1,\ldots,a_{t+1})$, $(a_0,b_1,a_{t+1})$ is an unshielded imperfect non-conductor, and $a_1$ is an ancestor of $b$ in $G^\star$.
\end{itemize}
\end{proposition}
The proof of this proposition can be found in Appendix~\ref{app:proofs_sec3}.

\begin{example}
For the conditional independence statements and the partially ordered partitions given in Example~\ref{ex: graphical score}, we have
\begin{align*}
    &E^{(6)}_{(\sP_1,\pi_1)}  = E^{(6)}_{(\sP_2,\pi_2)} = \{((1,4,2),(1,4,2)),((1,3,5,2),(1,4,2))\}, \\
    &E^{(6)}_{(\sP_3,\pi_3)}  = \{((1,3,2),(1,3,2)),((1,3,2),(1,4,2)), ((1,4,2),(1,4,2)),((1,4,2),(1,3,2))\}, \ E^{(6)}_{(\sP_4,\pi_4)} =\{\}.
\end{align*}
\end{example}

\begin{proof}[Proof of Theorem~\ref{thm:minimal_score_graph}] One can verify that $S_6 = \arg\min_{(\sP,\pi)\in \mathcal{S}} \score(\sP,\pi)$. On the other hand, considering that
\begin{align*}
    S_6 \subseteq S^{(n-2)}_5 \subseteq S^{(n-3)}_5 \subseteq \cdots \subseteq S^{(2)}_5 \subseteq S_4 \subseteq S_3 \subseteq S_2 \subseteq S_1,
\end{align*}
the propositions stated in the current subsection, along with Theorem~\ref{th: Markov equivalence characterization}, imply that for any $(\sP,\pi)\in S_6$, the sets $E^{(t)}_{(\sP,\pi)}, t\in \{1,2,3,4,6\}$,  and $D^{(t)}_{(\sP,\pi)}, t\in \{2,\ldots,n-2\}$ precisely characterize the Markov equivalence class of $G^\star$. Therefore, by finding a minimizer for the graphical score over $\mathcal{S}$, we are able to uniquely determine the Markov equivalence class of $G^\star$.
\end{proof}

\begin{remark} \label{remark: number of tests}
    Note that in order to find the graphical score of a partially ordered partition $(\sP,\pi)$, we only need to perform conditional independence testing to evaluate the first and third coordinates of $\score (\sP,\pi)$ with at most $n^2$ conditional independence tests needed in total.
\end{remark}

\subsection{Greedy optimization of the score} \label{subsection: greedy optimization}
Given that the number of partially ordered partitions on $[n]$ is exponential in $n$, it is not efficient in practice to check every partially ordered partition in order to find the ones with the optimal graphical score. This issue is why we propose a greedy method for optimizing the graphical score over $\mathcal{S}$. 
\begin{definition}
    Let $(\sP,\pi)$ be a partially ordered partition on $[n]$. For $C_1,C_2\in \sP$, we say $(C_1,C_2)$ is a \textit{pair of consecutive sections} if $C_1 \neq C_2$, $C_1\leq_\pi C_2$, and for any $C_3\in \sP\setminus \{C_1,C_2\}$, either $C_1 \not\leq_\pi C_3$ or $C_3\not\leq_\pi C_2$.
\end{definition}
\begin{figure}
    \centering
    \includegraphics[scale=0.4]{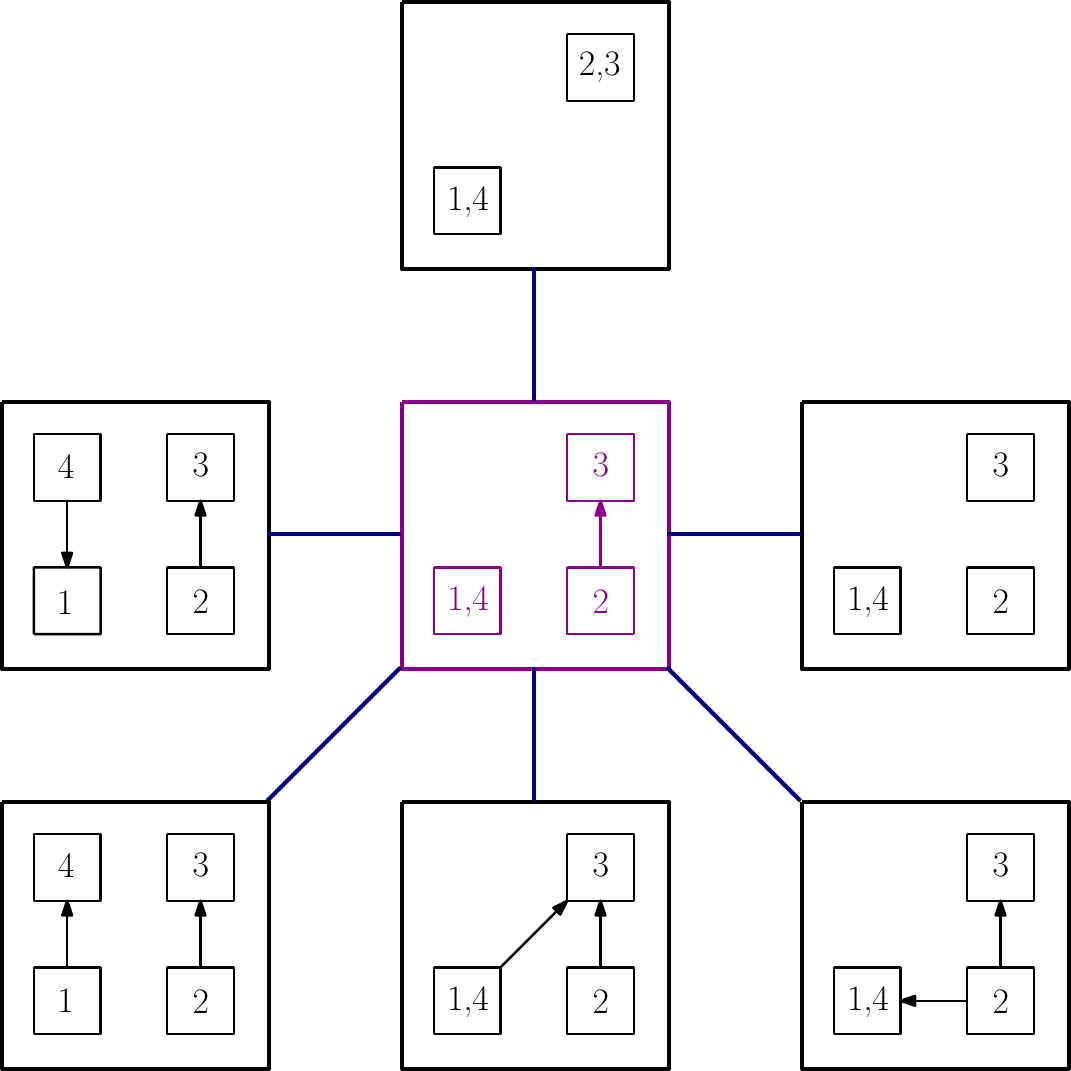}
    \caption{The neighbors of the partition $\left\{\{1,4\},\{2\},\{3\} \right\}$ obeying the partial order $\{2\}\leq\{3\}$ in  graph $\mathcal{G}$. The directions of the edges of $\mathcal{G}$ depend on the graphical scores of the shown partially ordered partitions. A directed arrow from one section of a partition to another indicates that the first section is less than the second one in the partial order.}
    \label{fig:greedy opt}
\end{figure}
In order to optimize the graphical score over $\mathcal{S}$, we perform a depth-first search on the directed graph $\mathcal{G} = (\mathcal{S} , \mathcal{E})$, where $\mathcal{E}\subseteq  \mathcal{S} \times \mathcal{S}$ such that for any two partially ordered partitions $(\sP_1,\pi_1)$ and $(\sP_2,\pi_2$), we have $\left((\sP_1,\pi_1),(\sP_2,\pi_2)\right)\in \mathcal{E}$ if and only if $\score(\sP_1,\pi_1) \geq \score(\sP_2,\pi_2)$ and one of the following happens: 
(See Figure~\ref{fig:greedy opt} for an example.)

\begin{itemize}
\item Partitions $\sP_1$ and $\sP_2$ are the same, but the partial orders differ in exactly one pair, i.e.,
\begin{itemize}
    \item $\sP_1=\sP_2$, $\pi_1\subsetneq \pi_2$ and $|\pi_2\setminus \pi_1|= 1$, or
    \item $\sP_1=\sP_2$, $\pi_2\subsetneq \pi_1$ and $|\pi_1\setminus \pi_2|= 1$.
\end{itemize}
\item $(\sP_2,\pi_2)$ is obtained by moving exactly one element between a pair of consecutive sections in $(\sP_1,\pi_1)$, i.e., 
\begin{itemize}
    \item there exists a consecutive pair $(C_1,C_2)$ in $(\sP_1,\pi_1)$ and $a \in C_1$ such that 
    \begin{align*}
        \sP_2 = \{C \in \sP_1 \mid C\neq C_1 , C\neq C_2\} \cup \{C_1\setminus \{a\}\} \cup \{C_2 \cup \{a\}\}.
    \end{align*}
    Also, $\pi_2$ is obtained by replacing $C_1$ by $C_1\setminus \{a\}$ and $C_2$ by $C_2\cup\{a\}$ in any pair in $\pi_1$ containing $C_1$ or $C_2$, or
    \item there exists a consecutive pair $(C_1,C_2)$ in $(\sP_1,\pi_1)$ and $a \in C_2$ such that 
    \begin{align*}
        \sP_2 = \{C \in \sP_1 \mid C\neq C_1 , C\neq C_2\} \cup \{C_1\cup \{a\}\} \cup \{C_2 \setminus \{a\}\}.
    \end{align*}
    Also, $\pi_2$ is obtained by replacing $C_1$ by $C_1\cup \{a\}$ and $C_2$ by $C_2\setminus\{a\}$ in any pair in $\pi_1$ containing $C_1$ or $C_2$.
\end{itemize}
\item $(\sP_2,\pi_2)$ is obtained by taking out one element from a section of $\sP_1$ and adding that element either right after or right before that section as another independent section, i.e.,
\begin{itemize}
    \item there exists $\hat{C}\in \sP_1$ and $a\in \hat{C}$ such that 
    \begin{align*}
        \sP_2 = \{C \in \sP_1 \mid C\neq \hat{C}\} \cup \{\hat{C}\setminus \{a\}\} \cup \{\{a\}\}.
    \end{align*}
    Also,
    \begin{align*}
        \pi_2 = \pi_1 \cup \{ (C,\{a\}) \mid (C,\hat{C}) \in \pi_1 , \hat{C} \neq C \} \cup \{ (\{a\},C) \mid (\hat{C},C) \in \pi_1 \},
    \end{align*}
    or
    \item there exists $\hat{C}\in \sP_1$ and $a\in \hat{C}$ such that 
    \begin{align*}
        \sP_2 = \{C \in \sP_1 \mid C\neq \hat{C}\} \cup \{\hat{C}\setminus \{a\}\} \cup \{\{a\}\}.
    \end{align*}
    Also,
    \begin{align*}
        \pi_2 = \pi_1 \cup \{ (C,\{a\}) \mid (C,\hat{C}) \in \pi_1\} \cup \{ (\{a\},C) \mid (\hat{C},C) \in \pi_1 , \hat{C} \neq C\}.
    \end{align*}
\end{itemize}
\end{itemize}

Relying on experimental evidence (see Section~\ref{subsection: simulation: greedy Markov equivalence discovery}), we conjecture that by following any of the longest (non-self-intersecting) directed paths in $\mathcal{G}$ starting from an arbitrary partially ordered partition, we eventually reach an optimal partially ordered partition.
\begin{conjecture} \label{conjecture: greedy optimization}
    For any $(\sP,\pi)\in \mathcal{S}$, there exist $(\sP_0,\pi_0)\in S_6$ and a directed path in $\mathcal{G}$ such that the path starts from $(\sP,\pi)$ and ends at $(\sP_0,\pi_0)$.
\end{conjecture}
Based on Conjecture~\ref{conjecture: greedy optimization}, we propose Algorithm~\ref{alg: Markov equivalence class discovery} for finding an optimal partially ordered partition.
\begin{algorithm}
\caption{Markov equivalence class discovery}\label{alg: Markov equivalence class discovery}
\begin{algorithmic}[1]
\Require{The set of all the conditional independence statements satisfied by $\mathbb{P}$ and an initial $(\sP_1,\pi_1)\in \mathcal{S}$.}
\Ensure{A partially ordered partition in $S_6$.}
\State Set $\hat{\sP} \coloneqq \sP_1$ and $\hat{\pi} \coloneqq \pi_1$.
\State Perform a depth-first search on $\mathcal{G}$ with root $(\hat{\sP},\hat{\pi})$ to find a directed path from $(\hat{\sP},\hat{\pi})$ to a partially ordered partition $(\tilde{\sP},\tilde{\pi})$ with $\score(\hat{\sP},\hat{\pi}) > \score(\tilde{\sP},\tilde{\pi})$. \label{alg: Markov equivalence class discovery repeat}
\If {$(\tilde{\sP},\tilde{\pi})$ is found} 
\State Set $\hat{\sP} \coloneqq \tilde{\sP}$ and $\hat{\pi} \coloneqq \tilde{\pi}$, and go back to step~\ref{alg: Markov equivalence class discovery repeat}.
\Else  
\State \Return $(\hat{\sP},\hat{\pi})$.
\EndIf
\end{algorithmic}
\end{algorithm}
In practice, Algorithm~\ref{alg: Markov equivalence class discovery} is still slow for graphs with a large number of vertices. The reason is that once the algorithm finds a partially ordered partition $(\hat{\sP},\hat{\pi})$ in $S_6$, it needs to perform a full depth-first search on $\mathcal{G}$ with root $(\hat{\sP},\hat{\pi})$ before it makes sure that there is no partially ordered partition with a lower score. So, it could end up traversing all the elements of $S_6$ before outputting $(\hat{\sP},\hat{\pi})$. Given that $|S_6|$ can be large, this makes the algorithm impractical. We propose making Algorithm~\ref{alg: Markov equivalence class discovery} even "greedier" by stopping the depth-first search performed in step~\ref{alg: Markov equivalence class discovery repeat} as soon as a directed path of length $N$ consisting of partially ordered partitions with the same score is observed, where $N$ is a threshold given to the algorithm as an input. One can choose to repeat the algorithm $M$ times with $M$ different initial partially ordered partitions to get better results, where $M$ is also part of the input. This results in a modified version of Algorithm~\ref{alg: Markov equivalence class discovery}, which is presented as Algorithm~\ref{alg: greedy Markov equivalence class discovery}. Subsection~\ref{subsection: simulation: greedy Markov equivalence discovery} shows the running time and success rate of testing Algorithm~\ref{alg: greedy Markov equivalence class discovery} on random graphs. 

\begin{algorithm}
\caption{Greedy Markov equivalence class discovery}\label{alg: greedy Markov equivalence class discovery}
\begin{algorithmic}[1]
\Require{The set of all the conditional independence statements satisfied by $\mathbb{P}$, two positive integers $N$ and $M$, and initial partially ordered partitions $(\sP_1,\pi_1),\ldots,(\sP_M,\pi_M)$.}
\Ensure{A partially ordered partition in $S_6$.}
\State Set $A\coloneqq \emptyset$. 
\For {$i$ in $1:M$}
\State Set $\hat{\sP} \coloneqq \sP_i$ and $\hat{\pi} \coloneqq \pi_i$.
\State Perform a depth-first search on $\mathcal{G}$ with root $(\hat{\sP},\hat{\pi})$ to find a directed path from $(\hat{\sP},\hat{\pi})$ to\label{alg: greedy Markov equivalence class discovery repeat}
\NoNumber{a partially ordered partition $(\tilde{\sP},\tilde{\pi})$ with $\score(\hat{\sP},\hat{\pi}) > \score(\tilde{\sP},\tilde{\pi})$. Stop the depth-first search}
\NoNumber{once a directed path of length $N$ consisting of partially ordered partitions of score $\score(\hat{\sP},\hat{\pi})$}
\NoNumber{is generated. }
\If {$(\tilde{\sP},\tilde{\pi})$ is found} 
\State Set $\hat{\sP} \coloneqq \tilde{\sP}$ and $\hat{\pi} \coloneqq \tilde{\pi}$, and go back to step~\ref{alg: greedy Markov equivalence class discovery repeat}.
\Else  
\State Set $A \coloneqq A \cup \{(\hat{\sP},\hat{\pi})\}$.
\EndIf
\EndFor
\State \Return $\arg\min_{(\sP,\pi)\in A} \score(\sP,\pi)$.
\end{algorithmic}
\end{algorithm}

\begin{remark}[Running time]
    For any partially ordered partition $(\sP,\pi) \in \mathcal{S}$, the time complexity of computing the sets \eqref{eq: sets in the graphical score} is respectively $O(n^2)$, $O(n^3)$, $O(n^3)$, $O(n^4)$, $O(n^2 2^n)$, and $O(n\  2^n)$. Therefore, the time complexity of computing $\score(\sP,\pi)$ is $O(n^2 2^n)$. Note that in the process of optimizing the graphical score in Algorithms~\ref{alg: Markov equivalence class discovery} and~\ref{alg: greedy Markov equivalence class discovery}, not the whole vector of the graphical score needs to be necessarily computed for each partially ordered partition since the graphical score vectors are being compared with respect to the lexicographical order.
\end{remark}

\section{Discovering a Markov equivalent directed graph} \label{graph discovery}
After  determining the Markov equivalence class of $G^\star$, a natural question is whether one can construct a Markov equivalent graph (or a few of them) to $G^\star$. In this section, we propose two methods for doing so. 

In this second step, our algorithm takes in the partially ordered partition produced in the first step and produces a directed graph. This graph is consistent with the $d$-separation statements given on input and induces the same partially ordered partition as the one produced in step 1. We note that some partially ordered partitions do not correspond to a directed graph which satisfies the required $d$-separation statements. If this happens, we return to step 1 and produce a different optimal partially ordered partition. We then repeat step 2. 
\begin{theorem}\label{thm: markov equivalent graph}
Suppose that $(\sP,\pi)\in S_6$, and $G=([n],E)$ is a directed graph with the following properties:
\begin{enumerate}
    \item For all $(a,b)\in [n]^2$, $a$ and $b$ are $p$-adjacent in $G$ if and only if $(a,b)\in E^{(1)}_{(\sP,\pi)}$. \label{thm: property1}
    \item For all $(a,b)\in [n]^2$, $a$ is an ancestor of $b$ in $G$ if and only if $C_{a,\sP}\leq_\pi C_{b,\sP}$.\label{thm: property2}
    \item If $(a,b,c)\in E^{(3)}_{(\sP,\pi)}$ and for all $b^\prime\in [n]$ with $C_{b^\prime,\sP} \leq_\pi C_{b,\sP}$ and $C_{b^\prime,\sP} \neq C_{b,\sP}$, $(a,b^\prime,c)\not\in E^{(3)}_{(\sP,\pi)}$, then $a$ and $c$ have a common child in $C_{b,\sP}$ in graph $G$.\label{thm: property3}
    \item If $(a,b),(c,b)\in E^{(1)}_{(\sP,\pi)}$, $(a,c)\not\in E^{(1)}_{(\sP,\pi)}$, $C_{b,\sP}\not\leq_\pi \max\{C_{a,\sP},C_{c,\sP}\}$, and $(a,b,c)\not\in E^{(3)}_{(\sP,\pi)}$, then $a$ and $c$ don't have a common child in $C_{b,\sP}$ in graph $G$.\label{thm: property4}
\end{enumerate}
Then $G$ is Markov equivalent to $G^\star$.

Moreover, $G^\star$ satisfies properties~\ref{thm: property1} to~\ref{thm: property4} with respect to the partially ordered partition associated with $G^\star$.
\end{theorem}
\begin{proof}
We prove that conditions~\ref{thm: markov con1} to~\ref{thm: markov con6} of Theorem~\ref{th: Markov equivalence characterization} are satisfied for $G$ and $G^\star$:
\begin{itemize}
    \item Proposition~\ref{prop: Rich1} and property~\ref{thm: property1} guarantee condition~\ref{thm: markov con1} of Theorem~\ref{th: Markov equivalence characterization}.
    \item Proposition~\ref{prop: Rich2} and properties~\ref{thm: property1} and~\ref{thm: property2} guarantee condition~\ref{thm: markov con2} of Theorem~\ref{th: Markov equivalence characterization}.
    \item Let $(a,b,c)$ be an unshielded perfect non-conductor in $G^\star$. Then by Proposition~\ref{prop: Rich3}, $(a,b,c)\in E^{(3)}_{(\sP,\pi)}$. Let $C\in \sP$ be a minimal section of partition $\sP$ such that $C\leq_\pi C_{b,\sP}$ and there exists $b^\prime \in C$ with $(a,b^\prime,c)\in E^{(3)}_{(\sP,\pi)}$. Then by property~\ref{thm: property3}, $a$ and $c$ have a common child $b_0$ in $C$ in $G$. Since $C\leq_\pi C_{b,\sP}$, property~\ref{thm: property2} guarantees that $b$ is the descendant of $b_0$ in $G$. Moreover, since $(a,b),(c,b)\in E^{(1)}_{(\sP,\pi)}$ and $(a,c)\not\in E^{(1)}_{(\sP,\pi)}$, property~\ref{thm: property1} implies that $a,b$ and $b,c$ are $p$-adjacent in $G$, and $a,c$ are not $p$-adjacent in $G$. 
    Hence, $(a,b,c)$ is an unshielded perfect non-conductor in $G$.

    Now let $(a,b,c)$ be an unshielded perfect non-conductor in $G$. Then there exists $b^\prime \in [n]$ such that $b$ is a descendant of $b^\prime$ and $b^\prime$ is the common child of $a$ and $c$ in $G$. Since $a$ and $c$ are not $p$-adjacent in $G$, $b^\prime$ cannot be an ancestor of $a$ or $c$ in $G$. So, $(a,b^\prime,c)\in D^{(1)}_{(\sP,\pi)}$ by properties~\ref{thm: property1} and~\ref{thm: property2}. On the other hand, if $(a,b,c)$ is not an unshielded perfect non-conductor in $G^\star$, then by Proposition~\ref{prop: Rich3}, $(a,b,c)\not\in E^{(3)}_{(\sP,\pi)}$.  So, by properties~\ref{thm: property1} and~\ref{thm: property2}, $\left((a,b^\prime,c),(a,b,c)\right)\in E^{(6)}_{(\sP,\pi)}$, which means that $b^\prime\in \an_{G^\star}(b)$ by Proposition~\ref{prop: Rich6}. Therefore, considering that $b$ is not a descendant of a common child of $a$ and $c$ in $G^\star$, $b^\prime$ cannot be a descendant of a common child of $a$ and $c$ in $G^\star$ either. However, $b^\prime$ being a common child of $a$ and $c$ in $G$, along with property~\ref{thm: property4}, implies that $(a,b^\prime,c)\in E^{(3)}_{(\sP,\pi)}$, which is a contradiction by Proposition~\ref{prop: Rich3}. This contradiction proves that $(a,b,c)$ is unshielded perfect non-conductor in $G^\star$ as well.
    \item Proposition~\ref{prop: Rich4}, along with the fact that $G$ and $G^\star$ have the same set of unshielded imperfect non-conductors and property~\ref{thm: property2}, guarantees condition~\ref{thm: markov con4} of Theorem~\ref{th: Markov equivalence characterization}.
    \item Proposition~\ref{prop: Rich5} and properties~\ref{thm: property1} and~\ref{thm: property2} guarantee condition~\ref{thm: markov con5} of Theorem~\ref{th: Markov equivalence characterization}.
    \item Proposition~\ref{prop: Rich6}, along with the fact that $G$ and $G^\star$ have the same set of unshielded imperfect non-conductors and property~\ref{thm: property2}, guarantees condition~\ref{thm: markov con6} of Theorem~\ref{th: Markov equivalence characterization}. 
\end{itemize}
Hence, $G$ and $G^\star$ are Markov equivalent.

Let $(\sP_0,\pi_0)$ be the partially ordered partition associated with $G^\star$. Then property~\ref{thm: property1} holds because of Proposition~\ref{prop: Rich1}, and property~\ref{thm: property2} holds due to the definition of $\pi_0$. 

Assume $(a,b,c)\in E^{(3)}_{(\sP_0,\pi_0)}$ and for all $b^\prime\in [n]$ with $C_{b^\prime,\sP_0} \leq_{\pi_0} C_{b,\sP_0}$ and $C_{b^\prime,\sP_0} \neq C_{b,\sP_0}$, $(a,b^\prime,c)\not\in E^{(3)}_{(\sP_0,\pi_0)}$. Then by Proposition~\ref{prop: Rich3}, in graph $G^\star$, $a$ and $c$ have a common child $b^\prime$ such that $b^\prime$ is an ancestor of $b$. By the same proposition, $(a,b^\prime,c)\in E^{(3)}_{(\sP_0,\pi_0)}$. So, given that $C_{b^\prime,\sP_0}\leq_{\pi_0} C_{b,\sP_0}$, we can conclude $C_{b^\prime,\sP_0}=C_{b,\sP_0}$, which proves $a$ and $c$ have a common child in $C_{b,\sP_0}$ in graph $G^\star$, and hence, property~\ref{thm: property3}.

If $(a,b),(c,b)\in E^{(1)}_{(\sP_0,\pi_0)}$, $(a,c)\not\in E^{(1)}_{(\sP_0,\pi_0)}$, $C_{b,\sP_0}\not\leq_{\pi_0} \max\{C_{a,\sP_0},C_{c,\sP_0}\}$, and $(a,b,c)\not\in E^{(3)}_{(\sP_0,\pi_0)}$, then by Propositions~\ref{prop: Rich1} and~\ref{prop: Rich3}, $(a,b,c)$ is an unshielded imperfect non-conductor in $G^\star$, which means $a$ and $c$ don't have a common child in $C_{b,\sP}$. Therefore, property~\ref{thm: property4} is also satisfied by $G^\star$.
\end{proof}
According to Theorem~\ref{thm: markov equivalent graph}, in order to construct a Markov equivalent graph to $G^\star$, we need to be able to recover the edges within a strongly connected component of $G$. This needs to be done so that all properties from Theorem~\ref{thm: markov equivalent graph} are preserved (for example, we have prescribed $p$-adjacencies that need to be the precise $p$-adjacencies in the graph we construct).

\begin{example}
    Consider the conditional independence statements and the partially ordered partitions presented in Example~\ref{ex: graphical score}. As discussed in that example, we have $(\sP_1, \pi_1), (\sP_2, \pi_2) \in S_6$. Corresponding to $(\sP_1, \pi_1)$, the graph $G^\star$ satisfies the properties listed in Theorem~\ref{thm: markov equivalent graph}. However, if there were a graph $G$ with these properties corresponding to $(\sP_2, \pi_2)$, then by property~\ref{thm: property2}, the set $\{1,2\}$ would form a strongly connected component in $G$, which makes vertices 1 and 2 adjacent, and thus, $p$-adjacent in $G$. But this violates property~\ref{thm: property1}. Therefore, although $(\sP_2,\pi_2)$ characterizes the Markov equivalence class of $G^\star$ by Theorem~\ref{thm:minimal_score_graph}, it does not admit a graph in this Markov equivalence class in the sense of Theorem~\ref{thm: markov equivalent graph}. 
\end{example}

\begin{figure}
        \centering
        \includegraphics[scale=0.5]{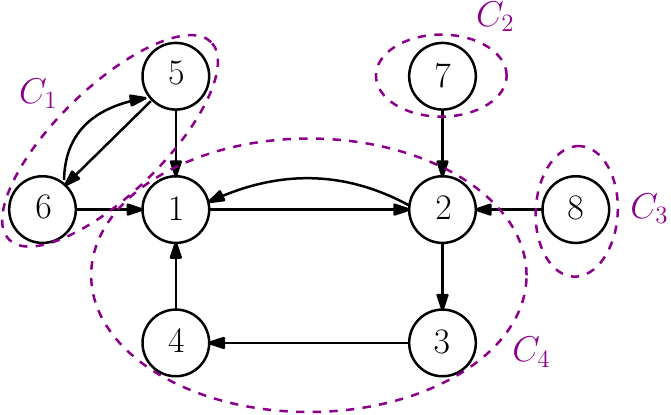}
        \caption{A directed graph with four strongly connected components $C_1,C_2,C_3,C_4$.}
        \label{fig:example for SCCR}
    \end{figure}

\begin{definition} \label{def: strongly connected component recovery algorithm}
    We call an algorithm a \textit{strongly connected component recovery algorithm}, or for short, an \textit{SCCR algorithm}, if for any given subset $C\subseteq[n]$ and sets $A_C\subseteq C^2$, $B_C \subseteq ([n]\setminus C) \times C$, $ComCh_C\subseteq ([n]\setminus C)^2$ and $NoComCh_C\subseteq ([n]\setminus C)^2$, the algorithm determines whether a directed graph $G_C = ([n],E_C)$ with the following properties exists and outputs $E_C$ if so:
\begin{itemize}
    \item The $p$-adjacent pairs in $G_C$ are exactly the pairs in $A_C\cup B_C$,
    \item $C$ is a strongly connected component in $G_C$, and $C \not \leq_{G_C} C_{a,G_C}$ for all $a\in [n]\setminus C$, i.e., all the edges connecting vertices of $C$ to the outside vertices are incoming to $C$,
    \item for any $(a,b)\in ComCh_C$, $a$ and $b$ have a common child in $C$, and 
    \item for any $(a,b)\in NoComCh_C$, $a$ and $b$ do not have a common child in $C$.
\end{itemize}
Note that in this definition without loss of generality, $A_C$ can be assumed to be symmetric (i.e., for any $(a,b)\in A_C$, $(b,a)\in A_C$) and not contain any pair of the form $(a,a)$. From now on, these two properties are assumed for $A_C$.

\begin{definition}
    For any set $D\subseteq [n]^2$, we define $\sym(D)\coloneqq \bigcup_{(a,b)\in D}\{(a,b),(b,a)\}$.
\end{definition}
\end{definition}
\begin{example}\label{example: SCCR}
    Suppose $G^\star$ is the graph in Figure~\ref{fig:example for SCCR}. Define 
    \begin{align*}
        &C_1\coloneqq \{5,6\}, \ A_{C_1}\coloneqq \sym\{(5,6) \}, \ B_{C_1} \coloneqq \{\}, \ ComCh_{C_1} \coloneqq \{\}, \ NoComCh_{C_1} \coloneqq \{ \},\\
        &C_2\coloneqq \{7\}, \ A_{C_2}\coloneqq \{\}, \ B_{C_2} \coloneqq \{\}, \ ComCh_{C_2} \coloneqq \{\}, \ NoComCh_{C_2} \coloneqq \{ \},\\
        &C_3\coloneqq \{8\}, \ A_{C_3}\coloneqq \{\}, \ B_{C_3} \coloneqq \{\}, \ ComCh_{C_3} \coloneqq \{\}, \ NoComCh_{C_3} \coloneqq \{ \},\\
        &C_4\coloneqq \{1,2,3,4\}, \ A_{C_4}\coloneqq \sym\{(1,2),(2,3),(3,4),(1,4),(2,4)\}, \\
        &B_{C_4} \coloneqq \{(5,1),(5,2),(5,4),(6,1),(6,2),(6,4),(7,1),(7,2),(8,1),(8,2)\}, \\
        &ComCh_{C_4} \coloneqq \{(7,8)\}, \ NoComCh_{C_4} \coloneqq \{ (5,7),(5,8),(6,7),(6,8)\}.
    \end{align*}
    If for each $i\in [4]$, $E_{C_i}$ is the output of an SCCR algorithm given the input $C_i$, $A_{C_i}$, $B_{C_i}$, $ComCh_{C_i}$ 
 and $NoComCh_{C_i}$, then by Theorem~\ref{thm: markov equivalent graph}, the graph $G=\left([8],\cup_{i=1}^4 E_{C_i}\right)$ is Markov equivalent to $G^\star$.
\end{example}
For a given $(\sP,\pi)\in S_6$, suppose $\sP = \{C_1,\ldots,C_{|\sP|}$\} is a linear extension of the partial order $\pi$, i.e., for all $i,j\in [|\sP|]$, we have $i\leq j$ if $C_i\leq_\pi C_j$. Then for each $i\in [|\sP|]$, appropriate $A_{C_i}, B_{C_i},ComCh_{C_i}$, and $NoComCh_{C_i}$ can be determined using conditions~\ref{thm: property1},~\ref{thm: property3} and~\ref{thm: property4} of Theorem~\ref{thm: markov equivalent graph}, and applying an SCCR algorithm to $C_i$ yields $E_{C_i}$, as described in Definition~\ref{def: strongly connected component recovery algorithm}, if such $E_{C_i}$ exists. Having chosen the inputs of the algorithm wisely, the graph $G=\left([n],\bigcup_{i=1}^{|\sP|} E_{C_i}\right)$ is Markov equivalent to $G^\star$ by Theorem~\ref{thm: markov equivalent graph}. See Algorithm~\ref{alg: graph discovery} for more details.

\begin{algorithm}[t]
\caption{Markov equivalent graph discovery}\label{alg: graph discovery}
\begin{algorithmic}[1]
\Require{The set of all the conditional independence statements satisfied by $\mathbb{P}$, a partially ordered partition $(\sP,\pi)$ in $S_6$ obtained by Algorithm~\ref{alg: Markov equivalence class discovery} or~\ref{alg: greedy Markov equivalence class discovery}, and an SCCR algorithm $\Alg()$.}
\Ensure{A graph Markov equivalent to $G^\star$.}
\State Set $\hat{\sP}\coloneqq \sP$ and $\hat{\pi} \coloneqq \pi$.
\If {there exist $a,b\in [n]$ such that $(a,b)\in E^{(1)}_{(\hat{\sP},\hat{\pi})}$, but $C_{a,\hat{\sP}} \not\leq_{\hat{\pi}} C_{b,\hat{\sP}}$ and $C_{b,\hat{\sP}} \not\leq_{\hat{\pi}} C_{a,\hat{\sP}}$}  \label{alg: graph discovery repeat2}
\State Go to step~\ref{alg: graph discovery repeat}.
\EndIf
\If {for a consecutive pair $(C_1,C_2)$ in $(\hat{\sP},\hat{\pi})$ and for all $a\in C_1,b\in C_2$, $(a,b)\not\in E^{(1)}_{(\hat{\sP},\hat{\pi})}$}
\State Go to step~\ref{alg: graph discovery repeat}.
\EndIf
\State Consider a linear extension $\hat{\sP} = \{C_1,\ldots,C_{|\hat{\sP}|}\}$ of the partial order $\hat{\pi}$.
\For{$i$ in $1:|\hat{\sP}|$} 
\State Set $A_{C_i} \coloneqq \Set{(a,b) \in E^{(1)}_{(\hat{\sP},\hat{\pi})} | a,b\in C_i}$.
\vspace{1mm}
\State Set $B_{C_i} \coloneqq \Set{(a,b) \in E^{(1)}_{(\hat{\sP},\hat{\pi})} | a\in \bigcup_{j=1}^{i-1}C_j,b\in C_i}$.
\vspace{1mm}
\State Set $ComCh_{C_i}$ and $NoComCh_{C_i}$ to be the sets of pairs that must have and must not have 
\NoNumber{common children in $C_i$ per conditions~\ref{thm: property3} and~\ref{thm: property4} of Theorem~\ref{thm: markov equivalent graph} respectively.}
\vspace{1mm}
\State Set $E_{C_i} \coloneqq \Alg\left(C_i,A_{C_i},B_{C_i},ComCh_{C_i},NoComCh_{C_i}\right)$.
\vspace{1mm}
\If {no $E_{C_i}$ is outputted}
\State Go to step~\ref{alg: graph discovery repeat}.
\EndIf
\EndFor
\State \Return $G=\left([n],\bigcup_{i=1}^{|\hat{\sP}|} E_{C_i}\right)$.
\State Perform a depth-first search on $\mathcal{G}$ with root $(\hat{\sP},\hat{\pi})$ to find a different partially ordered partition $(\tilde{\sP},\tilde{\pi})$ with $\score(\hat{\sP},\hat{\pi})=\score(\tilde{\sP},\tilde{\pi})$. Set $\hat{\sP}\coloneqq \tilde{\sP}$ and $\hat{\pi} \coloneqq \tilde{\pi}$, and go back to step~\ref{alg: graph discovery repeat2}.\label{alg: graph discovery repeat}
\end{algorithmic}
\end{algorithm}
Note that in step~\ref{alg: graph discovery repeat} of Algorithm~\ref{alg: graph discovery}, we might not need to actually perform a depth-first search. When running Algorithm~\ref{alg: Markov equivalence class discovery} or~\ref{alg: greedy Markov equivalence class discovery} for finding the initial $(\sP_1,\pi_1)\in S_6$, a list of partially ordered partitions with an equal score to $\score(\sP_1,\pi_1)$ is observed before outputting $(\sP_1,\pi_1)$. So, in step~\ref{alg: graph discovery repeat}, $(\tilde{\sP},\tilde{\pi})$ can be chosen from this list as long as the list is not exhausted.

In the remainder of this section, we discuss two SCCR algorithms. Subsection~\ref{subsection: SCCR2} proposes an efficient algorithm, which is conjectured to be correct based on experimental evidence. Subsection~\ref{subsection: SCCR1} proves the correctness of another SCCR algorithm which relies on \textit{submodular flow polyhedrons}. 
\subsection{SCCR algorithm 1: Construct and correct} \label{subsection: SCCR2}
In this subsection, we introduce an SCCR algorithm. Suppose a subset $C\subseteq[n]$ and sets $A_C\subseteq C^2$, $B_C \subseteq ([n]\setminus C) \times C$, $ComCh_C\subseteq ([n]\setminus C)^2$ and $NoComCh_C\subseteq ([n]\setminus C)^2$ are given. Our goal is to construct a set $E_C \subseteq [n]^2$ such that as described in Definition~\ref{def: strongly connected component recovery algorithm}, $E_C$ preserves $C$ as a strongly connected component, preserves the $p$-adjacencies in $A_C\cup B_C$, assigns a common child in $C$ to each pair in $ComCh_C$, and prevents the pairs in $NoComCh_C$ from having a common child in $C$. In order to do so, in each iteration Algorithm~\ref{alg: SCCR2} adds an edge from $A_C\cup B_C$ to $E_C$ and if the edge causes a pair not in $A_C\cup B_C$ to become $p$-adjacent (conditioned on $C$ becoming a strongly connected component in $([n],E_C)$ eventually) or causes a pair in $NoComCh_C$ to have a common child in $C$, the algorithm keeps correcting its previous choices of edges in $E_C$ until the problem is settled. The algorithm continues adding edges to $E_C$ until $A_C\cup B_C$ is exhausted.
\begin{definition}
    Suppose the sets $C\subseteq[n]$, $A_C\subseteq C^2$, $B_C \subseteq ([n]\setminus C) \times C$, and $NoComCh_C\subseteq ([n]\setminus C)^2$ are given. Let $\mathcal{F}$ denote the family of all multisubsets of $[n]$. An edge $(a,b)\in A_C \cup B_C$ is said to be \textit{safe} with respect to a function $Parents: C \to \mathcal{F}$ if
    \begin{itemize}
        \item for every $v\in Parents(b)\setminus \{a\}$, $(v,a)\not\in NoComCh_C$, and
        \item for every $v\in Parents(b) \setminus \{a\}$, if $v\in C$ or $a\in C$, then $(v,a)\in A_C\cup B_C$ or $(a,v)\in A_C \cup B_C$.
    \end{itemize}
\end{definition}

\begin{definition}
    Suppose the sets $C\subseteq[n]$, $A_C\subseteq C^2$, $B_C \subseteq ([n]\setminus C) \times C$, and $NoComCh_C\subseteq ([n]\setminus C)^2$ are given. Two edges $(a,b),(c,b)\in [n]\times C$ are said to be \textit{incompatible} if $a\neq c$ and one of the following happens:
    \begin{itemize}
        \item $a,c\in [n]\setminus C$ and $(a,c)\in NoComCh_C$, or
        \item $C \cap \{a,c\} \neq \emptyset$, $(a,c)\not\in A_C\cup B_C$, and $(c,a)\not\in A_C \cup B_C$.
    \end{itemize}
\end{definition}

Below, we describe the idea behind Algorithm~\ref{alg: SCCR2} briefly:

Algorithm~\ref{alg: SCCR2} maintains three (ordered) lists $D,D_{copy}$ and $Cause$ throughout its execution. $D$ records the edges that are added in each iteration. If there is need for correcting the previously added edges after adding a specific edge, $D$ will change accordingly. However, an edge is allowed to be corrected at most once after being added to $D$ for the first time. In order to record the status of each edge in $D$, $D_{copy}$ assigns a code to the edge in the following manner. For all $i\in \left[\left|D\right|\right]$,
\begin{itemize}
    \item $D_{copy}[i]=-1$ if the $i$th edge added to $D$ has been flipped once resulting in $D[i]$.
    \item $D_{copy}[i]=0$ if the $i$th edge added to $D$ has been removed.
    \item $D_{copy}[i]=1$ if the $i$th edge added to $D$ was added along with another edge to preserve the $p$-adjacency created by an edge that is removed from $D$.
    \item $D_{copy}[i]=2$ if the $i$th edge added to $D$ was added along with another edge to guarantee a pair in $ComCh_C$ has a common child in $C$.
    \item $D_{copy}[i]=D[i]$ otherwise.
\end{itemize}
Finally, the list $Cause$ records the reason for each edge having a certain code in $D_{copy}$. Specifically, for all $i\in \left[\left|D\right|\right]$, 
\begin{itemize}
    \item $Cause[i]$ is the index in $D$ of the edge the adding of which has caused corrections affecting $D[i]$ if $D_{copy}[i]\in \{-1,0,1\}$.
    \item $Cause[i]$ is $v\in [n]\setminus C$ such that $D[i]$ was added to guarantee the pair $(D[i]_1,v)$ has a common child in $C$ if $D_{copy}[i]=2$
    \item $Cause[i]=i$ otherwise.
\end{itemize}

In addition to $C,A_C,B_C,ComCh_C$ and $NoComCh_C$, the algorithm receives a positive integer $N$ as its input, which determines the maximum number of times allowed for the algorithm to erase part of its progress and rebuild $D$ starting from a certain point and avoiding certain choices. As long as the algorithm hasn't reached its maximum number of attempts, it begins by creating the list $Checked$ and the function $Parents: C\to \mathcal{F}$ based on $D,D_{copy}$ and $Cause$, which in the initial attempt, are empty lists. $Checked$ contains all $p$-adjacencies in $A_C\cup B_C$ that were added to $D$ as an edge at some point, and $Parents$ corresponds to each $v\in C$ a multiset containing all the parents of $v$ determined by the edges in $D$ that haven't been assigned code 0 in $D_{copy}$. Afterwards, the algorithm assigns a common child to any pair in $ComCh_C$ that has not been assigned one already.

Until $Checked$ contains all $p$-adjacencies in $A_C\cup B_C$, in each iteration, the algorithm adds an edge from $A_C\cup B_C$ to $D$ such that the edges added to $D$ form an \textit{almost} directed path. Lists $D$, $D_{copy}$, $Cause$, $Parents$ and $Checked$ are updated accordingly and if the added edge is not safe with respect to $Parents$, a correction process starts. In this process, the algorithm first tries to remove the unsafe edge. However, any edge ensures two vertices being $p$-adjacent. So, we are only allowed to remove an edge if instead two safe edges can be added to $D$ preserving that $p$-adjacency. Moreover, if the edge we are hoping to remove is a guarantee for a pair in $ComCh_C$ having a common child in $C$, a different common child in $C$ should be secured for this pair before removing the edge. 

If the removal is not successful due to any of the above-mentioned reasons, the algorithm flips the edge if both of its vertices lie in $C$ and otherwise, leaves it as it is for now. Call the flipped edge in the former case and the original edge in the latter case $e$. Unless $e$ is safe, the correction process continues by listing all the edges in $D$ incompatible with $e$. If any edge $\tilde{e}$ in this list has already been affected in a correction process, i.e., is assigned any of the codes $-1$ or $1$ in $D_{copy}$, the algorithm jumps back to the stage right before adding the cause assigned to $\tilde{e}$ in $Cause$, and avoids adding $\tilde{e}$ to $D$ in the first iteration right after that.  If no incompatible edge with $e$ has code $-1$ or $1$ in $D_{copy}$, the algorithm adds $e$ to a list containing all the potentially unsafe edges called $Potential.Problems$, and then chooses one of the incompatible edges and continues the correction process with that edge. 

After the correction process ends, all the edges in $Potential.Problems$ are checked to be safe. If for any edge $e\in Potential.Problems$, that's not the case, the algorithm chooses one of the edges incompatible with $e$ randomly and jumps back to the stage right before adding the cause listed for this  incompatible edge in $Cause$. Please see Algorithm~\ref{alg: SCCR2} in Appendix~\ref{app: SCCR2-alg} for more details.

The structure of Algorithm~\ref{alg: SCCR2} ensures the following:
\begin{proposition}
Suppose a subset $C\subseteq[n]$ and sets $A_C\subseteq C^2$, $B_C \subseteq ([n]\setminus C) \times C$, $ComCh_C\subseteq ([n]\setminus C)^2$ and $NoComCh_C\subseteq ([n]\setminus C)^2$ are given. For any $N\in \mathbb{N}$, if Algorithm~\ref{alg: SCCR2} outputs a set $E_C\subseteq [n]^2$, then $E_C$ satisfies all the properties listed in Definition~\ref{def: strongly connected component recovery algorithm}.
\end{proposition}
Based on experimental evidence presented in Subsection~\ref{subsection: simulation SCCR}, we conjecture the following:
\begin{conjecture}\label{conj: SCCR2}
    For all $n\in \mathbb{N}$, there exists $N\in \mathbb{N}$ that satisfies the following:
    Let $(\sP_0,\pi_0)$ be the partially ordered partition associated with $G^\star$. Suppose $C\in \sP_0$ and
    \begin{align*}
        &A_{C} \coloneqq \Set{(a,b) \in E^{(1)}_{(\sP_0,\pi_0)} | a,b\in C},
         B_{C} \coloneqq \Set{(a,b) \in E^{(1)}_{(\sP_0,\pi_0)} | a\not \in C,b\in C, C_{a,\sP_0}\leq_{\pi_0} C_{b,\sP_0}}.
    \end{align*}
    Moreover, let $ComCh_{C}$ and $NoComCh_{C}$ be the sets of pairs that must have and must not have common children in $C$ per conditions~\ref{thm: property3} and~\ref{thm: property4} of Theorem~\ref{thm: markov equivalent graph} respectively. Then given $C, A_C,B_C,ComCh_C,NoComCh_C$ and $N$, with high probability Algorithm~\ref{alg: SCCR2} outputs a set $E_C\subseteq [n]^2$ satisfying the properties described in Definition~\ref{def: strongly connected component recovery algorithm}.
\end{conjecture}

\begin{remark}[Running time]
    The structure of Algorithm~\ref{alg: SCCR2} ensures that in each iteration, any edge in $D$ can be corrected at most once after being added to $D$. This means that the running time of this algorithm is $O(N|A_C \cup B_C|)$. Note that $|A_C\cup B_C|$ is at most the number of $p$-adjacencies in the graph $G^\star$.
\end{remark}

An example of the execution of Algorithm~\ref{alg: SCCR2} is presented in Appendix~\ref{app: SCCR2}.

\subsection{SCCR algorithm 2: Submodular flow polyhedron} \label{subsection: SCCR1}
Suppose that $C\subseteq [n]$, $A_C\subseteq C^2$, $B_C \subseteq ([n]\setminus C) \times C$, $ComCh_C\subseteq ([n]\setminus C)^2$ and $NoComCh_C\subseteq ([n]\setminus C)^2$. In this subsection, we aim to recover a set $E_C\subseteq [n]^2$ with the properties outlined in Definition~\ref{def: strongly connected component recovery algorithm}, i.e, a set $E_C$ which preserves the $p$-adjacencies in $A_C\cup B_C$, makes $C$ a strongly connected component, assigns a common child in $C$ to any of the pairs in $ComCh_C$, and makes all the pairs in $NoComCh_C$ not have any common children in $C$. Let $\tilde{E}_C  \coloneqq \tilde{A}_C \cup B_C$, where $\tilde{A}_C\subseteq A_C$ contains exactly one of $(a,b)$ or $(b,a)$ for any pair $(a,b)\in A_C$. In other words, in addition to containing all the edges in $B_C$, the set $\tilde{E}_C$ assigns an arbitrary direction to each of the $p$-adjacencies in $A_C$. Our approach is to remove or change the directions of some of the edges in $\tilde{E}_C$ in order for it to satisfy the desired properties. 
\begin{example}
    If $A_{C}$ and $B_{C}$ are chosen to be the sets $A_{C_4}$ and $B_{C_4}$ defined in Example~\ref{example: SCCR} respectively, then one choice for $\tilde{A}_C$ is $\{(1,2),(2,3),(3,4),(1,4),(2,4)\}$.
\end{example}
\begin{definition}\label{def: removable subset}
    Let $Y\subseteq \tilde{E}_C$. The set $Y$ is said to be \textit{removable} if 
    \begin{itemize}
        \item for all $(a,c)\in ComCh_C$, there exists $b\in C$ with $(a,b),(c,b)\in B_C \setminus Y$, 
        \item for all $(a,c)\in NoComCh_C$, there does not exist $b\in C$ with $(a,b),(c,b)\in B_C \setminus Y$, and
        \item for all $(a,c)\in Y$, there exists vertex $b\in C$ such that $(a,b),(c,b) \in A_C \cup B_C$ and $Y\cap \{(a,b),(b,a),(c,b),(b,c)\}=\emptyset$.
    \end{itemize}
\end{definition}
For any removable set $Y\subseteq \tilde{E}_C$, we define a \textit{submodular flow polyhedron}. See Appendix~\ref{app: submodular flow polyhedron} for the definition of a submodular flow polyhedron. 
\begin{definition} \label{def: delta-+}
    For $E\subseteq [n]^2$ and $U\subseteq [n]$,
    \begin{align*}
        &\delta_{E}^-(U) \coloneqq \Set{(a,b)\in E | a \not\in U, b\in U}, \delta_{E}^+(U) \coloneqq \Set{(a,b)\in E | a \in U, b\not\in U}.
    \end{align*}
    Also, for any $x\in \mathbb{R}^E$ and any $E^\prime \subseteq E$, $x(E^\prime) \coloneqq \sum_{e\in E^\prime} x_e$.
\end{definition}
\begin{definition}
    Let $Y\subseteq \tilde{E}_C$ be a removable subset. Then the \textit{submodular flow polyhedron corresponding to $Y$} is denoted by $\mathcal{F}(Y)$ and defined as
    \begin{align*}
        \mathcal{F}(Y) \coloneqq \Set{x \in \mathbb{R}^{\tilde{A}_C \setminus Y} | x\left(\gtm(U)\right) -  x\left(\gtp(U)\right) \leq \left| \gtp(U) \right| - 1 \text{ for all } U\in 2^C \setminus \{\emptyset,C\}.}.
    \end{align*}
    One can verify that $2^C \setminus \{\emptyset,C\}$ is a crossing family and the function $U\mapsto \left| \gtp(U) \right| - 1$ is a crossing submodular function, and hence, $\mathcal{F}(Y)$ is a submodular flow polyhedron.
\end{definition}
\begin{definition}
    Let $Y\subseteq \tilde{E}_C$ be a removable set, and $J\subseteq \tilde{A}_C\setminus Y$. Define $\tilde{E}^{Y,J}_C \coloneqq \tilde{A}^{Y,J}_C \cup B^Y_C$, where 
    \begin{align*}
        &\ayj\coloneqq \Set{(b,a)\in [n]^2 | (a,b) \in J} \cup (\tilde{A}_C \setminus (Y\cup J)), \quad \by \coloneqq B_C\setminus Y,
    \end{align*}
    i.e., $\eyj$ is obtained from $\tilde{E}_C$ by removing the edges in $Y$ and changing any pair $(a,b)\in J$ to $(b,a)$.
\end{definition}

\begin{lemma} \label{lemma: submodular flow}
    Let $Y\subseteq \tilde{E}_C$ be a removable set, and $J\subseteq \tilde{A}_C\setminus Y$. 
    Then $C$ is a strongly connected component in the graph $\left([n],\eyj\right)$ if and only if $J=\Set{e\in \tilde{A}_C\setminus Y|x_e=-1}$ for some $x\in \mathcal{F}(Y)\cap \{-1,0\}^{\tilde{A}_C\setminus Y}$.
\end{lemma}
The proof of this lemma can be found in Appendix~\ref{app: proofs of submodular}.
\begin{definition}
Let $Y\subseteq \tilde{E}_C$ be a removable set. For any $(a,b,c)\in [n]^3$ with $c\in C$, $(a,b),(c,b)\in A_C\cup B_C$ and $Y\cap \{(a,b),(b,a),(c,b),(b,c)\}=\emptyset$, the \textit{weight vector corresponding to $(a,b,c)$} is a vector in $\mathbb{R}^{\tilde{A}_C\setminus Y}$ which is denoted by $w_{Y,(a,b,c)}$ and defined as follows:
\begin{align*}
  \text{ For all $e\in \tilde{A}_C \setminus Y$,}\quad 
    w_{Y,(a,b,c)}^{(e)} \coloneqq \begin{cases}
      -1 &  \text{if $e=(b,a)$ or $e=(b,c)$.} \\
      1 & \text{if $e=(a,b)$ or $e=(c,b)$.}\\
      0 & \text{otherwise.}
    \end{cases}
\end{align*}

\end{definition} 
\begin{theorem} \label{thm: SCCR1}
    Suppose that $C\subseteq [n]$, $A_C\subseteq C^2$, $B_C \subseteq ([n]\setminus C) \times C$, $ComCh_C\subseteq ([n]\setminus C)^2$ and $NoComCh_C\subseteq ([n]\setminus C)^2$ are given. Then there exists $E_C\subseteq [n]^2$ with the properties outlined in Definition~\ref{def: strongly connected component recovery algorithm} and containing at most one of $(a,b)$ or $(b,a)$ for any $a,b\in [n]$ if and only if there exist a removable set $Y\subseteq \tilde{E}_C$ and a vertex $x$ of the polyhedron $\mathcal{F}(Y)\cap [-1,0]^{\tilde{A}_C\setminus Y}$ such that the following hold:
    \begin{enumerate}
        \item For all $(a,b,c)\in [n]^3$ with $c\in C$, $(a,b),(c,b)\in A_C\cup B_C$, $Y\cap \{(a,b),(b,a),(c,b),(b,c)\}=\emptyset$ and $(a,c)\not\in A_C\cup B_C$, \label{thm: SCCR1: property 1}
        \begin{align}\label{eq: SCCR1 property 1}
            w_{Y,(a,b,c)}^Tx < \max_{y\in \{-1,0\}^{\tilde{A}_C\setminus Y}} w_{Y,(a,b,c)}^T y.
        \end{align}
        \item For all $(a,c)\in Y$, there exists $b\in C$ with $(a,b),(c,b)\in A_C \cup B_C$ and $Y\cap \{(a,b),(b,a),(c,b),(b,c)\}=\emptyset$ such that \label{thm: SCCR1: property 2}
        \begin{align}\label{eq: SCCR1 property 2}
            w_{Y,(a,b,c)}^Tx = \max_{y\in \{-1,0\}^{\tilde{A}_C\setminus Y}} w_{Y,(a,b,c)}^T y.
        \end{align}
    \end{enumerate}
    Moreover, if a removable set $Y$ and a vertex $x$ of $\mathcal{F}(Y)\cap [-1,0]^{\tilde{A}_C\setminus Y}$ satisfying the above properties exist, then 
    $\eyjx$ satisfies the properties listed for $E_C$ in Definition~\ref{def: strongly connected component recovery algorithm}, where 
    \begin{align*}
        J(x) =\Set{e\in \tilde{A}_C\setminus Y| x_e=-1}.
    \end{align*}
\end{theorem}
The proof of this theorem can be found in Appendix~\ref{app: proofs of submodular}. 
Theorem~\ref{thm: SCCR1} leads us to our second SCCR algorithm, Algorithm~\ref{alg: SCCR1}. 
\begin{algorithm}[h]
\caption{SCCR algorithm 2: Submodular flow polyhedron}\label{alg: SCCR1}
\begin{algorithmic}[1]
\Require{Sets $C\subseteq[n]$, $A_C\subseteq C^2$, $B_C \subseteq ([n]\setminus C) \times C$, $ComCh_C\subseteq ([n]\setminus C)^2$ and $NoComCh_C\subseteq ([n]\setminus C)^2$.}
\Ensure{Either a declaration of failure, or a subset of  $[n]^2$ with the properties described in Definition~\ref{def: strongly connected component recovery algorithm}.}    
\State Let $\tilde{E}_C  \coloneqq \tilde{A}_C \cup B_C$, where $\tilde{A}_C\subseteq A_C$ contains exactly one of $(a,b)$ or $(b,a)$ for any pair $(a,b)\in A_C$.
\For {removable subsets $Y\subseteq \tilde{E}_C$}
    \For {vertices $x$ of the polyhedron $\mathcal{F}(Y)\cap [-1,0]^{\tilde{A}_C\setminus Y}$}
        \If {$x$ satisfies properties~\ref{thm: SCCR1: property 1} and~\ref{thm: SCCR1: property 2} of Theorem~\ref{thm: SCCR1}}
            \State \Return $\eyjx$.
        \EndIf
    \EndFor
\EndFor
\State \Return "failed".
\end{algorithmic}
\end{algorithm}

\begin{remark}[Running time]
There are $O\left(2^{|A_C \cup B_C|}\right)$ choices for the removable set $Y$. On the one hand, for any removable set $Y\subseteq \tilde{E}_C$ and any $(a,b,c)\in [n]^3$, the time complexity of computing the right hand side of \eqref{eq: SCCR1 property 1} and \eqref{eq: SCCR1 property 2} is $O\left( |\tilde{A}_C\setminus Y|\ 2^{|\tilde{A}_C \setminus Y}|\right)  \leq  O \left( |A_C|\  2^{|A_C|} \right)$. On the other hand, the number of defining inequalities of $\mathcal{F}(Y)\cap [-1,0]^{\tilde{A}_C\setminus Y}$ is $2^{|C|} + 2|\tilde{A}_C \setminus Y|-2$. Therefore, by \cite{Avis1992APA}, the time complexity of enumerating all of the vertices of $\mathcal{F}(Y)\cap [-1,0]^{\tilde{A}_C\setminus Y}$ is 
\begin{align*}
   & O\left(\left(2^{|C|} +2 |\tilde{A}_C \setminus Y| - 2\right) |\tilde{A}_C \setminus Y | \left( 2^{|C|} + 3|\tilde{A}_C \setminus Y| - 2\right) \binom{2^{|C|} + 2|\tilde{A}_C \setminus Y| - 2}{|\tilde{A}_C\setminus Y|} \right) \\
   \leq & O \left(|A_C|\left(2^{|C|} + 3|A_C| \right)^2 \left(2^{|C|} + 2|A_C|\right)^{|A_C|} \right).
\end{align*}
Now given that for each vertex, the time complexity of computing the left hand side of \eqref{eq: SCCR1 property 1} and \eqref{eq: SCCR1 property 2} is $O\left(|A_C|\right)$, the overall  time complexity of Algorithm~\ref{alg: SCCR1} is 
\begin{align*}
    O\left(|A_C|^2 {|A_C\cup B_C|} ^2\left(2^{|C|} + 3|A_C| \right)^2 \left(2^{|C|} + 2|A_C|\right)^{|A_C|}  2^ {|A_C \cup B_C|} \right).
\end{align*}
\end{remark}

\section{Simulations} \label{simulations}
In this section, we provide experimental results of Algorithm~\ref{alg: greedy Markov equivalence class discovery}, Algorithm~\ref{alg: SCCR2}, and the  combination of Algorithms~\ref{alg: greedy Markov equivalence class discovery},~\ref{alg: graph discovery} and~\ref{alg: SCCR2}. Our simulation studies compare the performance of these algorithms on graphs with $7$ vertices and different levels of sparsity as well as graphs with the same levels of sparsity and different numbers of vertices. The code related to these simulations is available in~\cite{githubrepo}.

\subsection{Discovering the Markov equivalence class} \label{subsection: simulation: greedy Markov equivalence discovery}
\begin{figure} [H]
    \centering
    \begin{subfigure}{1\textwidth}
    \centering
        \includegraphics[width=0.8\textwidth]{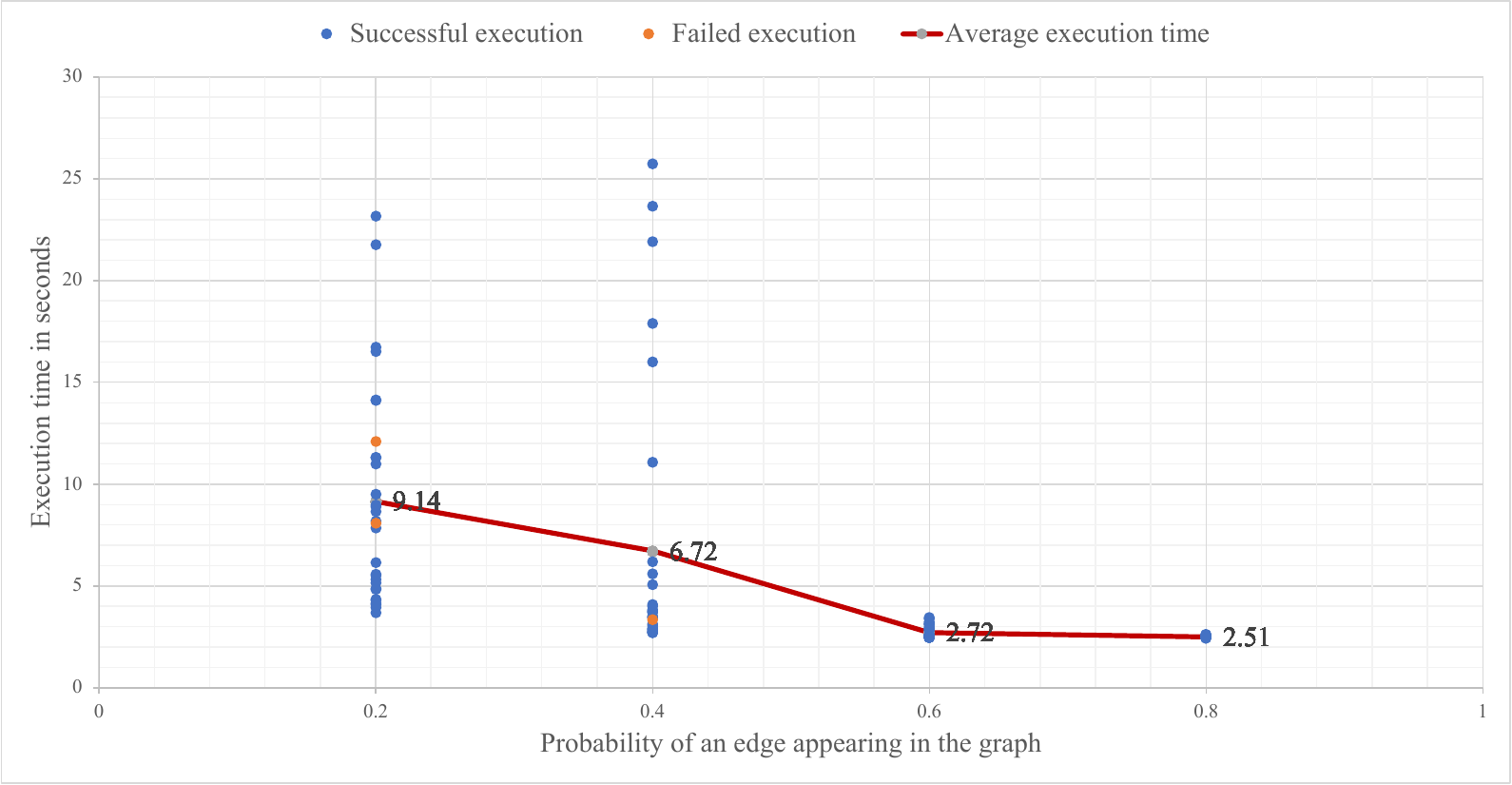}
        \caption{Results of running Algorithm~\ref{alg: greedy Markov equivalence class discovery} on 120 random graphs with 7 vertices and different levels of sparsity. Thirty graphs were tested for each of the probabilities 0.2, 0.4, 0.6, and 0.8. From left to right, the success rates are 0.93, 0.97, 1 and 1. }
    \label{fig:opt-n=7}
    \end{subfigure}
    \begin{subfigure}{1\textwidth}
    \centering
        \includegraphics[width=0.8\textwidth]{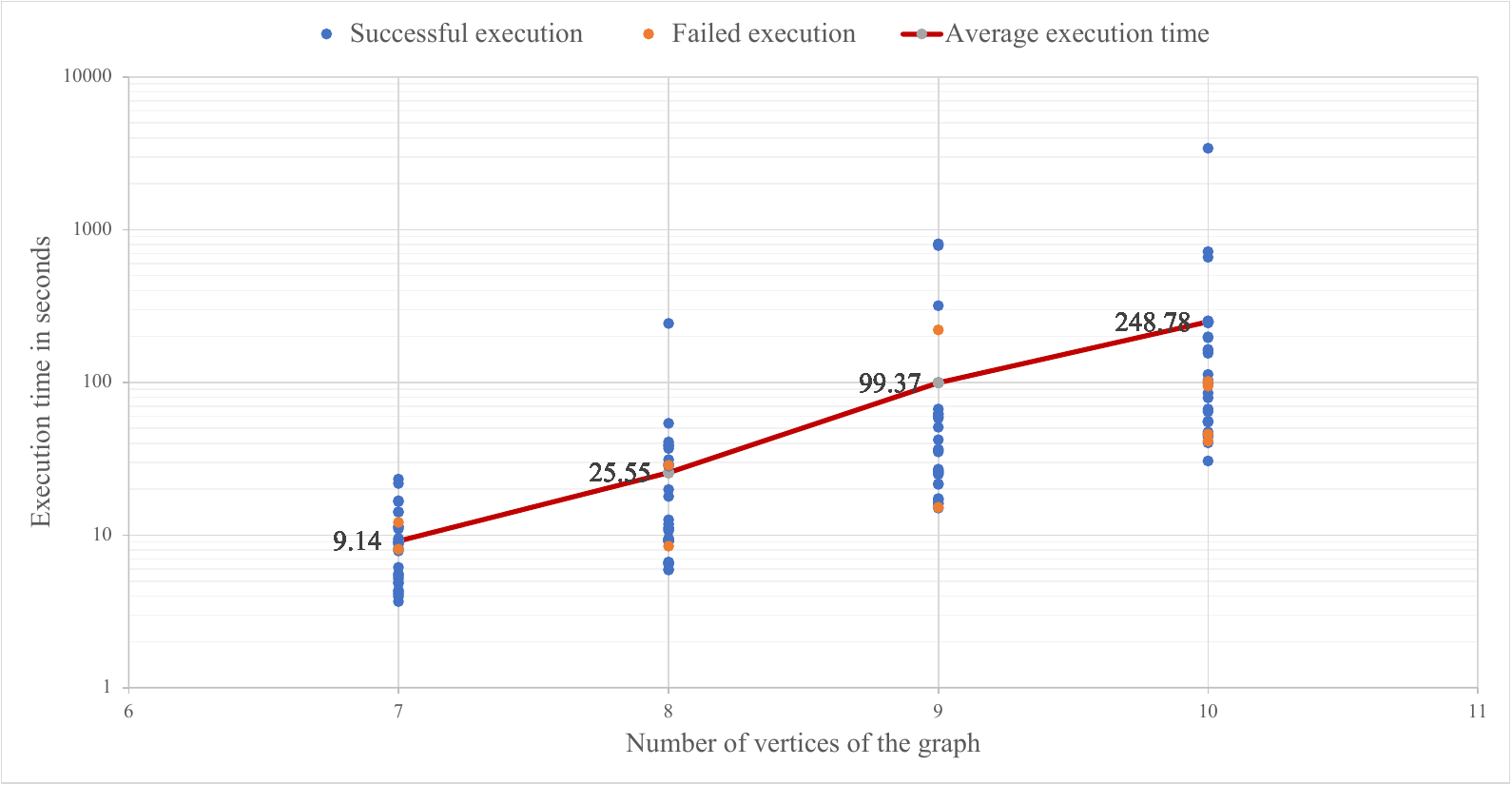}
        \caption{Results of running Algorithm~\ref{alg: greedy Markov equivalence class discovery} on 120 random graphs with different numbers of vertices $n$ where each edge appears in the graph with probability 0.2. Thirty graphs were tested for each $n$. From left to right, the success rates are 0.93, 0.93, 0.93 and 0.87.}
    \label{fig:opt-p=0.2}
    \end{subfigure}
    \begin{subfigure}{1\textwidth}
    \centering
        \includegraphics[width=0.8\textwidth]{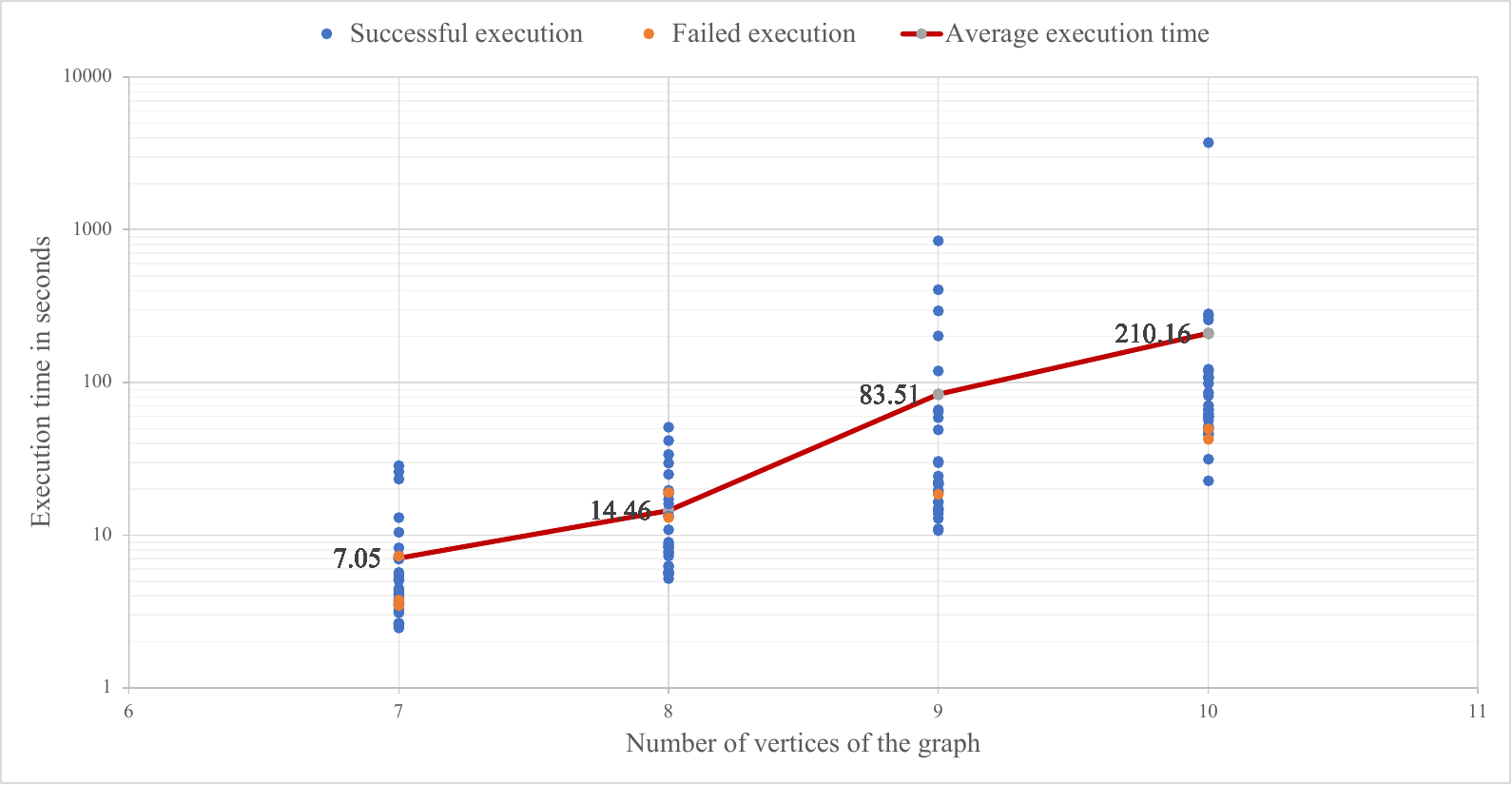}
        \caption{
        Results of running Algorithm~\ref{alg: greedy Markov equivalence class discovery} on 120 random graphs with different numbers of vertices $n$ where each edge appears in the graph with probability 0.3. Thirty graphs were tested for each $n$. From left to right, the success rates are 0.90, 0.93, 0.97 and 0.93.}
    \label{fig:opt-p=0.3}
    \end{subfigure}
    \vspace{-4mm}
    \caption{Simulation results on Algorithm~\ref{alg: greedy Markov equivalence class discovery}.}
     \label{fig:opt}
\end{figure}

This subsection describes our simulation results for Algorithm~\ref{alg: greedy Markov equivalence class discovery}. The tests were performed as follows:

For each $(n,p) \in \left(\{7\} \times \{0.2,0.4,0.6,0.8\}\right) \cup \left(\{7,8,9,10\}\times \{0.2,0.3\}\right)$, we generated 30 random graphs with $n$ vertices, where the probability of each edge appearing in the graph is $p$. These graphs were generated according to the Erdős–Rényi model and using the \texttt{R} library \texttt{igraph}. The set of all $d$-separations satisfied by each graph was then generated and given to  Algorithm~\ref{alg: greedy Markov equivalence class discovery} as part of its input. We chose $M$ (i.e., the number of different initial partially ordered partitions) to be $3$ with the following initial partially ordered partitions: (Non-trivial pairs in each partial order (if any) are shown below.)
\begin{align*}
    &\sP_1 = \left\{ \{1,\ldots,n\} \right\},\\
    &\sP_2 = \left\{\left\{1,\ldots,\left\lfloor\frac{n}{2}\right\rfloor \right\} , \left\{\left \lfloor\frac{n}{2} \right \rfloor +1,\ldots,n  \right\} \right\}, \left\{1,\ldots,\left\lfloor\frac{n}{2}\right\rfloor \right\} \leq_{\pi_2} \left\{\left \lfloor\frac{n}{2} \right \rfloor +1,\ldots,n  \right\} ,\\
    &\sP_3 = \left\{ \{1\}, \{2\},\ldots,\{n\} \right\}.
\end{align*}
The integers $N$ corresponding to $n=7$, $n=8$, $n=9$ and $n=10$ were chosen to be $30$, $30$, $40$ and $50$ respectively, where $N$ is the maximum number of adjacent partially ordered partitions with the same score which the algorithm is allowed to observe before returning a partially ordered partition. Figure~\ref{fig:opt} shows our simulation results. 

\subsection{SCCR algorithm 1: Construct and correct} \label{subsection: simulation SCCR}
In this subsection, we show the experimental results related to Algorithm~\ref{alg: SCCR2}, which reconstructs the edges within and pointing to a given strongly connected component. 

\begin{figure}[H]
    \centering
    \begin{subfigure}{1\textwidth}
    \centering
        \includegraphics[width=0.8\textwidth]{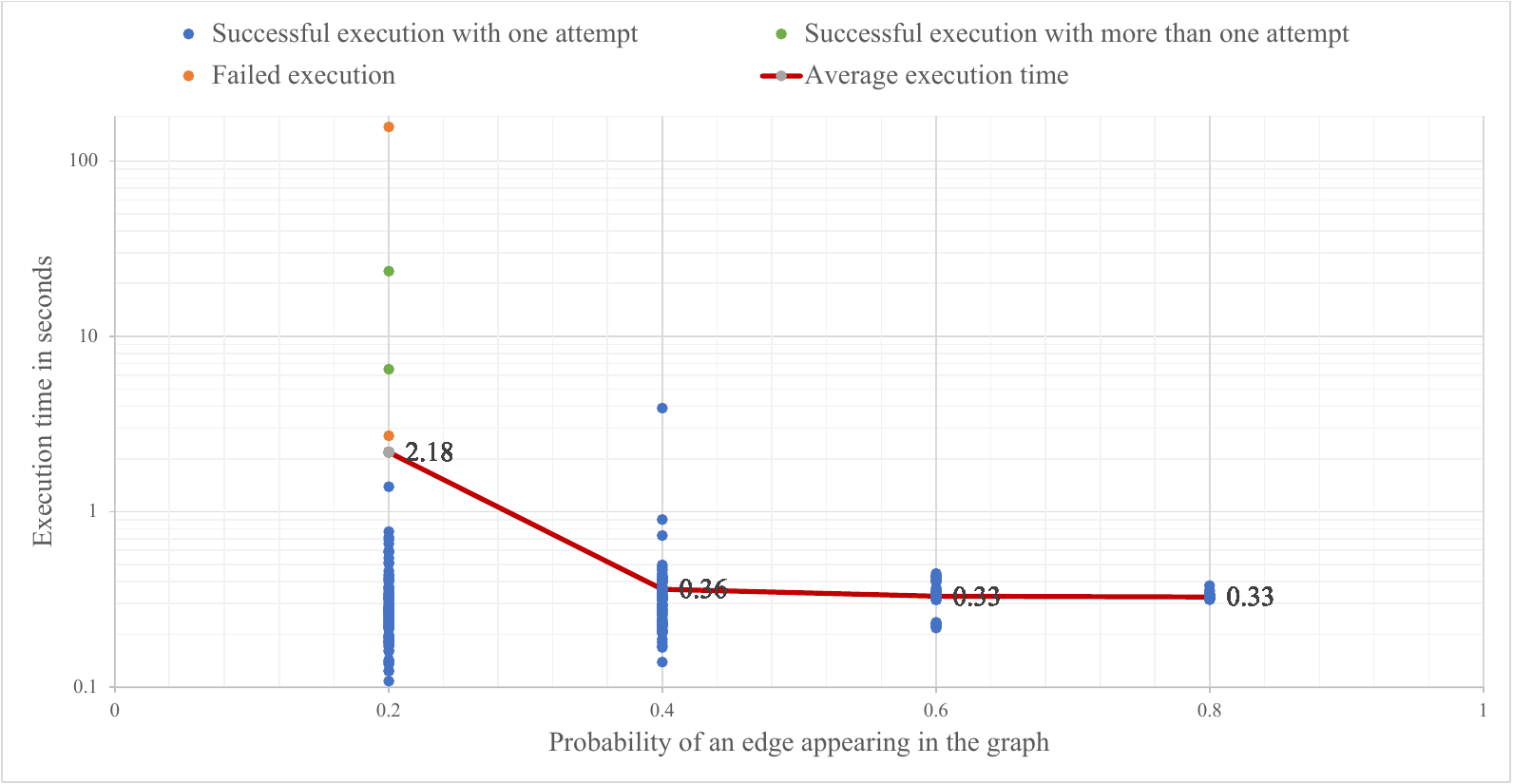}
        \caption{Results of running Algorithm~\ref{alg: SCCR2} on 400 random graphs with 7 vertices and different levels of sparsity. One-hundred graphs were tested for each of the probabilities 0.2, 0.4, 0.6, and 0.8. From left to right, the success rates are 0.98, 1, 1 and 1. }
    \label{fig:SCCR-n=7}
    \end{subfigure}
    \begin{subfigure}{1\textwidth}
    \centering
        \includegraphics[width=0.8\textwidth]{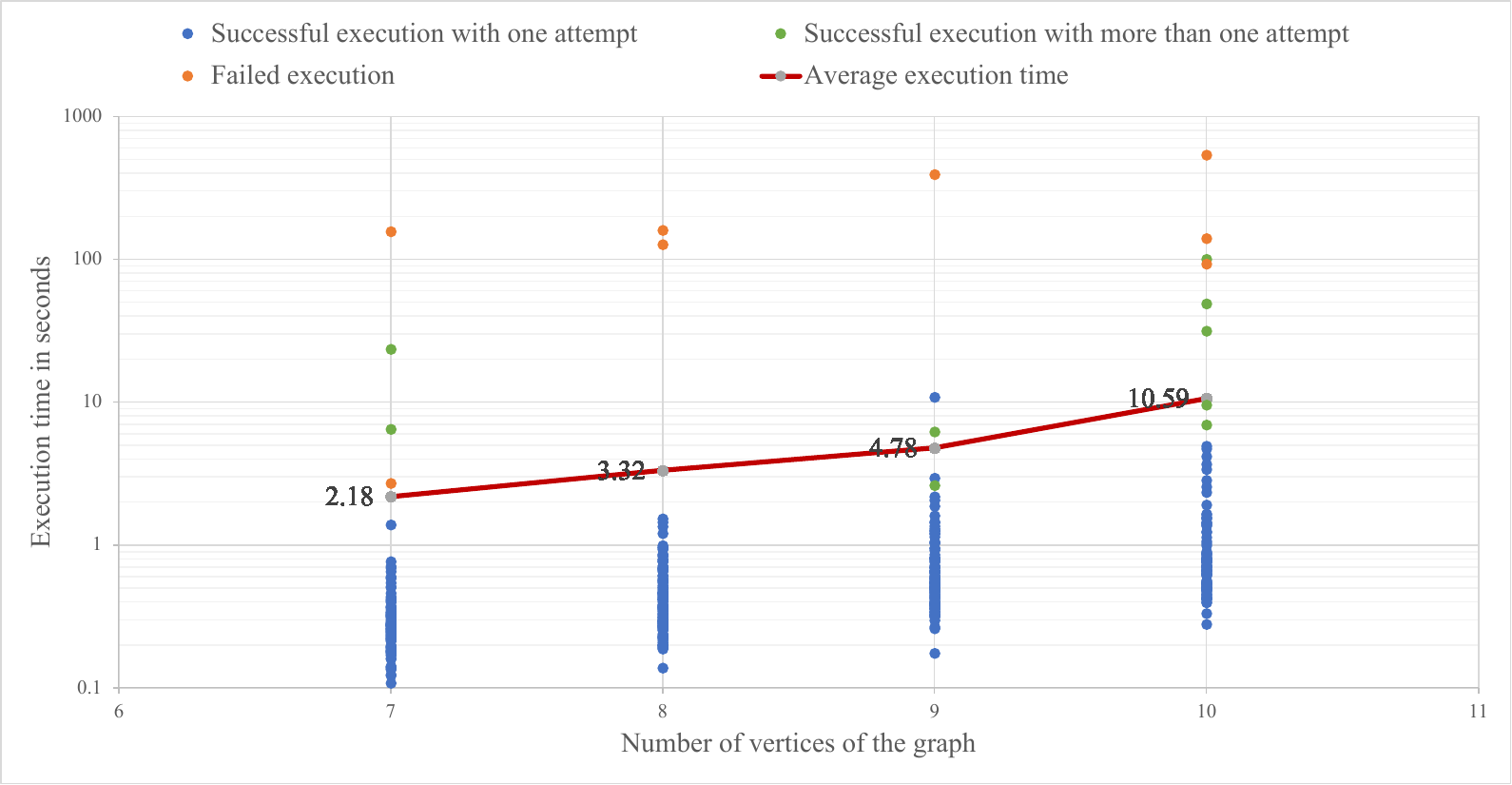}
        \caption{Results of running Algorithm~\ref{alg: SCCR2} on 400 random graphs with different numbers of vertices $n$ where each edge appears in the graph with probability 0.2. One-hundred graphs were tested for each $n$. From left to right, the success rates are 0.98, 0.98, 0.99 and 0.97.}
    \label{fig:SCCR-p=0.2}
    \end{subfigure}
    \begin{subfigure}{1\textwidth}
    \centering
        \includegraphics[width=0.8\textwidth]{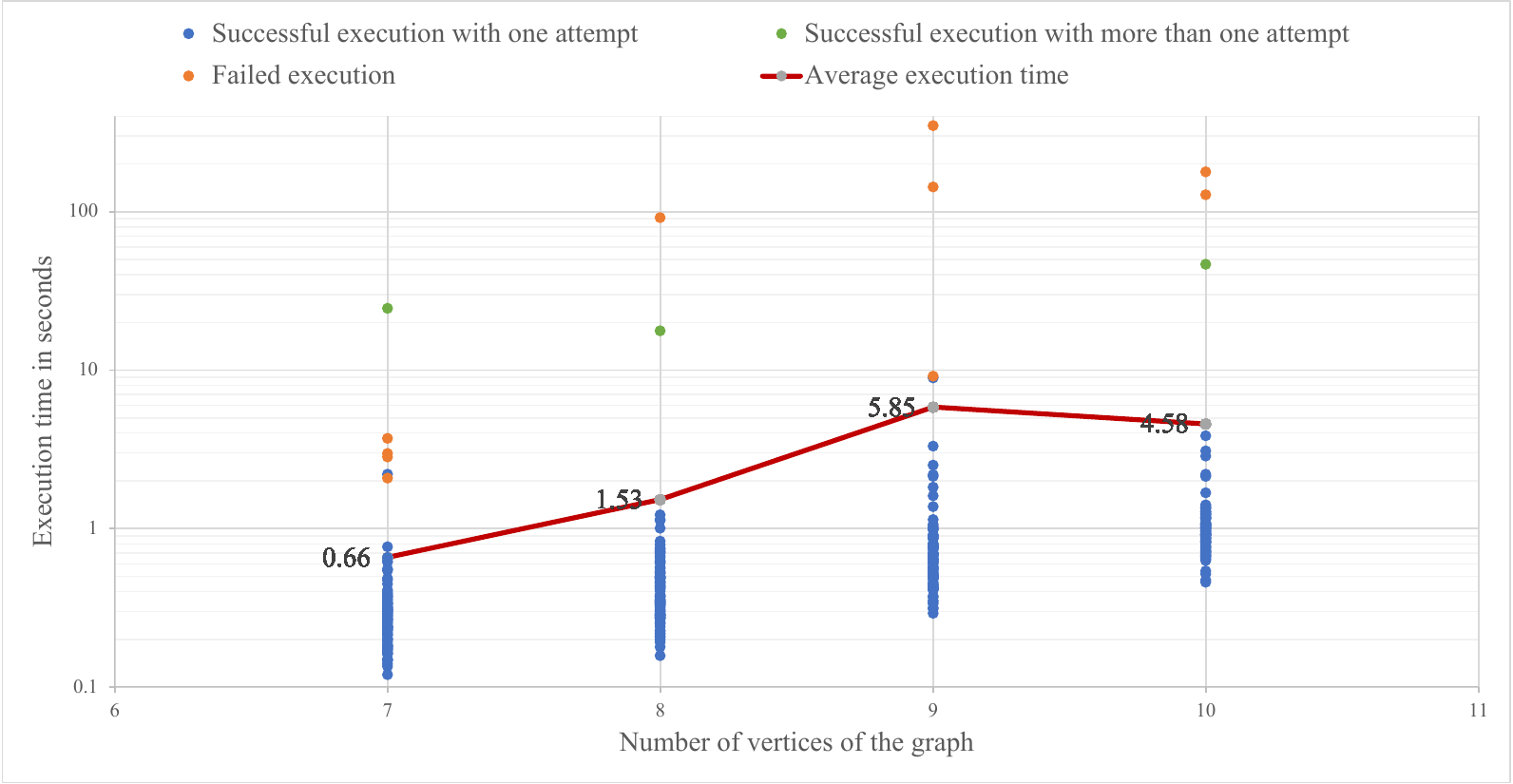}
        \caption{Results of running Algorithm~\ref{alg: SCCR2} on 400 random graphs with different numbers of vertices $n$ where each edge appears in the graph with probability 0.3. One-hundred graphs were tested for each $n$. From left to right, the success rates are 0.96, 0.99, 0.97 and 0.98.}
    \label{fig:SCCR-p=0.3}
    \end{subfigure}
    \vspace{-4mm}
    \caption{Simulation results on Algorithm~\ref{alg: SCCR2}.}
    \label{fig:SCCR}
\end{figure}

\begin{figure}[H]
    \centering
    \includegraphics[width=0.8\textwidth]{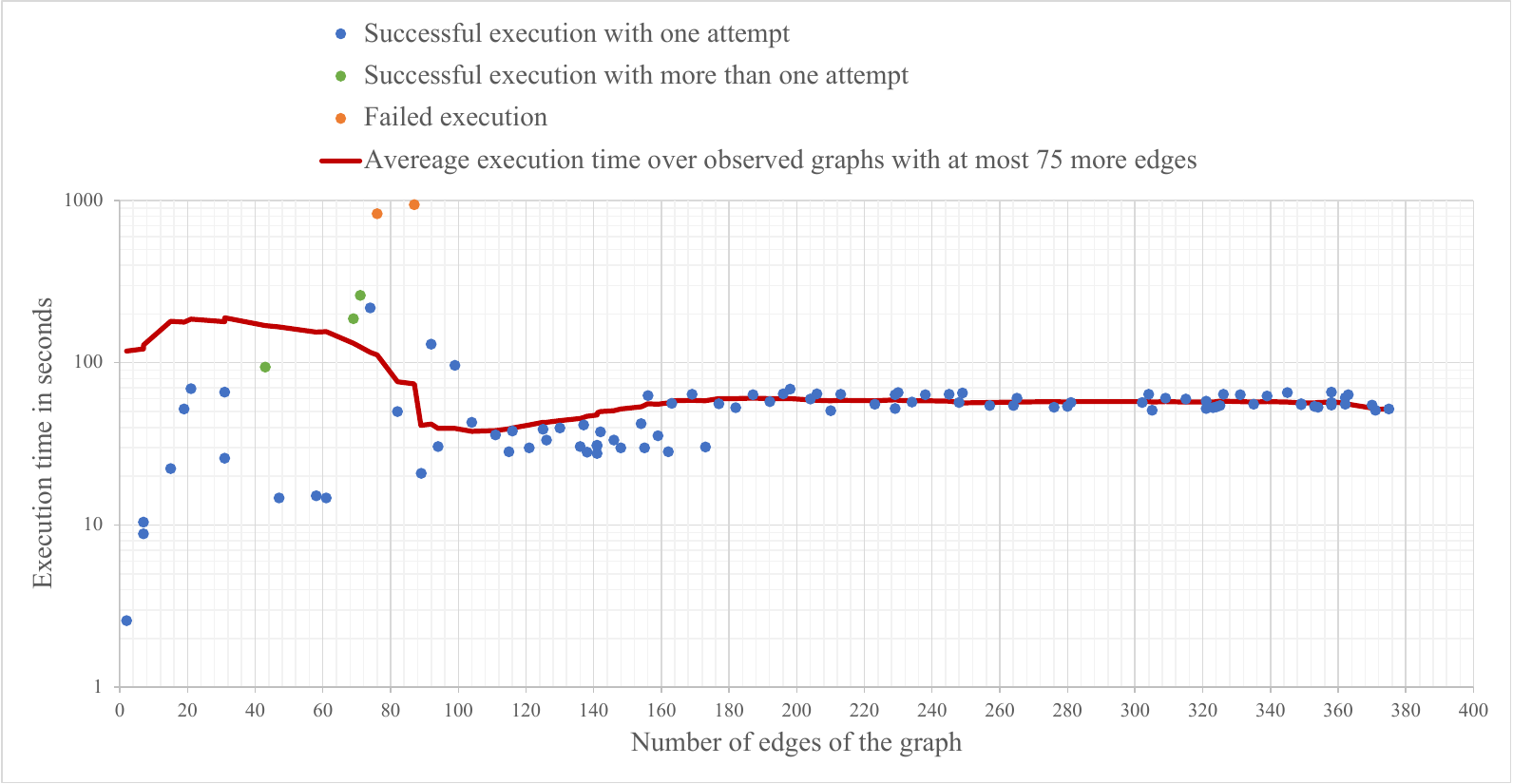}
    \caption{Results of running Algorithm~\ref{alg: SCCR2} on 100 random graphs with 20 vertices. The success rate is 0.98 and the average execution time is 71.90 seconds.}
    \label{fig:SCCR-2}
\end{figure}

For each $(n,p) \in \left(\{7\} \times \{0.2,0.4,0.6,0.8\}\right) \cup \left(\{7,8,9,10\}\times \{0.2,0.3\}\right)$, we generated 100 random graphs with $n$ vertices, where the probability of each edge appearing in the graph is $p$. These graphs were generated according to the Erdős–Rényi model and using the \texttt{R} library \texttt{igraph}. For each random graph, we considered the partially ordered partition $(\sP_0,\pi_0)$ associated with the graph. For each $C\in \sP_0$, we defined $A_C$, $B_C$, $ComCh_C$ and $NoComCh_C$ as described in Conjecture~\ref{conj: SCCR2}, and then fed Algorithm~\ref{alg: SCCR2} the sets $C,A_C,B_C,ChomCh_C,NoComCh_C$ and the integer $N=100$ as the input. Due to the random nature of this algorithm, different executions may lead to different outputs. We allowed for at most 20 attempts on each graph.
In cases where in one of the attempts the algorithm outputted a set $E_C$ for all $C\in \sP_0$, we declared success and otherwise, declared failure.
See Figure~\ref{fig:SCCR} for the results of these simulations.
We also performed the same experiment on 100 random graphs with 20 vertices and different numbers of edges. The results are shown in Figure~\ref{fig:SCCR-2}.

For a fixed number of vertices, there seems to be a decreasing trend in the execution time of the algorithm as the graph gets denser while for a fixed probability $p$, the execution time tends to increase as the number of vertices increases. In all our experiments, the success rate has been 0.96 or above.
\subsection{Discovering a Markov equivalent graph} \label{subsection: simulation: greedy Markov equivalent graph discovery}
\begin{figure}[H]
    \centering
    \begin{subfigure}{1\textwidth}
    \centering
        \includegraphics[width=0.8\textwidth]{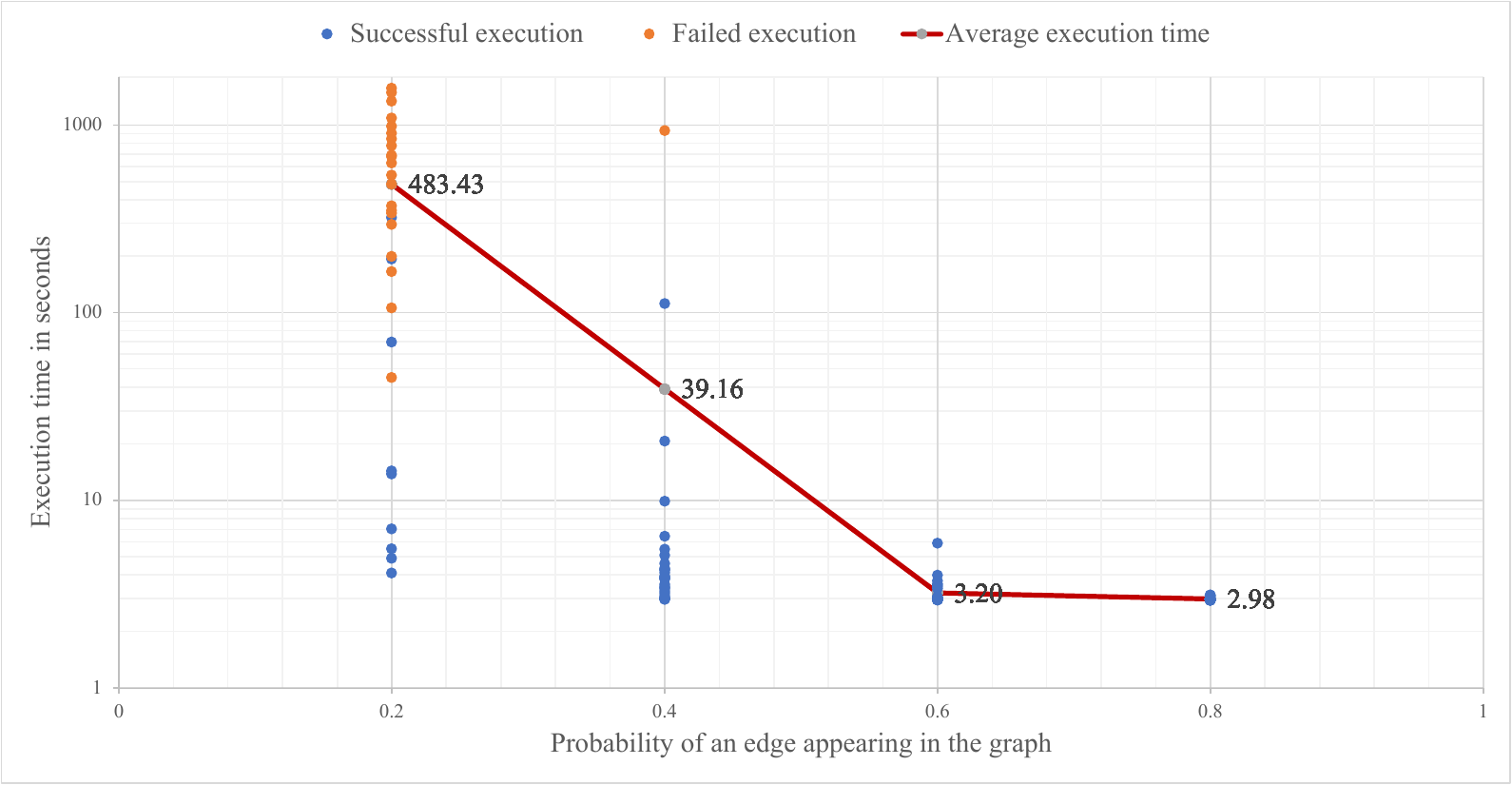}
        \caption{Results of running a combination of Algorithms~\ref{alg: greedy Markov equivalence class discovery} and~\ref{alg: graph discovery} on the sets of $d$-separations implied by 120 random graphs with 7 vertices and different levels of sparsity. Thirty graphs were tested for each of the probabilities 0.2, 0.4, 0.6, and 0.8. From left to right, the success rates are 0.3, 0.97, 1 and 1. }
    \label{fig:graph discovery-n=7}
    \end{subfigure}
   
    \begin{subfigure}{1\textwidth}
    \centering
        \includegraphics[width=0.8\textwidth]{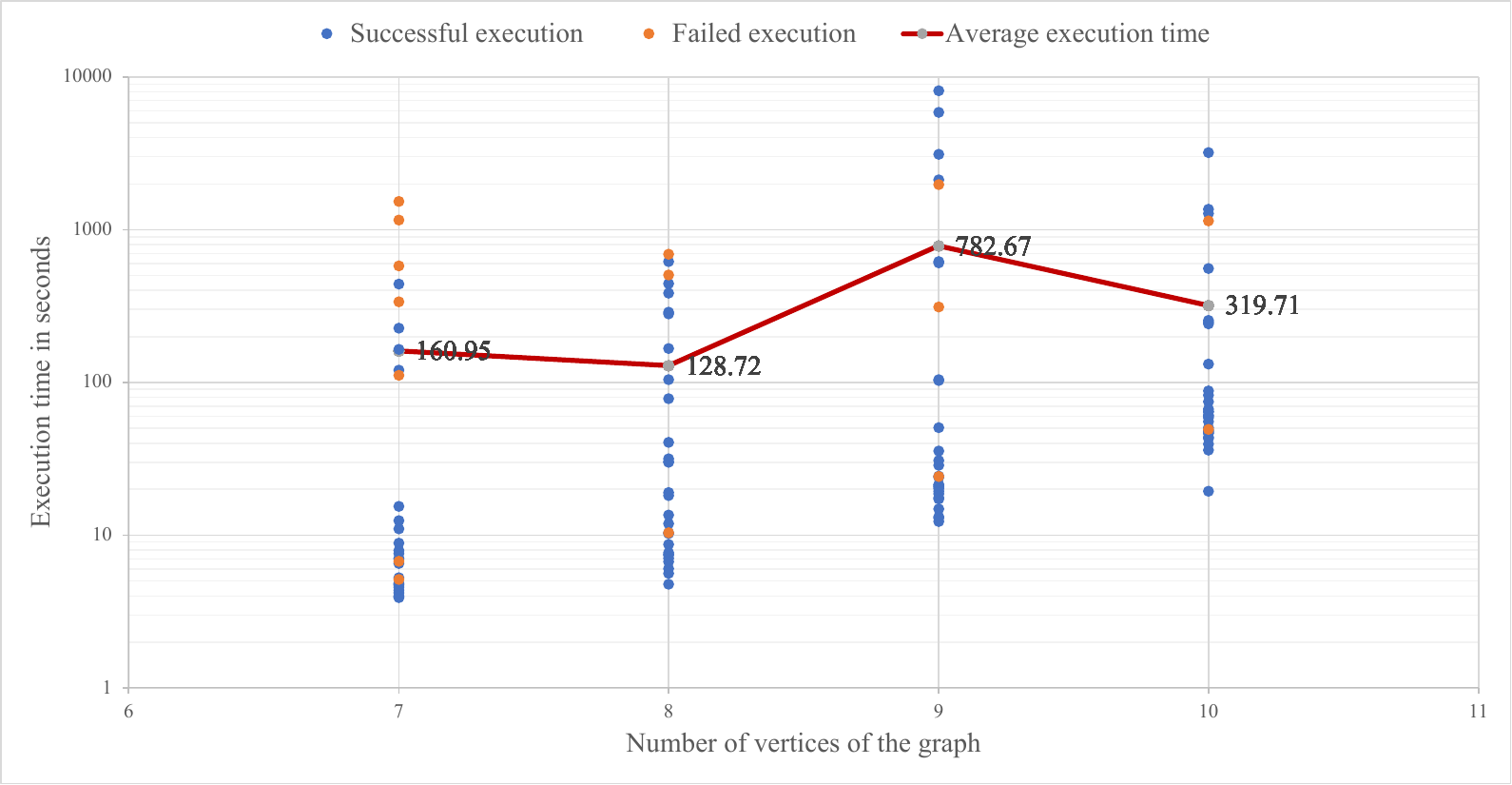}
        \caption{Results of running a combination of Algorithms~\ref{alg: greedy Markov equivalence class discovery} and~\ref{alg: graph discovery} on the sets of $d$-separations implied by 120 random graphs with different numbers of vertices $n$ where each edge appears in the graph with probability 0.3. Thirty graphs were tested for each $n$. From left to right, the success rates are 0.77, 0.90, 0.90 and 0.93.}
    \label{fig:graph discovery-p=0.3}
    \end{subfigure}
    \vspace{-7mm}
    \caption{Simulation results on the combination of Algorithms~\ref{alg: greedy Markov equivalence class discovery} and~\ref{alg: graph discovery}.}
    \label{fig:graph discovery}
\end{figure}

In this subsection, we show the experimental results of the combination of Algorithms~\ref{alg: greedy Markov equivalence class discovery} and~\ref{alg: graph discovery}. Note that this combination takes the set of conditional independence statements inferred from the observed samples as the input and outputs a Markov equivalent graph to $G^\star$. Our experiment was performed as follows:

For each $(n,p) \in \left(\{7\} \times \{0.2,0.4,0.6,0.8\}\right) \cup \left(\{7,8,9,10\}\times \{0.2,0.3\}\right)$, we generated 30 random graphs with $n$ vertices, where the probability of each edge appearing in the graph is $p$. These graphs were generated according to the Erdős–Rényi model and using the \texttt{R} library \texttt{igraph}. For each random graph, the set of all $d$-separations arising from the graph were generated and given to Algorithm~\ref{alg: greedy Markov equivalence class discovery}. The other parameters of this algorithm were chosen as outlined in Subsection~\ref{subsection: simulation: greedy Markov equivalence discovery}. The SCCR algorithm used in Algorithm~\ref{alg: graph discovery} was SCCR algorithm 1, i.e. Algorithm~\ref{alg: SCCR2}, with $N=100$. Note that only one attempt of this algorithm was allowed on each optimal partially ordered partition. While Algorithm~\ref{alg: graph discovery} continues testing optimal partially ordered partitions until it finds one which has a corresponding graph, in our experiment, we restricted ourselves to testing at most 300 optimal partially ordered partitions before declaring failure. 
See Figure~\ref{fig:graph discovery} for the results.


Figure~\ref{fig:graph discovery-n=7} suggests for sparser graphs, one seems to need to test more optimal partially ordered partitions before being able to recover a Markov equivalent graph from one of them. One alternative method for sparser graphs would be to first perform Algorithm~\ref{alg: greedy Markov equivalence class discovery} to obtain an optimal partially ordered partition $(\sP,\pi)$, and hence, a characterization of the Markov equivalence class, and then simply run a brute-force search over all graphs whose edges are subsets of $E^{(1)}_{(\sP,\pi)}$ to find a graph in the Markov equivalence class.

\section{Discussion} \label{discussion}
We considered the question of causal discovery given the conditional independence statements of a distribution $\mathbb{P}$ which is Markov and faithful to an unknown graph $G^\star$. Assuming there are no latent variables, we proposed two algorithms; one recovers the Markov equivalence class of $G^\star$ via a hybrid approach and the other uses the output of the first algorithm to recover a graph in this class. 

In addition to proving the two conjectures proposed throughout the paper, i.e. Conjectures~\ref{conjecture: greedy optimization} and~\ref{conj: SCCR2}, we believe a future direction of study could involve relaxing the assumption of $\mathbb{P}$ being Markov and faithful to $G^\star$. Unlike SEMs whose causal graphs are acyclic, for an arbitrary SEM, it is only known  that in linear and discrete cases the corresponding distribution is in the graphical model of its causal graph~\cite[Theorem 6.3]{Bongers_2021}. On the other hand, faithfulness is generally considered to be a restrictive assumption. Hence, replacing any of the above-mentioned assumptions with weaker ones would be interesting. 

Another modification to the algorithm that could be studied is to move through graphs rather than through partially ordered partitions in order to find a Markov equivalent graph to $G^\star$.

Furthermore, extending the methods proposed in the paper to the case where latent variables might be present is another question of interest for future work.

\vspace{-6mm}

\appendix
\section{Proofs for Section~\ref{poset discovery}}\label{app:proofs_sec3}
\begin{proof}[Proof of Proposition~\ref{prop: Rich1}]
Let $a,b\in [n]^2$ and $a\padj b$. Then by \cite[Lemma 1]{richardson1997char}, $a\not\independent b \mid Z$ for any set $Z\subseteq [n] \setminus \{a,b\}$. For every $(\sP,\pi)\in \sS$, take $Z\coloneqq \bigcup \Set{C\in \sP |  C\leq_\pi \max\{C_{a,\sP}, C_{b,\sP}\}}\setminus \{a,b\}$. So, the set of $p$-adjacencies in $G^\star$ is contained in $E^{(1)}_{(\sP,\pi)}$.

Now let $(\sP_0,\pi_0)$ be the partially ordered partition associated with $G^\star$. We will proceed to show that for all $(a,b)\in [n]^2$, if $(a,b)\in E^{(1)}_{(\sP_0,\pi_0)}$, then $a\padj b$. Combined with the first paragraph of the proof, this implies that $E^{(1)}_{(\sP_0,\pi_0)}$ is exactly equal to the set of $p$-adjacencies in $G^\star$.

Suppose that $(a,b)\in E^{(1)}_{(\sP_0,\pi_0)}$. Then $a\not\independent b \mid Z$, where 
\begin{align*}
    Z \coloneqq \bigcup \Set{C\in \sP_0 |  C\leq_{\pi_0} \max\{C_{a,\sP_0}, C_{b,\sP_0}\}}\setminus \{a,b\}.
\end{align*}
So, there exists a path $P: a=v_1,e_1,\ldots,e_{t-1},v_t=b$ in $G^\star$ such that every collider on $P$ either is in $Z$ or have a descendant in $Z$, and for all $\ell\in\{2,\ldots,t-1\}$, if $v_\ell\in Z$, then $v_\ell$ is a collider on the path $P$.

We claim that for all $\ell\in \{2,\ldots,t-1\}$, $v_\ell\in Z$. By contradiction, assume the claim does not hold. Define
\begin{align*}
    p\coloneqq \min\Set{\ell \in \{2,\ldots,t-1\} | v_\ell \not \in Z }.
\end{align*}
Note that for all $k_1\in [n]\setminus \{a,b\}$, if $k_1$ has a descendant $k_2$ in $Z$, then $C_{k_1,G^\star} \leq_{G^\star} C_{k_2,G^\star}$, and thus, $C_{k_1,\sP_0} \leq_{\pi_0} C_{k_2,\sP_0}$. On the other hand, since $k_2\in Z$, $C_{k_2,\sP_0} \leq_{\pi_0} \max\{C_{a,\sP_0},C_{b,\sP_0}\}$. So, $C_{k_1,\sP_0} \leq_{\pi_0}  \max\{C_{a,\sP_0},C_{b,\sP_0}\}$, and as a result, $k_1 \in Z$. So, considering that $v_{p-1} \in Z \cup \{a\}$, if $e_{p-1} = (v_p,v_{p-1})$, then $v_p \in Z$, which contradicts the definition of $p$. Hence, $e_{p-1} = (v_{p-1}, v_p)$. Now if $v_p$ is a collider on $P$, then either $v_p \in Z$ or $v_p$ has a descendant in $Z$, both of which mean that $v_p\in Z$, and are impossible. So, $e_p=(v_p,v_{p+1})$. Let
\begin{align*}
q \coloneqq \max \Set{\ell \in \{p,\ldots,t-1\} | e_{\ell^\prime} = (v_{\ell^\prime},v_{\ell^\prime+1}) \text{ and $v_{\ell^\prime} \not \in Z$ for all $\ell^\prime \in \{p,\ldots,\ell\}$.} }.
\end{align*}
Then if $v_{q+1} \in Z\cup \{b\}$, it means that $v_q$ has a descendant in $Z\cup \{b\}$, and thus, is in $Z$, which is not true. So, $v_{q+1}\not\in Z\cup \{b\}$ and given the definition of $q$, one can conclude that $e_{q+1} = (v_{q+2},v_{q+1})$. Therefore, $v_{q+1}$ is a collider on $P$, which implies that $v_{q+1} \in Z$. This contradiction proves the initial claim.

Every inner vertex on the path $P$ belonging to $Z$ means that every vertex on the path, except for $v_1$ and $v_t$, is a collider. So, either $t=2$ or $t=3$. In the first case, $a$ and $b$ are adjacent, and thus, $p$-adjacent. In the second case, $v_2$ is a common child of $a$ and $b$. So, $v_2$ is a collider on $P$, as a result of which, $v_2\in Z$. Thus, $C_{v_2,\sP_0} \leq_{\pi_0} \max\{C_{a,\sP_0},C_{b,\sP_0}\}$. So, either $C_{v_2,\sP_0} \leq_{G^\star} C_{a,\sP_0}$ or  $C_{v_2,\sP_0} \leq_{G^\star} C_{b,\sP_0}$. This means that $v_2$ is an ancestor of $a$ or $b$. Hence, $a\padj b$.

Therefore, $E^{(1)}_{(\sP_0,\pi_0)}$ is contained in $E^{(1)}_{(\sP,\pi)}$ for every $(\sP,\pi) \in \sS$. So, for every $(\sP,\pi)\in S_1$, we have $E^{(1)}_{(\sP,\pi)} = E^{(1)}_{(\sP_0,\pi_0)}$, and thus, $(\sP_0,\pi_0)\in S_1$ and $E^{(1)}_{(\sP,\pi)}$ is equal to the set of $p$-adjacencies in $G^\star$.
\end{proof}
\bigskip

\begin{proof}[Proof of Proposition~\ref{prop: Rich2}]
Let $(\sP,\pi)\in S_1$, and $(a,b,c)$ be an unshielded conductor in $G^\star$. Then $b\in \an_{G^\star}(a,c)$. So, for any set $Z\subseteq [n]\setminus \{a,b,c\}$, we have $b\in \an_{G^\star}(a,c,Z)$. On the other hand, $a\padj b$ and $c\padj b$. So, $a,b,c$ is a path in $\moral\left(G^\star(\an_{G^\star} (a,c,Z)) \right)$ that does not contain any vertex in $Z$. Thus, $a \not\independent c \mid Z$ in $\moral\left(G^\star(\an_{G^\star} (a,c,Z)) \right)$, and a result, $a\not\independent c\mid Z$ in $G^\star$.
Now if $C_{b,\sP} \not\leq_\pi \max\{C_{a,\sP},C_{c,\sP}\}$, the set $Z\coloneqq \bigcup\Set{C\in \sP |  C\leq_\pi \max\{C_{a,\sP}, C_{c,\sP}\}}\setminus \{a,c\} $ does not contain $b$, and therefore, $a \not\independent c \mid Z$. So, $(a,c)\in E^{(1)}_{(\sP,\pi)}$, which means that $a \padj c$ since $(\sP,\pi)\in S_1$. This is a contradiction. Thus, $C_{b,\sP} \leq_\pi \max\{C_{a,\sP},C_{c,\sP}\}$, and $(a,b,c)\in E^{(2)}_{(\sP,\pi)}$. So, the set of unshielded conductors in $G^\star$ is contained in $E^{(2)}_{(\sP,\pi)}$ for all $(\sP,\pi)\in S_1$.

Let $(\sP_0,\pi_0)$ be the partially ordered partition associated with $G^\star$. Then for every $(a,b,c)\in E^{(2)}_{(\sP_0,\pi_0)}$, we have $a\padj b$, $c\padj b$, and $a\not\padj c$ since $(\sP_0,\pi_0)\in S_1$. Also, since $C_{b,\sP_0} \leq_{\pi_0} \max\{C_{a,\sP_0},C_{c,\sP_0}\}$, we have $C_{b,\sP_0}\leq_{G^\star} C_{a,\sP_0}$ or $C_{b,\sP_0}\leq_{G^\star} C_{c,\sP_0}$. So, $b$ is an ancestor of $a$ or $c$ in $G^\star$, and $(a,b,c)$ is unshielded conductor in $G^\star$. Hence, $E^{(2)}_{(\sP_0,\pi_0)}$ is the set of unshielded conductors in $G^\star$, and as a result, is contained in $E^{(2)}_{(\sP,\pi)}$ for every $(\sP,\pi) \in S_1$. So, for every $(\sP,\pi)\in S_2$, we have $E^{(2)}_{(\sP,\pi)} = E^{(2)}_{(\sP_0,\pi_0)}$, and thus, $(\sP_0,\pi_0)\in S_2$ and $E^{(2)}_{(\sP,\pi)}$ is equal to the set of unshielded conductors in $G^\star$.
\end{proof}
\bigskip

\begin{proof}[Proof of Proposition~\ref{prop: Rich3}]
Let $(\sP,\pi)\in S_2$, and $(a,b,c)$ be an unshielded perfect non-conductor in $G^\star$. Then $a$ and $c$ have a common child $d$ in $G^\star$ such that $d\in \an_{G^\star}(b)$. Consider that path $P : a,(a,d),d,(c,d),c$ in $G^\star$. The only collider on $P$ is $d$, whose descendant, vertex $b$, is in set $Z$, where 
\begin{align*}
    Z\coloneqq \bigcup \Set{C\in \sP |  C\leq_\pi \max\{ C_{a,\sP},C_{b,\sP} , C_{c,\sP}}\} \setminus \{a,c\}.
\end{align*}
So, $a \not\independent c \mid Z$, and as a result, $(a,b,c)\in E^{(3)}_{(\sP,\pi)}$. So, all unshielded perfect non-conductors of $G^\star$ are contained in $E^{(3)}_{(\sP,\pi)}$.

Now let $(\sP_0,\pi_0)$ be the partially ordered partition associated with $G^\star$, and $(a,b,c)$ be in $E^{(3)}_{(\sP_0,\pi_0)}$. We will prove that $(a,b,c)$ is an unshielded perfect non-conductor in $G^\star$. Define
\begin{align*}
    Z\coloneqq \bigcup \Set{C\in \sP_0 |  C\leq_{\pi_0} \max\{ C_{a,\sP_0},C_{b,\sP_0} , C_{c,\sP_0}}\} \setminus \{a,c\}.
\end{align*}
Since $(a,b,c)\in E^{(3)}_{(\sP_0,\pi_0)}$, we have $a\not\independent c\mid Z$. So, there exists a path $P:a=v_1,e_1,\ldots,e_{t-1},v_t=c$ in $G^\star$ such that all colliders on $P$ belong to $Z$ or have a descendant in $Z$, and for all $\ell\in\{2,\ldots,t-1\}$, if $v_\ell\in Z$, then $v_\ell$ is a collider on $P$. We claim that for all $\ell\in\{2,\ldots,t-1\}$, $v_\ell \in Z$. By contradiction, assume otherwise, and let $p\coloneqq \min \Set{\ell\in\{2,\ldots,t-1\} | v_\ell \not\in Z}$. If $e_{p-1}=(v_p,v_{p-1})$, then since $v_{p-1}\in Z\cup \{a\}$, by definition, $v_{p} \in Z$. So, $e_{p-1}= (v_{p-1},v_p)$. On the other hand, $v_p$ being a collider on $P$ means that $v_p \in Z$ or $v_p$ has a descendant in $Z$. Both cases result in $v_p \in Z$, and thus, are not possible. So, $e_{p}=(v_p,v_{p+1})$. Let
\begin{align*}
    q \coloneqq \max \Set{\ell \in \{p,\ldots,t-1\} | e_{\ell^\prime} = (v_{\ell^\prime},v_{\ell^\prime+1}) \text{ and $v_{\ell^\prime} \not \in Z$ for all $\ell^\prime \in \{p,\ldots,\ell\}$.} }.
\end{align*}
If $v_{q+1} \in Z\cup \{c\}$, then $v_q$  has a child in $Z\cup \{c\}$, and thus, $v_q\in Z$, which is impossible. So, $v_{q+1} \not \in Z\cup \{c\}$, and given the definition of $v_q$, $e_{q+1} = (v_{q+2},v_{q+1})$. Therefore, $v_{q+1}$ is a collider on $P$. But this means that $v_{q+1} \in Z$ or $v_{q+1}$ has a descendant in $Z$, both of which imply that $v_{q+1} \in Z$. This contradiction proves our claim. 

Now all the inner vertices of the path $P$ belonging to $Z$ implies that for all $\ell\in \{2,\ldots,t-1\}$, $v_\ell$ is a collider on the path $P$. So, either $t=2$, or $t=3$ and $v_2$ is a common child of $a$ and $c$ in $Z$. The former case does not happen because $a$ and $c$ are not $p$-adjacent, and thus, not adjacent. In the latter case, if $C_{v_2,G^\star} \leq_{G^\star} \max\{C_{a,G^\star},C_{c,G^\star} \}$, it means that $v_2$ is an ancestor of $a$ or $c$, and therefore, $a$ and $c$ are $p$-adjacent. So, since this is not the case and $v_2\in Z$, we have $C_{v_2,G^\star} \leq_{G^\star} C_{b,G^\star}$. This means that $b$ is a descendant of a common child of $a$ and $c$, and as a result, $(a,b,c)$ is an unshielded perfect non-conductor in $G^\star$.

Hence, for every $(\sP,\pi)\in S_2$, the set $E^{(3)}_{(\sP_0,\pi_0)}$ is contained in $E^{(3)}_{(\sP,\pi)}$. So, for every $(\sP,\pi)\in S_3$, we have $E^{(3)}_{(\sP,\pi)} = E^{(3)}_{(\sP_0,\pi_0)}$, and thus, $(\sP_0,\pi_0)\in S_3$ and $E^{(3)}_{(\sP,\pi)}$ is equal to the set of unshielded perfect non-conductors in $G^\star$.
\end{proof}
\bigskip

\begin{proof}[Proof of Proposition~\ref{prop: Rich4}]
    Let $(a,b_1,c),(a,b_2,c)\in [n]^3$ be two unshielded imperfect non-conductors, and $b_1$ be an ancestor of $b_2$ in $G^\star$. Then for any $(\sP,\pi)\in S_3$, by Proposition~\ref{prop: Rich1}, we have $(a,b_1),(c,b_1),(a,b_2),(c,b_2) \in E^{(1)}_{(\sP,\pi)}$ and $(a,c)\not\in E^{(1)}_{(\sP,\pi)}$, by Proposition~\ref{prop: Rich2}, $(a,b_1,c),(a,b_2,c) \not \in E^{(2)}_{(\sP,\pi)}$, and by Proposition~\ref{prop: Rich3}, $(a,b_1,c),(a,b_2,c) \not \in E^{(3)}_{(\sP,\pi)}$. We proceed to show $C_{b_1,\sP}\leq_\pi C_{b_2,\sP}$. Assume otherwise, and define
    \begin{align*}
       Z\coloneqq \bigcup \Set{C\in \sP |  C\leq_\pi \max\{ C_{a,\sP},C_{b_2,\sP} , C_{c,\sP}}\} \setminus \{a,c\}.
    \end{align*}
    Since $(a,b_1,c) \not \in  E^{(2)}_{(\sP,\pi)}$, $(a,b_1),(c,b_1)\in E^{(1)}_{(\sP,\pi)}$, and $(a,c)\not\in E^{(1)}_{(\sP,\pi)}$, one can conclude $C_{b_1,\sP}\not\leq_\pi \max\{C_{a,\sP},C_{c,\sP}\}$. So, if $C_{b_1,\sP} \not\leq_\pi C_{b_2,\sP}$, then $b_1\not \in Z$. On the other hand, $b_1\in \an_{G^\star}(b_2)$. So, $b_1\in \an_{G^\star}(a,c,Z)$. Considering that $a \padj b_1$ and $b_1 \padj c$, this shows that the graph $\moral\left(G^\star(\an_{G^\star} (a,c,Z)) \right)$ admits the path $P:a,b_1,c$, which does not contain any vertex in $Z$. So, $a \not\independent c \mid Z$ in $\moral\left(G^\star(\an_{G^\star} (a,c,Z)) \right)$, and thus, $a\not\independent c\mid Z$ in $G^\star$. However, this implies that $(a,b_2,c)\in E^{(3)}_{(\sP,\pi)}$, which is a contradiction. So, $C_{b_1,\pi}\leq_\pi C_{b_2,\pi}$, and as a result, $\left((a,b_1,c),(a,b_2,c)\right) \in E^{(4)}_{(\sP,\pi)}$.

    Now let $(\sP_0,\pi_0)$ be the partially ordered partition associated with $G^\star$, and $\left((a,b_1,c),(a,b_2,c)\right)\in E^{(4)}_{(\sP_0,\pi_0)}$.
    Then by Proposition~\ref{prop: Rich1}, $a\padj b_1$, $c\padj b_1$, $a\padj b_2$, $c\padj b_2$, and $a\not\padj c$. By Proposition~\ref{prop: Rich2}, $(a,b_1,c)$ and $(a,b_2,c)$ are unshielded non-conductors, and by Proposition~\ref{prop: Rich3}, $(a,b_1,c)$ and $(a,b_2,c)$ are imperfect. Since $C_{b_1,\pi_0} \leq_{\pi_0} C_{b_2,\pi_0}$, $b_1$ is an ancestor of $b_2$ in $G^\star$. Hence, $E^{(4)}_{(\sP_0,\pi_0)}$ is contained in $E^{(4)}_{(\sP,\pi)}$ for every $(\sP,\pi)\in S_3$. So, $(\sP_0,\pi_0)\in S_4$, and for every $(\sP,\pi)\in S_4$, $\left((a,b_1,c),(a,b_2,c)\right)\in E^{(4)}_{(\sP,\pi)}$ if and only if $(a,b_1,c),(a,b_2,c)$ are unshielded imperfect non-conductors and $b_1$ is an ancestor of $b_2$ in $G^\star$ as desired.
\end{proof}
\bigskip
\begin{proof}[Proof of Lemma~\ref{lemma: mutually exclusive}]
    Assume that case 1 does not happen. Then there exists $i\in [t]$ such that $C_{a_i,G^\star} \not \leq_{G^\star} C_{a_{i-1},G^\star}$.  Let 
    $
        k\coloneqq \min \Set{i\in [t] | C_{a_i,G^\star} \not\leq_{G^\star} C_{a_{i-1},G^\star} }.
    $
     If for all $i\in \{k,\ldots,t\}$, $C_{a_i,G^\star}\leq_{G^\star} C_{a_{i+1},G^\star}$, then $a_1,\ldots,a_t\in \an_{G^\star} (a_0,a_{t+1})$, which is not true, and thus, 
    \begin{align*}
        \ell \coloneqq \min \Set{i\in \{k,\ldots,t\} | C_{a_i,G^\star} \not\leq_{G^\star} C_{a_{i+1},G^\star}}
    \end{align*}
    is well defined. Let 
    \begin{align*}
        j\coloneqq \max \Set{i\in \{k,\ldots,\ell\} | C_{a_i,G^\star} \not\leq_{G^\star} C_{a_{i-1},G^\star}}. 
    \end{align*}
Then one can confirm that $(a_{j-1},a_j,a_{j+1})$ and $(a_{\ell-1},a_{\ell},a_{\ell+1})$ are mutually exclusive over the uncovered itinerary $(a_{j-1},\ldots,a_{\ell+1})$.
\end{proof}
\bigskip
\begin{proof}[Proof of Proposition~\ref{prop: Rich5}]
    We prove the proposition by induction on $t$.

    By Proposition~\ref{prop: Rich4}, the partially ordered partition associated with $G^\star$ is in $S_4=S^{(1)}_5$. Also, $(a_0,a_1,a_2)$ and $(a_0,a_1,a_2)$ are mutually exclusive with respect to the uncovered itinerary $(a_0,a_1,a_2)$ in $G^\star$ if and only if $a_0\padj a_1$, $a_1\padj a_2$, $a_0\not\padj a_2$, and $a_1$ is not an ancestor of $a_0$ or $a_2$ in $G^\star$. By Propositions~\ref{prop: Rich1} and~\ref{prop: Rich2}, for all $(\sP,\pi)\in S_4=S^{(1)}_5$, this is equivalent to $(a_0,a_1),(a_1,a_2)\in E^{(1)}_{(\sP,\pi)}$, $(a_0,a_2)\not\in E^{(1)}_{(\sP,\pi)}$, and $(a_0,a_1,a_2)\not\in E^{(2)}_{(\sP,\pi)}$. One can confirm that the aforementioned is equivalent to $(a_0,a_1,a_2)\in D^{(1)}_{(\sP,\pi)}$.

    Now assume the proposition holds for all $t\in [\tilde{t}-1]$. Suppose $(a_0,a_1,\ldots,a_{\tilde{t}+1})\in D^{(\tilde{t})}_{(\sP,\pi)}$ for some $(\sP,\pi)\in S^{(\tilde{t}-1)}_5$. We will show that $(a_0,a_1,a_2)$ and $(a_{\tilde{t}-1},a_{\tilde{t}},a_{\tilde{t}+1})$ are mutually exclusive with respect to the uncovered itinerary $(a_0,\ldots,a_{\tilde{t}+1})$. Assume otherwise. Then by Lemma~\ref{lemma: mutually exclusive}, one of the following happens:
    \begin{enumerate}
        \item Each of the vertices $a_1,\ldots,a_{\tilde{t}}$ is the ancestor of $a_0$ or $a_{\tilde{t}+1}$, or
        \item there exists $j,\ell\in [\tilde{t}]$ with $j\leq \ell$ such that $(a_{j-1},a_j,a_{j+1})$ and $(a_{\ell-1},a_\ell,a_{\ell+1})$ are mutually exclusive with respect to the uncovered itinerary $(a_{j-1},\ldots,a_{\ell+1})$ in $G^\star$. 
    \end{enumerate}
    If case 1 happens, then $P: a_0,a_1,\ldots,a_{\tilde{t}+1}$ is a path in $\moral\left(G^\star(\an_{G^\star} (a_0,a_{\tilde{t}+1},Z)) \right)$ for all $Z\subseteq [n]\setminus \{a_0,a_{\tilde{t}+1}\}$. So, if $a_i\not\in Z$ for all $i\in [\tilde{t}]$, we get $a_0 \not\independent a_{\tilde{t}+1} \mid Z$ in $\moral\left(G^\star(\an_{G^\star} (a_0,a_{\tilde{t}+1},Z)) \right)$, which means that $a_0 \not\independent a_{\tilde{t}+1} \mid Z$ in $G^\star$. On the other hand, since $(a_{\tilde{t}+1},a_0)\not\in E^{(1)}_{(\sP,\pi)}$, we have 
    \begin{align*}
        a_0 \independent a_{\tilde{t}+1} \mid \bigcup \Set{C\in \sP |  C\leq_\pi \max\{C_{a_0,\sP}, C_{a_{\tilde{t}+1},\sP}\}}\setminus \{a_0,a_{\tilde{t}+1}\}.
    \end{align*}
    Hence, there exists $i\in [\tilde{t}]$ such that $a_i \in \bigcup \Set{C\in \sP |  C\leq_\pi \max\{C_{a_0,\sP}, C_{a_{\tilde{t}+1},\sP}\}}\setminus \{a_0,a_{\tilde{t}+1}\}$. But this implies that $C_{a_i,\sP} \leq_{\pi} \max \{C_{a_0,\sP}, C_{a_{\tilde{t}+1},\sP}\}$, which is a contradiction. So, case 1 is impossible to happen.

    Assume that case 2 happens. Then since $(a_{j-1},a_j,a_{j+1})$ and $(a_{\ell-1},a_\ell,a_\ell)$ are mutually exclusive with respect to the uncovered itinerary $(a_{j-1},\ldots,a_{\ell+1})$, either $j-1\neq 0$ or $\ell+1\neq \tilde{t}+1$. Thus, by the induction hypothesis, $(a_{j-1},\ldots,a_{\ell+1})\in D^{(t)}_{(\sP,\pi)}$ for some $t\in [\tilde{t}-1]$. Hence, given the definition of $D^{(\tilde{t})}_{(\sP,\pi)}$, we get $(a_0,\ldots,a_{\tilde{t}+1})\not \in D^{(\tilde{t})}_{(\sP,\pi)}$. This contradiction proves that $(a_0,a_1,a_2)$ and $(a_{\tilde{t}-1},a_{\tilde{t}},a_{\tilde{t}+1})$ are mutually exclusive with respect to the uncovered itinerary $(a_0,\ldots,a_{\tilde{t}+1})$.

    Now let $(\sP_0,\pi_0)$ be the partially ordered partition associated with $G^\star$. If $(a_0,a_1,a_2)$ and $(a_{\tilde{t}-1},a_{\tilde{t}},a_{\tilde{t}+1})$ are mutually exclusive with respect to the uncovered itinerary $(a_0,a_1,\ldots,a_{\tilde{t}+1})$ in $G^\star$, then by Proposition~\ref{prop: Rich1}, $(a_i,a_{i+1}) \in E^{(1)}_{(\sP_0,\pi_0)}$ for all $i \in \{0,\ldots,\tilde{t}\}$, and $(a_i,a_j) \not \in E^{(1)}_{(\sP,\pi)}$ for all $i \in \{2,\ldots,\tilde{t}+1\}$ and $j \in \{0,\ldots,i-2\}$. Also, $C_{a_1,G^\star}=\cdots=C_{a_{\tilde{t}},G^\star}$, $C_{a_1,G^\star}\not\leq_{G^\star} C_{a_0,G^\star}$, and $C_{a_1,G^\star}\not\leq_{G^\star} C_{a_{\tilde{t}+1},G^\star}$. So, $(a_0,\ldots,a_{\tilde{t}+1})\in D^{(\tilde{t})}_{(\sP_0,\pi_0)}$. This proves that $D^{(\tilde{t})}_{(\sP,\pi)} \subseteq D^{(\tilde{t})}_{(\sP_0,\pi_0)}$ for all $(\sP,\pi)\in S^{(\tilde{t}-1)}_5$. 

    Therefore, for all $(\sP,\pi)\in S^{(\tilde{t})}_5$, we get $D^{(\tilde{t})}_{(\sP,\pi)} = D^{(\tilde{t})}_{(\sP_0,\pi_0)}$. So, the proposition also holds for $t=\tilde{t}$. This concludes our proof.
\end{proof}
\bigskip

\begin{proof}[Proof of Proposition~\ref{prop: Rich6}]
Let $(a_0,a_1,a_2)$ and $(a_{t-1},a_t,a_{t+1})$ be mutually exclusive with respect to the uncovered itinerary $(a_0,a_1,\ldots,a_{t+1})$, $(a_0,b_1,a_{t+1})$ be an unshielded imperfect non-conductor, and $a_1$ be an ancestor of $b$ in $G^\star$. Then by Propositions~\ref{prop: Rich1},~\ref{prop: Rich2},~\ref{prop: Rich3}, and~\ref{prop: Rich5}, for any $(\sP,\pi)\in S^{(n-2)}_{5}$, $(a_0,b),(a_{t+1},b) \in E^{(1)}_{(\sP,\pi)}$, $(a_0,b,a_{t+1}) \not \in E^{(2)}_{(\sP,\pi)}$, $(a_0,b,a_{t+1})\not \in E^{(3)}_{(\sP,\pi)}$, and $(a_0,a_1,\ldots,a_{t},a_{t+1})\in D^{(t)}_{(\sP,\pi)}$. We will show that $C_{a_1,\sP} \leq_{\pi} C_{b,\sP}$. Define
\begin{align*}
       Z\coloneqq \bigcup \Set{C\in \sP |  C\leq_\pi \max\{ C_{a_0,\sP},C_{a_{t+1},\sP} , C_{b,\sP}}\} \setminus \{a_0,a_{t+1}\}.
\end{align*}
Since $C_{a_1,G^\star}= C_{a_2,G^\star}=\cdots=C_{a_t,G^\star}$ and $a_1$ is an ancestor of $b$ in $G^\star$, we have $a_1,\ldots,a_t\in \an_{G^\star}(a_0,a_{t+1},Z)$. So, considering that $a_i\padj a_{i+1}$ for all $i\in \{0,\ldots,t\}$, 
    $P:a_0,a_1,\ldots,a_t,a_{t+1}$
is a path in $\moral\left(G^\star(\an_{G^\star} (a_0,a_{t+1},Z)) \right)$. If $a_i\not\in Z$ for all $i\in [t]$, then $P$ is a path with no vertex in $Z$, which means that $a_0\not\independent a_{t+1} \mid Z$ in $\moral\left(G^\star(\an_{G^\star} (a_0,a_{t+1},Z)) \right)$, and therefore, $a_0\not\independent a_{t+1}\mid Z$ in $G^\star$. However, since $(a_0,b,a_{t+1})\not\in E^{(3)}_{(\sP,\pi)}$, this statement cannot be true. So, there exists $i\in [t]$ such that $a_i\in Z$, which implies that $C_{a_i,\sP} \leq_\pi \max\{ C_{a_0,\sP},C_{a_{t+1},\sP} , C_{b,\sP}\}$. On the other hand, since $(a_0,\ldots,a_{t+1})\in D^{(t)}_{(\sP,\pi)}$, $C_{a_1 ,\sP}=C_{a_i,\sP}$ and $C_{a_1,\sP} \not\leq_\pi \max\{C_{a_0,\sP},C_{a_{t+1},\sP}\}$. So, $C_{a_1,\sP}\leq_\pi C_{b,\sP}$ as desired. This means that $\left((a_0,a_1,\ldots,a_{t},a_{t+1}),(a_0,b,a_{t+1})\right)\in E^{(6)}_{(\sP,\pi)}$.

Now let $(\sP_0,\pi_0)$ be the partially ordered partition associated with $G^\star$, and suppose that $\left((a_0,a_1,\ldots,a_{t},a_{t+1}),(a_0,b,a_{t+1})\right)\in E^{(6)}_{(\sP_0,\pi_0)}$. Then Propositions~\ref{prop: Rich1},~\ref{prop: Rich2},~\ref{prop: Rich3}, and~\ref{prop: Rich5} immediately imply that $(a_0,a_1,a_2)$ and $(a_{t-1},a_t,a_{t+1})$ are mutually exclusive with respect to the uncovered itinerary $(a_0,a_1,\ldots,a_{t+1})$, $(a_0,b_1,a_{t+1})$ is an unshielded imperfect non-conductor, and $a_1$ is an ancestor of $b$ in $G^\star$. This shows that $E^{(6)}_{(\sP_0,\pi_0)}\subseteq E^{(6)}_{(\sP,\pi)}$ for all $(\sP,\pi)\in S^{(n-2)}_{5}$, and thus, $E^{(6)}_{(\sP_0,\pi_0)} = E^{(6)}_{(\sP,\pi)}$ for all $(\sP,\pi)\in S_6$. This concludes the proof.
\end{proof}
\section{Pseudocode of SCCR algorithm 1: Construct and correct} \label{app: SCCR2-alg}
This section presents the details of our first SCCR algorithm.
\floatstyle{ruled}
\restylefloat{algorithm}
\begin{algorithm}[H]
\caption{Subalgorithm of Algorithm~\ref{alg: SCCR2}}\label{alg: subproc of SCCR2}
\begin{algorithmic}[1]
\Require{An integer $j\in \mathbb{N}$, sets $A_C,B_C$ and $ComCh_C$, and lists $D$, $D_{copy}$ and $Cause$.}
                \State Set $avoid\coloneqq D[j]$.
                \If {$D_{copy}[j]=2$} 
                    \If {$j$ is odd} 
                        Set $T\coloneqq \{1,\ldots,j-1\} \cup \{j+2,\ldots,2\left|ComCh_C \right|\}$.
                    \Else 
                        \ Set $T\coloneqq \{1,\ldots,j-2\} \cup \{j+1,\ldots,2\left| ComCh_C\right|\}$.
                    \EndIf
                    \State {After setting the following, go to step~\ref{alg: SCCR2 step 1} of Algorithm~\ref{alg: SCCR2}:
                        \begin{align*}
                            &A_C \coloneqq \text{a shuffling of $A_C$}, && B_C \coloneqq \text{a shuffling of $B_C$}, && D \coloneqq D[T], \\
                            &D_{copy} \coloneqq D_{copy}[T], &&Cause\coloneqq Cause[T], && a\coloneqq 0.
                        \end{align*}
                        }
                \EndIf
                \State Set $j\coloneqq Cause[j]$.
                \State Set $D_{new}\coloneqq D[1,\ldots,j-1]$, $D_{copy,new} \coloneqq D_{copy}[1,\ldots,j-1]$, and 
                \NoNumber{$Cause_{new}\coloneqq Cause[1,\ldots,j-1]$.}
                \For {$k$ in $1:j-1$}
                    \If {$D_{copy,new}[k]\in \{-1,0\}$ and $Cause[k]\geq j$} 
                        \State Set $D_{new}[k] \coloneqq (D_{new}[k]_2,D_{new}[k]_1)$. \Comment{Flip $D_{new}[k]$.}
                        \State Set $D_{copy,new}[k] \coloneqq D_{new}[k]$, and $Cause[k] \coloneqq k$.
                    \EndIf
                \EndFor
                \If {$D_{copy}[j]\in \{-1,0\}$ and $D[j]\in A_C$} 
                    Set $a\coloneqq D[j]_2$.
                \ElsIf {$D_{copy}[j]\in \{-1,0\}$}
                    Set $a\coloneqq D[j]_1$.
                \ElsIf {$D[j]\in A_C$} 
                    Set $a\coloneqq D[j]_1$.
                \Else 
                    \ Set $a\coloneqq D[j]_2$.
                \EndIf
                \State {After setting the following, go to step~\ref{alg: SCCR2 step 1} of Algorithm~\ref{alg: SCCR2}:
                \vspace{-3mm}
                        \begin{align*}
                            &A_C \coloneqq \text{a shuffling of $A_C$}, &&B_C \coloneqq \text{a shuffling of $B_C$}, &&D \coloneqq D_{new},\\
                            &D_{copy} \coloneqq D_{copy,new}, &&Cause\coloneqq Cause_{new}.
                \end{align*}}
\end{algorithmic}
\end{algorithm}
\floatstyle{nobottomruled}
\restylefloat{algorithm}
\begin{algorithm}[H]
\caption{SCCR algorithm 1: Construct and correct}\label{alg: SCCR2}
\begin{algorithmic}[1]
\Require{Sets $C\subseteq[n]$, $A_C\subseteq C^2$, $B_C \subseteq ([n]\setminus C) \times C$, $ComCh_C\subseteq ([n]\setminus C)^2$ and $NoComCh_C\subseteq ([n]\setminus C)^2$, a positive integer $N$.}
\Ensure{Either a declaration of failure, or a set $E_C \subseteq [n]^2$ with the properties described in Definition~\ref{def: strongly connected component recovery algorithm}.}
\State Set $I\coloneqq 0$.
\State Consider $D$, $D_{copy}$ and $Cause$ to be empty lists, and set $a\coloneqq 0$ and $avoid\coloneqq 0$.
\State Set $D_{initial} \coloneqq D$, $D_{copy,initial}\coloneqq D_{copy}$, $Cause_{initial}\coloneqq Cause$, $avoid_{initial}\coloneqq avoid$, and $a_{initial}\coloneqq a$. \label{alg: SCCR2 step 1}
\If {$I=N$} 
        \Return "failed".
\EndIf
\algstore{SCCR2-part1}
\end{algorithmic}
\end{algorithm}
\floatstyle{nobottomandcaptionruled}
\restylefloat{algorithm}
\begin{algorithm}[H]
\begin{algorithmic} [1]
\algrestore{SCCR2-part1}
\State Set $I\coloneqq I+1$. 
\State Set $first.edge \coloneqq \TRUE$.
\LineComment{Define $Checked$ and $Parents$ based on the input.}
\State Consider $Check$ to be an empty list, and $Parents:C\to \mathcal{A}$ to be a function with $Parents(v)=\emptyset$ for all $v\in C$, where $\mathcal{A}$ is the family of multisubsets of $[n]$.

\For{any edge $e=(e_1,e_2)$ in $D$}
    \State Add $\left(\min\{e_1,e_2\},\max\{e_1,e_2\}\right)$ to $Checked$ if it is not already in $Checked$.
    \LineComment{The edges in $D$ that have been later removed have code 0 in $D_{copy}$.}
    \If {$e$ is the $i$th element in $D$ and $D_{copy}[i]\neq 0$}
        \State Set $Parents(e_2)\coloneqq Parents(e_2) \cup \{e_1\}$.
    \EndIf
\EndFor
\LineComment{Choosing common children in $C$ for the pairs in $ComCh_C$.} 
\LineComment{The edges which are added to $D$ to make sure the pairs in $ComCh_C$ have a common child in $C$ have code 2 in $D_{copy}$.}
\If {the number of 2's in $D_{copy}$ is less that $2\left|ComCh_C\right|$}
    \For {pairs $p=(p_1,p_2)$ in $ComCh_C$} 
        \If {a common child hasn't yet been specified for $p$ in $C$} 
            \State Find $v\in C$ such that $(p_1,v),(p_2,v)\in B_C$, $avoid\not\in \{(p_1,v),(p_2,v)\}$, and 
            \NoNumber{the edges $(p_1,v)$ and $(p_2,v)$ are safe with respect to $Parents$.}
            \If {no such $v$ is found} 
                \State {\Return "failed".}
            \Else 
                \State Set $D\coloneqq D,(p_1,v),(p_2,v)$ and $D_{copy}\coloneqq D_{copy},2,2$.
                \State Set $Cause\coloneqq Cause,p_2,p_1$ and $Parents(v)\coloneqq Parents(v) \cup \{p_1,p_2\}$.
                \State Add to $Checked$ any of the pairs $\left(\min\{p_1,v\},\max\{p_1,v\}\right)$ and 
                \NoNumber{$\left(\min\{p_2,v\},\max\{p_2,v\}\right)$ which is not already in $Checked$.} 
            \EndIf
        \EndIf
    \EndFor
\EndIf

\If {$a=0$} 
    Choose $a\in C$ to be a minimizer for $$\left| \Set{v\in C | (a,v),(v,a)\in A_C}\right| + \left|\Set{v\in [n]\setminus C | (v,a)\in B_C}\right|.$$
\EndIf
\State Let $f:C\to \mathbb{Z}$ be a function with $f(v)=0$ for all $v\in V$. 
\State Set $R\coloneqq 0$. \Comment{$R$ controls the number of consecutive iterations of the loop of line~\ref{alg: SCCR2 loop} in which no new edge is added to the construction.}
\State Set $R_{compare} \coloneqq |D|$. \Comment{$R_{compare}$ stores the maximum number of such consecutive iterations allowed.}
\While {$\left|Checked\right| \neq \tfrac{\left|A_C\right|}{2}+\left|B_C\right|$ or $\big(C=\{v_1,v_2\}$ for some $v_1,v_2\in [n]$ and $A_C\neq \emptyset$ and $\{(v_1,v_2),(v_2,v_1)\} \not\subseteq \Set{e\in [n]^2 | e=D[i],D_{copy}[i]\neq 0\text{ for some } i\in \mathbb{N}}\big)$} \label{alg: SCCR2 loop}
    \If {$R>R_{compare}$}
        After setting the following, go back to step~\ref{alg: SCCR2 step 1}:
            \begin{align*}
                            &A_C \coloneqq \text{a shuffling of $A_C$}, && B_C \coloneqq \text{a shuffling of $B_C$}, && D \coloneqq D_{initial}, \\
                            &D_{copy} \coloneqq D_{copy,initial}, &&Cause\coloneqq Cause_{initial}, && a\coloneqq a_{initial}\\
                            &avoid\coloneqq avoid_{initial}.
            \end{align*}
    \EndIf
    \LineComment{Add edges to the construction in the following manner.}
    \State Let $S$ be an ordered set containing all vertices $b\in [n]$ with the following properties: 
    \begin{itemize}
        \item $(a,b) \in A_C$ or $(b,a)\in B_C$ 
        \item if $first.edge=\TRUE$ and $b\in C$, then $avoid\neq (a,b)$, and
        \item if $first.edge=\TRUE$ and $b\in [n]\setminus C$, then $avoid\neq (b,a)$.
    \end{itemize}
\algstore{SCCR2-part2}
\end{algorithmic}
\end{algorithm}
\floatstyle{nobottomandcaptionruled}
\restylefloat{algorithm}
\begin{algorithm}[H]
\begin{algorithmic} [1]
\algrestore{SCCR2-part2}
    \If {$S = \emptyset$} 
        \State \Return "failed".
    \Else
        \State Choose $b\in S$ such that $\left(\min\{a,b\},\max\{a,b\}\right)\not\in Checked$.
        \If {no such $b$ exists} 
            \State Set $f(a)\coloneqq f(a)+1$ and choose $b$ to be element number $f(a)$ (modulus $|S|$) in  $S$.
        \EndIf
    \EndIf
    \State $Cause.num \coloneqq |D|+1$. \Comment{This is the code assigned to any edge added, flipped, or removed}
    \NoNumber{in this iteration of the loop. This code is stored in $Cause$.}
    \If {$a,b\in C$} 
        \State Set $c\coloneqq a$ and $d\coloneqq b$.
    \Else 
        \State Set $c\coloneqq b$ and $d\coloneqq a$.
    \EndIf
    \State Set $D\coloneqq D,(c,d)$ and $D_{copy}\coloneqq D_{copy},(c,d)$.
    \State Set $Cause \coloneqq Cause,Cause.num$ and $Parents(d)\coloneqq Parents(d)\cup \{c\}$. 
    \If {$\left(\min\{c,d\},\max\{c,d\}\right) \in Checked$} 
        \State Set $R\coloneqq R+1$.
    \Else
        \State Set $R\coloneqq 0$ and $R_{compare} = |D|$, and add $\left(\min\{c,d\},\max\{c,d\}\right)$ to $Checked$.
    \EndIf
    \If {edge $(c,d)$ is safe with respect to $Parents$}
        \State Set $a\coloneqq d$ and go to step~\ref{alg: SCCR2 next edge}.
    \EndIf
    \LineComment{Remove or flip the edges in $D$.}
    \State Set $i\coloneqq |D|$, and let $Potential.Problems$ be an empty list.
    \State Set $R_{vertex}\coloneqq 0$ and $R_{vertex,compare}\coloneqq \infty$. \Comment{$R_{vertex}$ counts the number of consecutive}
    \NoNumber{iterations of the loop of line~\ref{alg: SCCR2 loop2} with edges incident to the same vertex, and $R_{vertex,compare}$}
    \NoNumber{sets the maximum number of such iterations allowed.}
    \While{\TRUE} \label{alg: SCCR2 loop2}
        \If {$R_{vertex}>R_{vertex,compare}$} 
            \State \Return "failed".
        \EndIf
        \State Set $e \coloneqq D[i]$.
        \If {$e_1\in C$} 
            \State Find $v\in C$ such that $(e_1,v),(e_2,v)\in A_C$, $(e_1,v),(e_2,v)$ are safe with respect to 
            \NoNumber{$Parents$, and if $first.edge=\TRUE$, then $avoid\not\in \{(e_1,v),(e_2,v)\}$.}
            \If {such $v$ is found} 
                \State Set $D[i]\coloneqq (e_2,e_1)$, $D_{copy}[i]\coloneqq 0$, and $Cause[i]\coloneqq Cause.num$. 
                \LineComment{When an edge is deleted, it is flipped in $D$ and assigned code 0 in $D_{copy}$.}
                \State Set $D\coloneqq D,(e_1,v),(e_2,v)$, $D_{copy} \coloneqq D_{copy},1,1$, and 
                \NoNumber{$Cause\coloneqq Cause.num,Cause.num$. }
                \LineComment{Edges added to preserve a $p$-adjacency are assigned code 1 in $D_{copy}$.}
                \State Set $Parents(e_2) \coloneqq Parents(e_2) \setminus \{e_1\}$ and $Parents(v) \coloneqq Parents(v) \cup \{e_1,e_2\}$.
                \If {$\left(\min\{e_1,v\},\max\{e_1,v\} \right) \not \in Checked$} 
                    \State {Set $R \coloneqq 0$ and $Checked \coloneqq Checked, \left(\min\{e_1,v\},\max\{e_1,v\} \right)$.}
                \EndIf
                \If {$\left(\min\{e_2,v\},\max\{e_2,v\} \right) \not \in Checked$} 
                    \State {Set $R \coloneqq 0$ and $Checked \coloneqq Checked, \left(\min\{e_2,v\},\max\{e_2,v\} \right)$.}
                \EndIf
                \State Set $a\coloneqq v$ and go to step~\ref{alg: SCCR2 end of loop}.
            \EndIf
            \LineComment{Since $e$ couldn't be removed, we flip it. A flipped edge is assigned code -1 in $D_{copy}$.}
            \State Set $D[i]\coloneqq (e_2,e_1)$, $D_{copy}[i] \coloneqq -1$, and $Cause[i] \coloneqq Cause.num$.
            \State Set $Parents(e_1)\coloneqq Parents(e_1)\cup \{e_2\}$ and $Parents(e_2) \coloneqq Parents(e_2) \setminus \{e_1\}$.
\algstore{SCCR2-part3}
\end{algorithmic}
\end{algorithm}
\floatstyle{nobottomandcaptionruled}
\restylefloat{algorithm}
\begin{algorithm}[H]
\begin{algorithmic} [1]
\algrestore{SCCR2-part3}
            \If {$(e_2,e_1)$ is safe with respect to $Parents$}
                \State Set $a\coloneqq e_1$ and go to step~\ref{alg: SCCR2 end of loop}.
            \EndIf
            \State Set $S \coloneqq \big\{k\in \mathbb{N} \mid D[k] = (v,e_1) \text{ such that } v\in Parents(e_1) 
                \text{ and $D[k]$ is}$
            \NoNumber{incompatible with $(e_2,e_1)$.\big\}.}
            \If {there exists $k\in S$ with $D_{copy}[k] \in \{-1,1\}$} 
                \State Choose $j\in S$ such that $D[j]$ is the first appearance of  
                \NoNumber{$D\left[\max(S\cap \Set{k\in \mathbb{N}|D_{copy}[k]\in \{-1,1\}})\right]$ in $D$.}
                \State Run Algorithm~\ref{alg: subproc of SCCR2} with $j,A_C,B_C,ComCh_C,D_{copy}$ and $Cause$.
            \EndIf
            \State Set $Potential.Problems \coloneqq Potential.Problems,(e_2,e_1)$. 
            \State Set $i \coloneqq \max(S\cap \Set{k \in \mathbb{N} |  D_{copy}[k] \neq 0})$.
            \State Set $R_{vertex} \coloneqq 0$.
            \State Go back to step~\ref{alg: SCCR2 loop2}.
        \Else 
            \State Find $v\in C$ such that $(e_1,v)\in B_C$ and $(e_2,v)\in A_C$, $(e_1,v),(e_2,v)$ are safe with 
            \NoNumber{respect to $Parents$, and if $first.edge=\TRUE$, then $avoid\not\in \{(e_1,v),(e_2,v)\}$.}
            \If {such $v$ is found} 
                \If {$D_{copy}[i] \neq 2$}
                    \State Set $D[i]\coloneqq (e_2,e_1)$, $D_{copy}[i]\coloneqq 0$, and $Cause[i]\coloneqq Cause.num$. 
                \EndIf
                \State Set $D\coloneqq D,(e_1,v),(e_2,v)$, $D_{copy} \coloneqq D_{copy},1,1$, and 
                \NoNumber{$Cause\coloneqq Cause.num,Cause.num$. }
                \State Set $Parents(e_2) \coloneqq Parents(e_2) \setminus \{e_1\}$ and $Parents(v) \coloneqq Parents(v) \cup \{e_1,e_2\}$.
                \If {$\left(\min\{e_1,v\},\max\{e_1,v\} \right) \not \in Checked$} 
                    \State Set $R \coloneqq 0$ and $R_{compare}\coloneqq |D|$.
                    \State Set $Checked \coloneqq Checked, \left(\min\{e_1,v\},\max\{e_1,v\} \right)$.
                \EndIf
                \If {$\left(\min\{e_2,v\},\max\{e_2,v\} \right) \not \in Checked$} 
                    \State Set $R \coloneqq 0$ and $R_{compare}\coloneqq |D|$.
                    \State Set $Checked \coloneqq Checked, \left(\min\{e_2,v\},\max\{e_2,v\} \right)$.
                \EndIf
                \State Set $S\coloneqq \Set{k\in \mathbb{N} | D[k]=e \text{ and } D_{copy}[k]=2}$.
                \If {$S=\emptyset$}
                    \State Set $a\coloneqq v$ and go to step~\ref{alg: SCCR2 end of loop}.
                \Else 
                    \For {$k \in S$} 
                         \State Find $w\in C$ such that $(e_1,w),(Cause[k],w)\in B_C$ and $(e_1,w)$ and 
                         \NoNumber{$(Cause[k],w)$  are safe with respect to $Parents$.}
                         \If {such $w$ is found}
                            \State Set $D[k]=(e_1,w)$.
                            \State If $k$ is odd, set $D[k+1]\coloneqq (Cause[k],w)$. Otherwise, set $D[k-1]\coloneqq$
                            \NoNumber{$(Cause[k],w).$}
                            \State Set $Parents(w) \coloneqq Parents(w) \cup \{e_1, Cause[k]\}$.
                            \State Set $Parents(e_2) \coloneqq Parents(e_2) \setminus \{e_1, Cause[k]\}$.
                            \If {$\left(\min\{e_1,w\},\max\{e_1,w\} \right) \not \in Checked$} 
                                \State Set $R \coloneqq 0$ and $R_{compare}\coloneqq |D|$.
                                \State Set $Checked \coloneqq Checked, \left(\min\{e_1,w\},\max\{e_1,w\} \right)$.
                            \EndIf
                            \If {$\left(\min\{Cause[k],w\},\max\{Cause[k],w\} \right) \not \in Checked$} 
                                \State Set $R \coloneqq 0$ and $R_{compare}\coloneqq |D|$.
                                \State Set $Checked \coloneqq Checked, \left(\min\{Cause[k],w\},\max\{Cause[k],w\} \right)$.
                            \EndIf
                        \Else 
                            \State Go to step~\ref{alg: SCCR2 couldn't delete the incoming edge.}.
                        \EndIf
                    \EndFor
\algstore{SCCR2-part4}
\end{algorithmic}
\end{algorithm}
\floatstyle{nocaptionruled}
\restylefloat{algorithm}
\begin{algorithm}[H]
\begin{algorithmic} [1]
\algrestore{SCCR2-part4}
                    \State Set $a\coloneqq v$ and go to step~\ref{alg: SCCR2 end of loop}.
                \EndIf
            \EndIf
            \State Set $S \coloneqq \big\{k\in \mathbb{N} \mid D[k] = (v,e_2) \text{ such that } v\in Parents(e_2) 
                \text{ and $D[k]$ is}$ \label{alg: SCCR2 couldn't delete the incoming edge.}
            \NoNumber{incompatible with $e$.\big\}.}
            \If {there exists $k\in S$ with $D_{copy}[k] \in \{-1,1\}$} 
                \State Choose $j$ randomly from $S\cap \Set{k\in \mathbb{N} | D_{copy}[k]\in \{-1,1\}}$.
                \State Run Algorithm~\ref{alg: subproc of SCCR2} with $j,A_C,B_C,ComCh_C,D_{copy}$ and $Cause$.
            \EndIf
            \State Set $Potential.Problems \coloneqq Potential.Problems,e$. 
            \State Set $i \coloneqq \max(S\cap \Set{k \in \mathbb{N} |  D_{copy}[k] \neq 0})$.
            \If {$R_{vertex}=0$}
                \State Set $R_{vertex,compare} \coloneqq \left| \Set{\tilde{e}\in D | \tilde{e}_2=e_2} \right|$.
            \EndIf
            \State Set $R_{vertex} \coloneqq R_{vertex}+1$.
            \State Go back to step~\ref{alg: SCCR2 loop2}.
        \EndIf
    \EndWhile
    \For {$e\in Potential.Problems$} 
        \If {$e$ is not safe with respect to $Parents$}
            \State Set $S \coloneqq \Set{v\in Parents(e_2) | (v,e_2) \text{ is incompatible with } e}$. \label{alg: SCCR2 end of loop}
            \State Choose $v$ randomly from $S$ and let $D[j]$ be the first appearance of $(v,e_2)$ in $D$ with 
            \NoNumber{$D_{copy}[j]\neq 0$.}
            \State Run Algorithm~\ref{alg: subproc of SCCR2} with $j,A_C,B_C,ComCh_C,D_{copy}$ and $Cause$.
        \EndIf
    \EndFor
    \State Set $first.edge \coloneqq \FALSE$. \label{alg: SCCR2 next edge}
\EndWhile
\State Set $E_C = \Set{D[i]|D_{copy}[i]\neq 0}$.
\If {$C$ is a strongly connected component in the graph $([n],E_C)$}
    \State \Return $E_C$.
\Else 
    \State \Return "failed".
\EndIf
\end{algorithmic}
\end{algorithm}
\section{An example of the execution of Algorithm~\ref{alg: SCCR2}} \label{app: SCCR2}
\begin{figure}[H]
     \centering
     \begin{subfigure}[b]{0.27\textwidth}
         \centering
         \includegraphics[width=\textwidth]{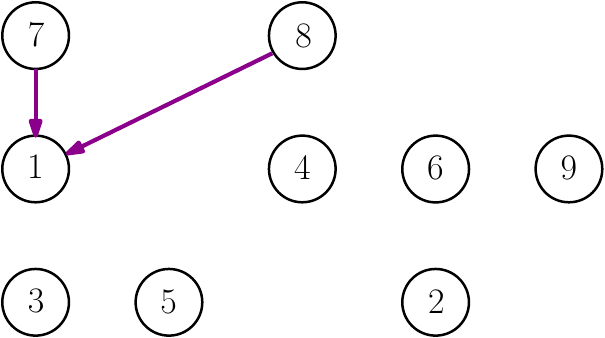}
     \end{subfigure}
     \hfill
     \begin{subfigure}[b]{0.27\textwidth}
         \centering
         \includegraphics[width=\textwidth]{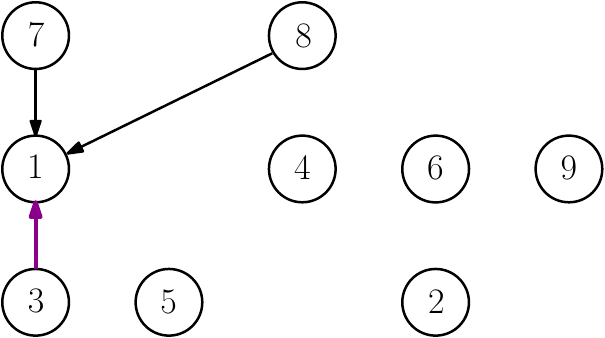}
     \end{subfigure}
     \hfill
     \begin{subfigure}[b]{0.27\textwidth}
         \centering
         \includegraphics[width=\textwidth]{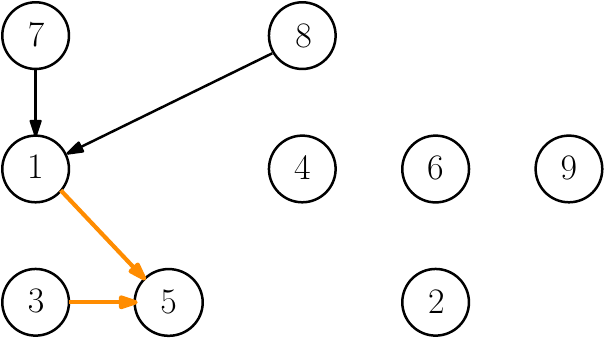}
     \end{subfigure}
     \hfill
     \vspace{7mm}

     \begin{subfigure}[b]{0.27\textwidth}
         \centering
         \includegraphics[width=\textwidth]{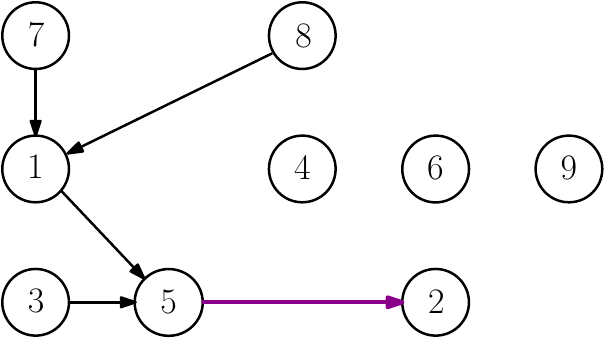}
     \end{subfigure}
     \hfill
     \begin{subfigure}[b]{0.27\textwidth}
         \centering
         \includegraphics[width=\textwidth]{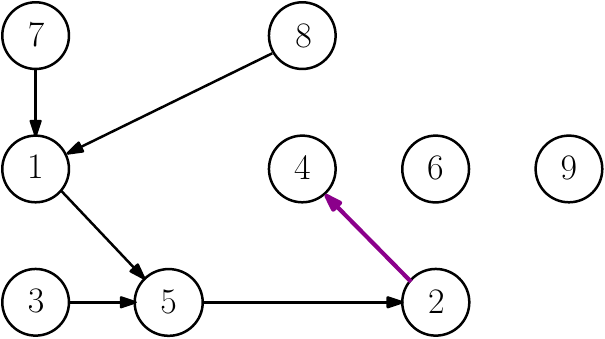}
     \end{subfigure}
     \hfill
     \begin{subfigure}[b]{0.27\textwidth}
         \centering
         \includegraphics[width=\textwidth]{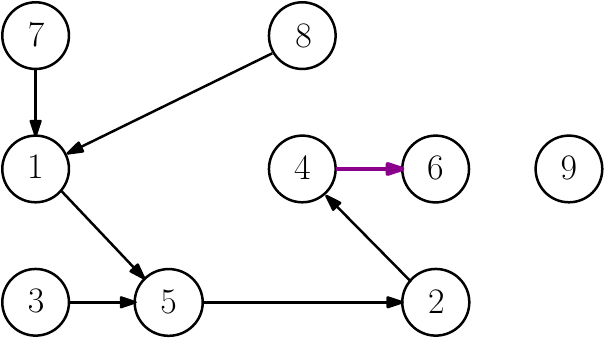}
     \end{subfigure}
     \hfill
     \vspace{7mm}

     \begin{subfigure}[b]{0.27\textwidth}
         \centering
         \includegraphics[width=\textwidth]{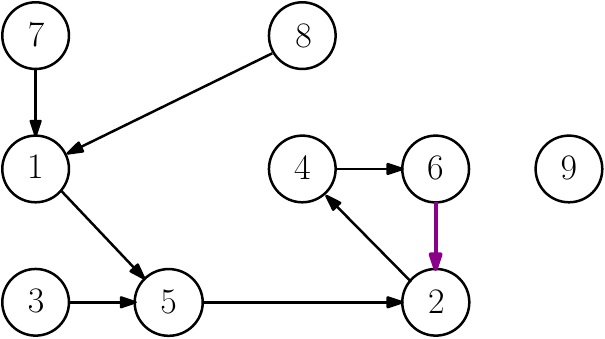}
     \end{subfigure}
     \hfill
     \begin{subfigure}[b]{0.27\textwidth}
         \centering
         \includegraphics[width=\textwidth]{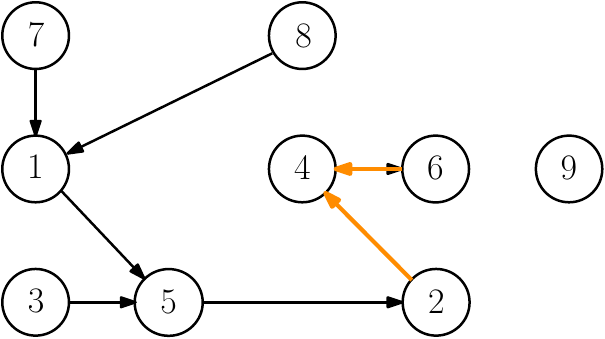}
     \end{subfigure}
     \hfill
     \begin{subfigure}[b]{0.27\textwidth}
         \centering
         \includegraphics[width=\textwidth]{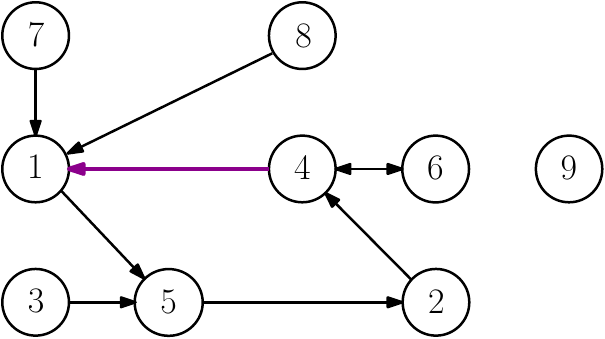}
     \end{subfigure}
     \hfill
     \vspace{7mm}
     \begin{subfigure}[b]{0.27\textwidth}
         \centering
         \includegraphics[width=\textwidth]{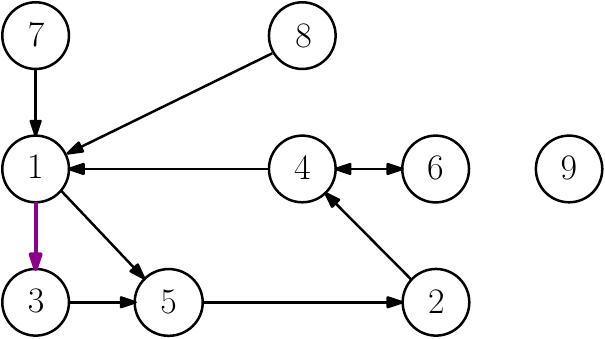}
     \end{subfigure}
     \hfill
     \begin{subfigure}[b]{0.27\textwidth}
         \centering
         \includegraphics[width=\textwidth]{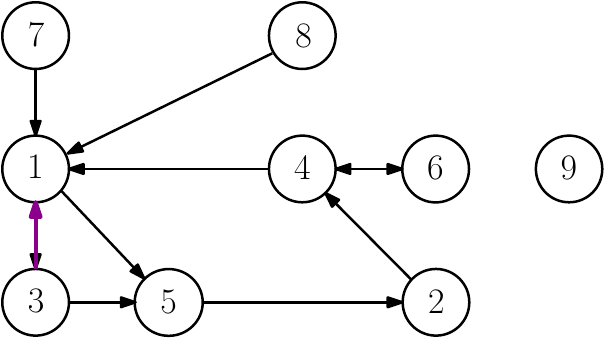}
     \end{subfigure}
     \hfill
     \begin{subfigure}[b]{0.27\textwidth}
         \centering
         \includegraphics[width=\textwidth]{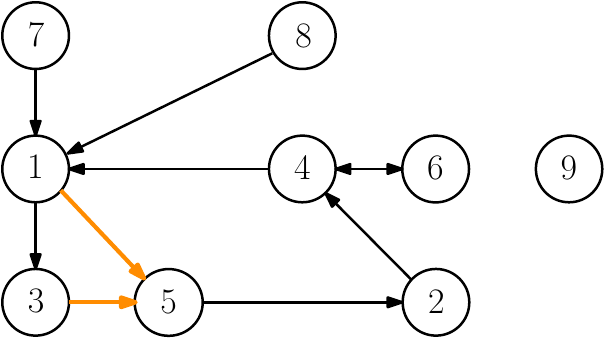}
     \end{subfigure}
     \hfill
     \vspace{7mm}

     \begin{subfigure}[b]{0.27\textwidth}
         \centering
         \includegraphics[width=\textwidth]{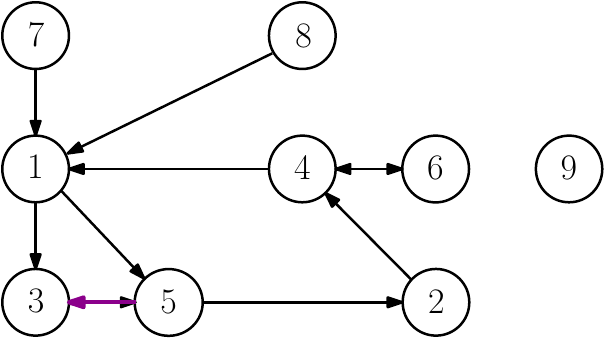}
     \end{subfigure}
     \hfill
     \begin{subfigure}[b]{0.27\textwidth}
         \centering
         \includegraphics[width=\textwidth]{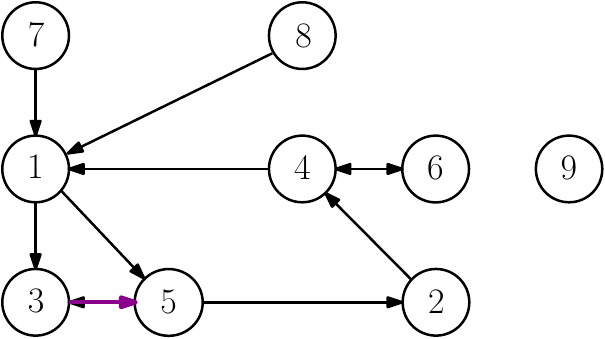}
     \end{subfigure}
     \hfill
     \begin{subfigure}[b]{0.27\textwidth}
         \centering
         \includegraphics[width=\textwidth]{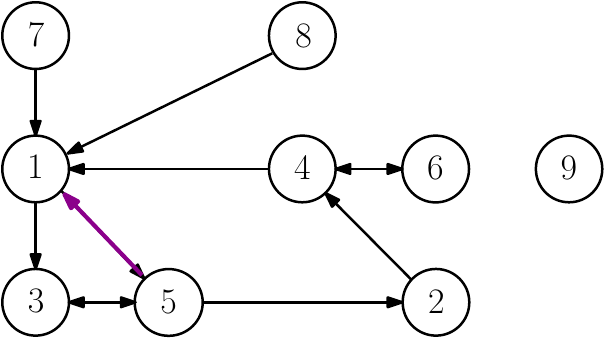}
     \end{subfigure}
     \hfill
     \vspace{7mm}

     \begin{subfigure}[b]{0.27\textwidth}
         \centering
         \includegraphics[width=\textwidth]{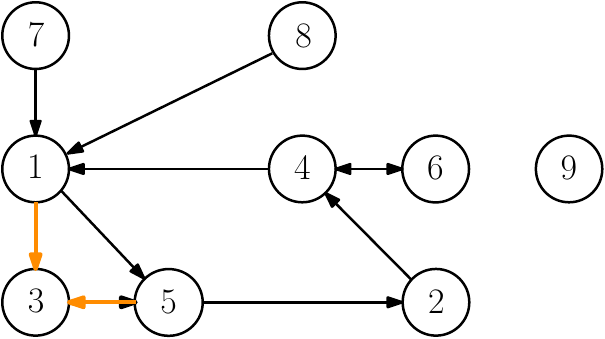}
     \end{subfigure}
     \hfill
     \begin{subfigure}[b]{0.27\textwidth}
         \centering
         \includegraphics[width=\textwidth]{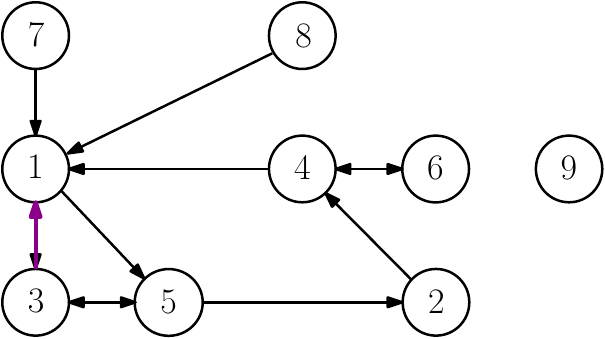}
     \end{subfigure}
     \hfill
     \begin{subfigure}[b]{0.27\textwidth}
         \centering
         \includegraphics[width=\textwidth]{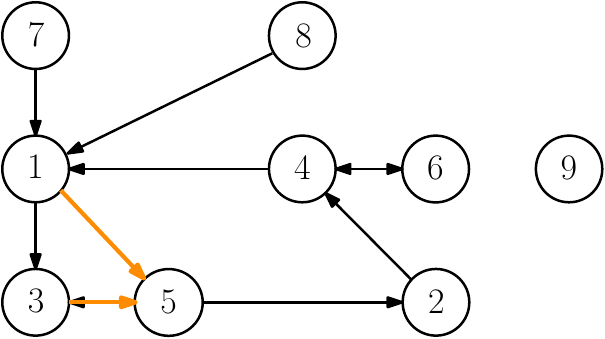}
     \end{subfigure}
     \hfill
     \vspace{7mm}
     \begin{subfigure}[b]{0.27\textwidth}
         \centering
         \includegraphics[width=\textwidth]{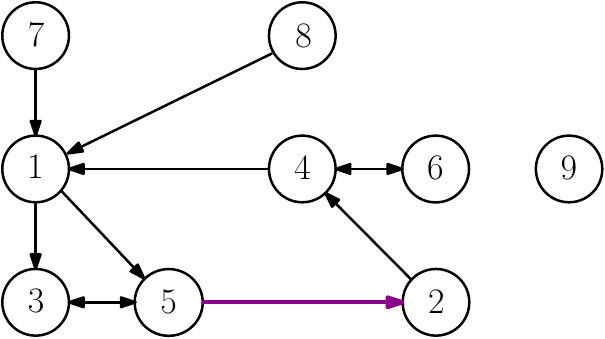}
     \end{subfigure}
     \hfill
     \begin{subfigure}[b]{0.27\textwidth}
         \centering
         \includegraphics[width=\textwidth]{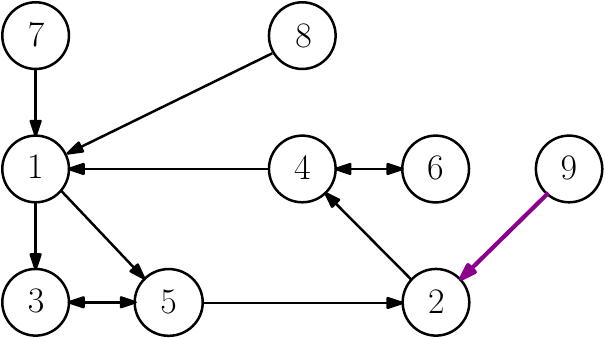}
     \end{subfigure}
     \hfill
     \begin{subfigure}[b]{0.27\textwidth}
         \centering
         \includegraphics[width=\textwidth]{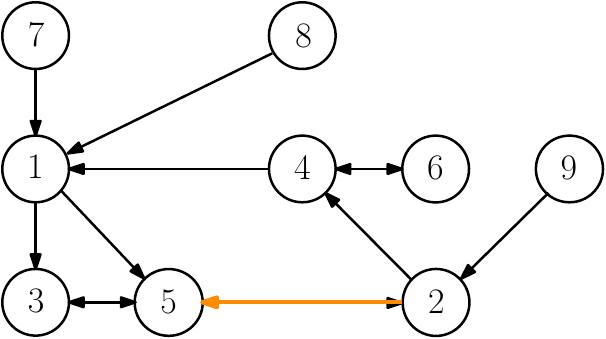}
     \end{subfigure}
     \hfill
     \vspace{3mm}
\caption{Attempt 1 in the execution of Algorithm~\ref{alg: SCCR2} 
. Each purple edge indicates the edge that is added to the construction at each step. If an unsafe edge is added, the algorithm initiates a correction process. Edges in orange are the ones added or flipped in such process. Attempt 1 is halted because after adding the unsafe edge $(9,2)$, this edge cannot be removed in the correction process, which results in the last $(5,2)$ added to the construction being flipped. Now $(2,5)$ is incompatible with $(1,5)$, which is an edge previously added in a correction process.}
\label{fig:example of SCCR2-1}
\end{figure}

\begin{figure}[H]
     \centering
     \begin{subfigure}[b]{0.27\textwidth}
         \centering
         \includegraphics[width=\textwidth]{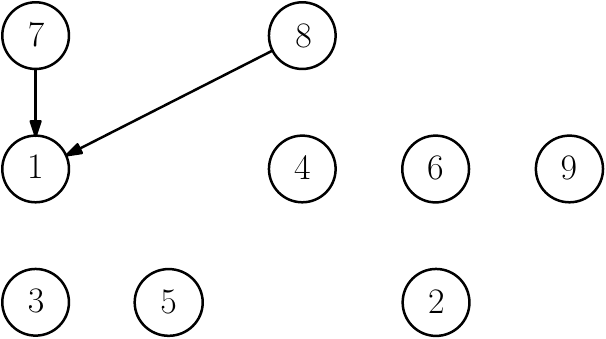}
     \end{subfigure}
     \hfill
     \begin{subfigure}[b]{0.27\textwidth}
         \centering
         \includegraphics[width=\textwidth]{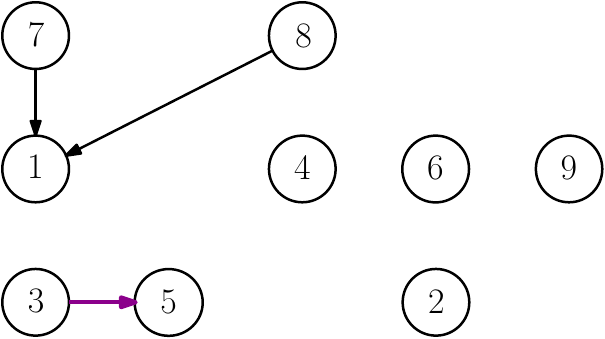}
     \end{subfigure}
     \hfill
     \begin{subfigure}[b]{0.27\textwidth}
         \centering
         \includegraphics[width=\textwidth]{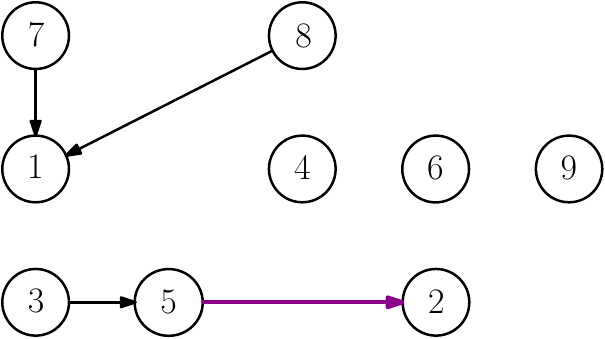}
     \end{subfigure}
     \hfill
     \vspace{7mm}

     \begin{subfigure}[b]{0.27\textwidth}
         \centering
         \includegraphics[width=\textwidth]{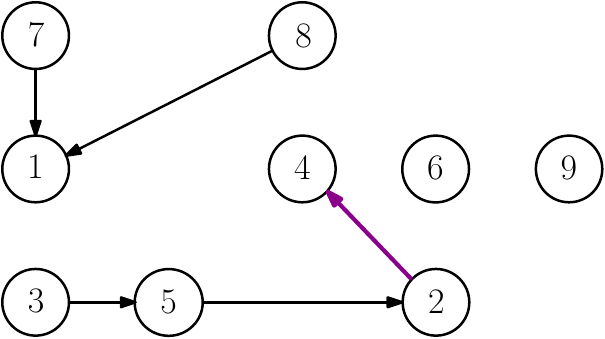}
     \end{subfigure}
     \hfill
     \begin{subfigure}[b]{0.27\textwidth}
         \centering
         \includegraphics[width=\textwidth]{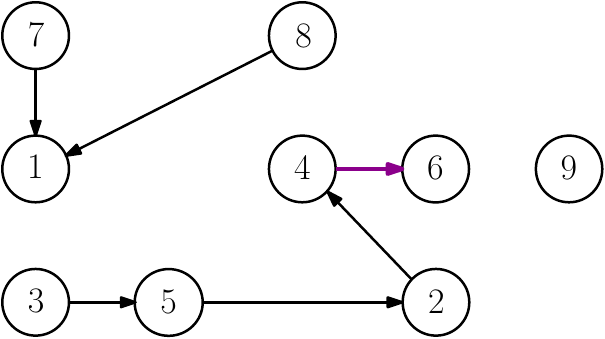}
     \end{subfigure}
     \hfill
     \begin{subfigure}[b]{0.27\textwidth}
         \centering
         \includegraphics[width=\textwidth]{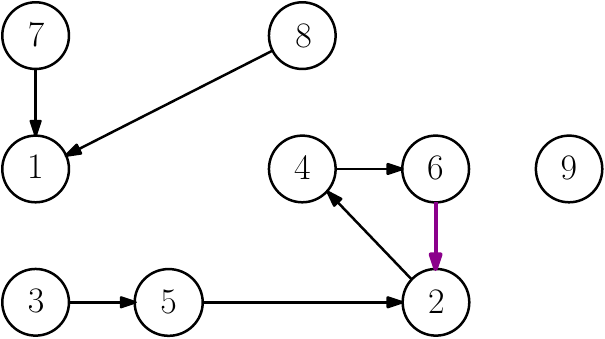}
     \end{subfigure}
     \hfill
     \vspace{7mm}

     \begin{subfigure}[b]{0.27\textwidth}
         \centering
         \includegraphics[width=\textwidth]{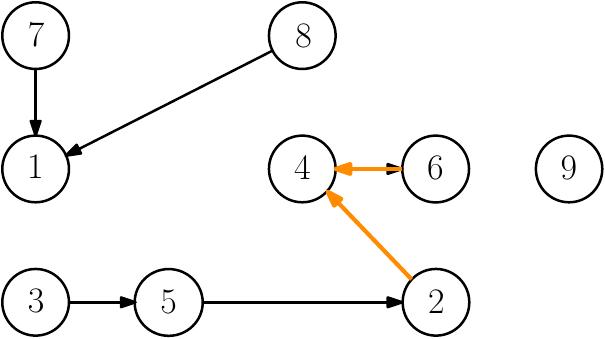}
     \end{subfigure}
     \hfill
     \begin{subfigure}[b]{0.27\textwidth}
         \centering
         \includegraphics[width=\textwidth]{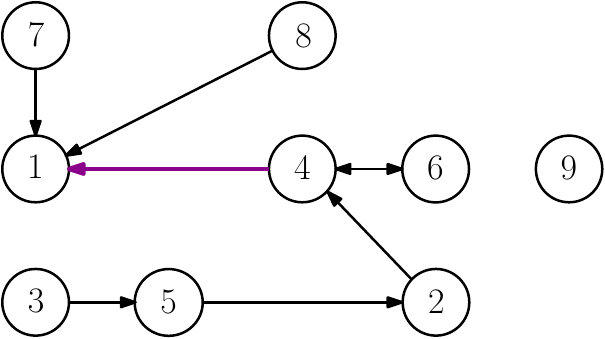}
     \end{subfigure}
     \hfill
     \begin{subfigure}[b]{0.27\textwidth}
         \centering
         \includegraphics[width=\textwidth]{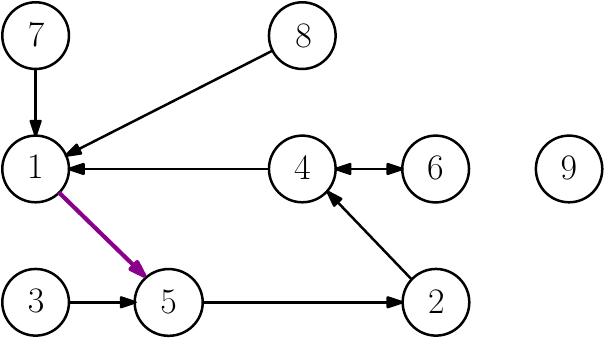}
     \end{subfigure}
     \hfill
     \vspace{7mm}

     \begin{subfigure}[b]{0.27\textwidth}
         \centering
         \includegraphics[width=\textwidth]{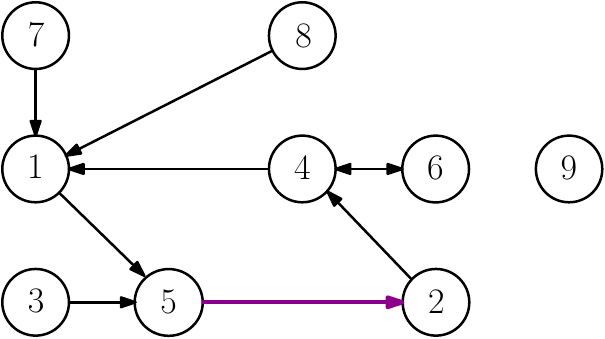}
     \end{subfigure}
     \hfill
     \begin{subfigure}[b]{0.27\textwidth}
         \centering
         \includegraphics[width=\textwidth]{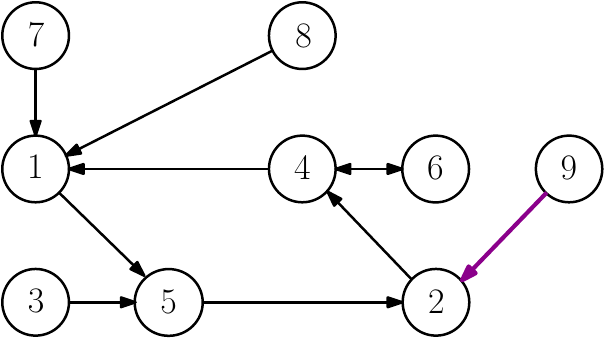}
     \end{subfigure}
     \hfill
     \begin{subfigure}[b]{0.27\textwidth}
         \centering
         \includegraphics[width=\textwidth]{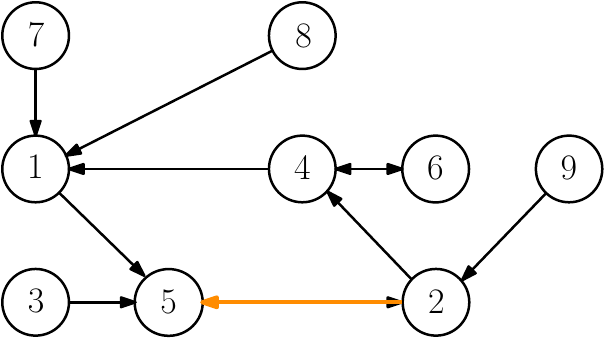}
     \end{subfigure}
     \hfill
     \vspace{7mm}
     \begin{subfigure}[b]{0.27\textwidth}
         \centering
         \includegraphics[width=\textwidth]{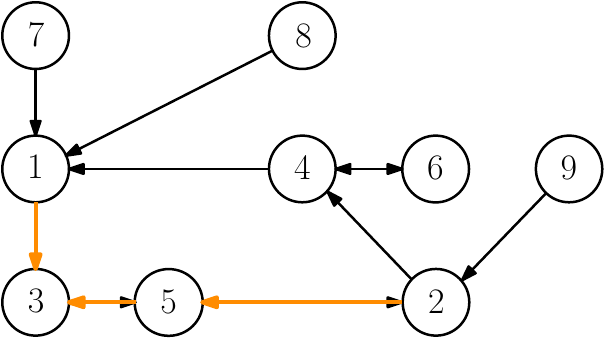}
     \end{subfigure}
     \hfill
     \begin{subfigure}[b]{0.27\textwidth}
         \centering
     \end{subfigure}
     \hfill
     \begin{subfigure}[b]{0.27\textwidth}
         \centering
     \end{subfigure}
     \hfill
     \vspace{3mm}
\caption{Attempt 2 in the execution of Algorithm~\ref{alg: SCCR2}
. Each purple edge indicates the edge that is added to the construction at each step. If an unsafe edge is added, the algorithm initiates a correction process. Edges in orange are the ones added or flipped in such process. Attempt 2 is halted because by the end of the last correction process initiated by adding the unsafe edge $(9,2)$, $(9,2)$ is still incompatible with the edge $(5,2)$ in the construction.}
\label{fig:example of SCCR2-2}
\end{figure}

In this section, we illustrate the step by step execution of Algorithm~\ref{alg: SCCR2} on the following input:
\begin{align*}
&C \coloneqq \{1,2,3,4,5,6\}, \\
&A_C\coloneqq \sym\{(1,3),(3,5),(1,5),(2,4),(2,6),(4,6),(1,4),(2,5)\},\\
&B_C \coloneqq \{(7,1),(8,1),(7,4),(8,4),(9,6),(9,2) \}, \\
&ComCh_C\coloneqq \sym\{(7,8)\}, \quad NoComCh_C \coloneqq \sym\{(7,9)\},\quad N\coloneqq 100.
\end{align*}
The algorithm outputs the set 
\begin{align*}
    E_C = \{(1,3),(3,5),(5,3),(1,5),(4,1),(5,2),(2,4),(2,6),(6,4),(7,1),(8,1),(9,6) \}.
\end{align*}
$E_C$ is outputted after three attempts, i.e., by the end of the execution, $I=3$. Figures~\ref{fig:example of SCCR2-1},~\ref{fig:example of SCCR2-2} and~\ref{fig:example of SCCR2-3} show the steps when $I=1$, $I=2$ and $I=3$ respectively. 

\begin{figure}[H]
     \centering
     \begin{subfigure}[b]{0.27\textwidth}
         \centering
         \includegraphics[width=\textwidth]{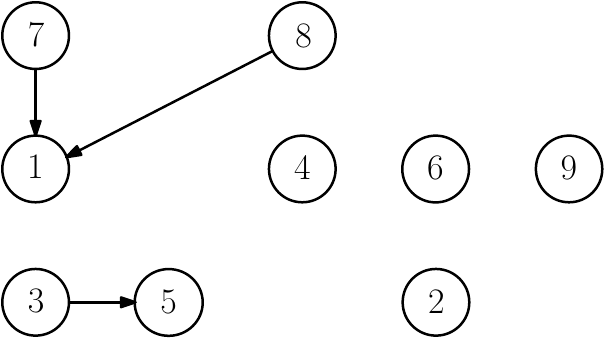}
     \end{subfigure}
     \hfill
     \begin{subfigure}[b]{0.27\textwidth}
         \centering
         \includegraphics[width=\textwidth]{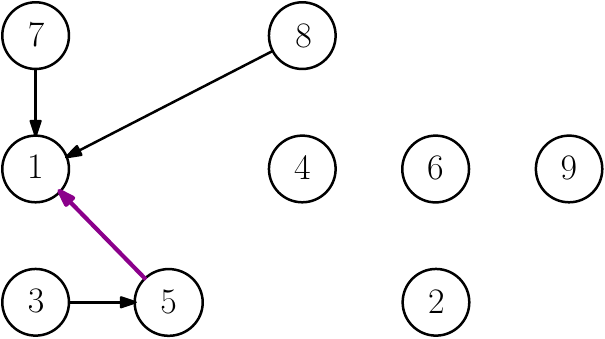}
     \end{subfigure}
     \hfill
     \begin{subfigure}[b]{0.27\textwidth}
         \centering
         \includegraphics[width=\textwidth]{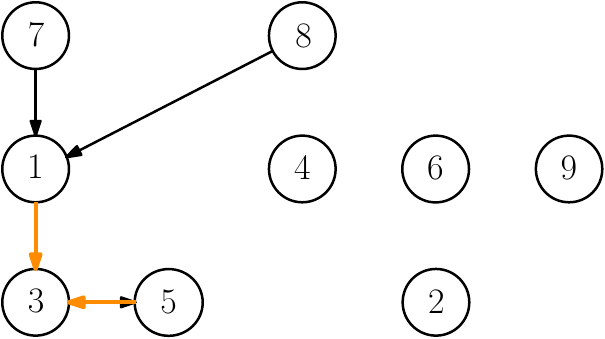}
     \end{subfigure}
     \hfill
\vspace{3mm}
\caption{Attempt 3 in the execution of Algorithm~\ref{alg: SCCR2}.}
\label{fig:example of SCCR2-3}
\end{figure}
\begin{figure}[H]
\ContinuedFloat
\centering
     \begin{subfigure}[b]{0.27\textwidth}
         \centering
         \includegraphics[width=\textwidth]{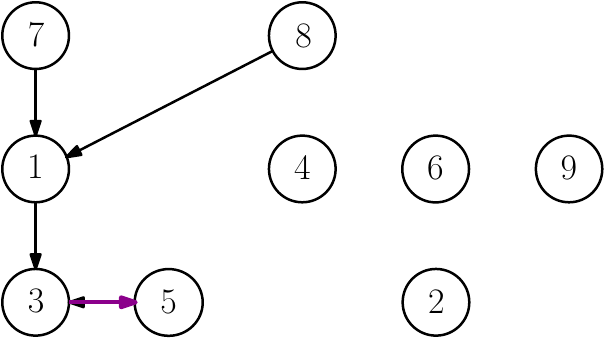}
     \end{subfigure}
     \hfill
     \begin{subfigure}[b]{0.27\textwidth}
         \centering
         \includegraphics[width=\textwidth]{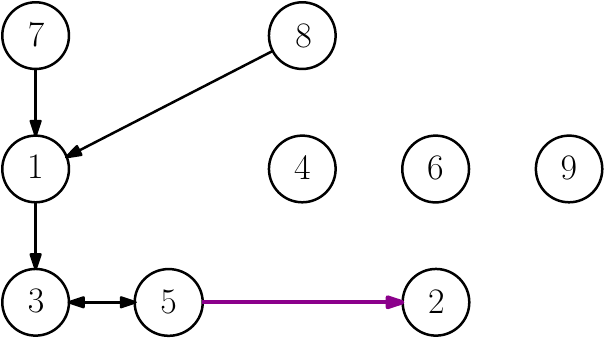}
     \end{subfigure}
     \hfill
     \begin{subfigure}[b]{0.27\textwidth}
         \centering
         \includegraphics[width=\textwidth]{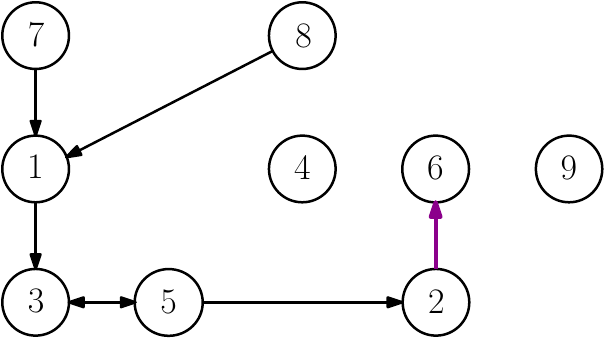}
     \end{subfigure}
     \hfill
     \vspace{7mm}

     \begin{subfigure}[b]{0.27\textwidth}
         \centering
         \includegraphics[width=\textwidth]{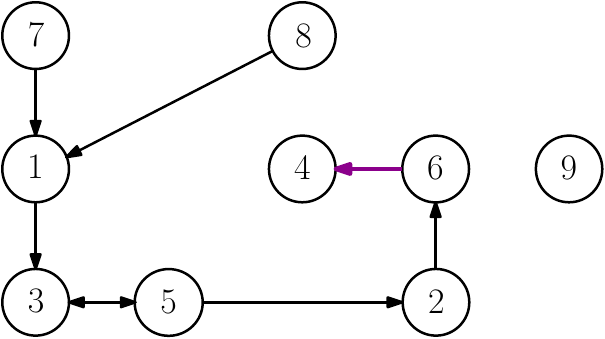}
     \end{subfigure}
     \hfill
     \begin{subfigure}[b]{0.27\textwidth}
         \centering
         \includegraphics[width=\textwidth]{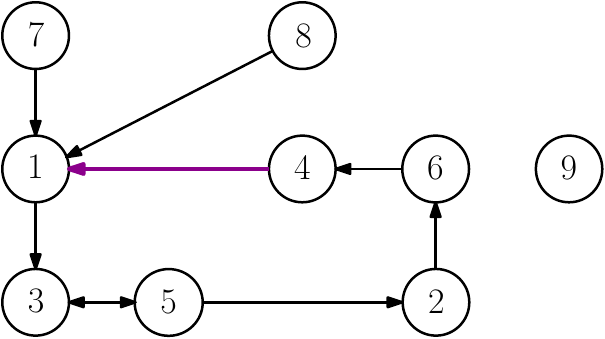}
     \end{subfigure}
     \hfill
     \begin{subfigure}[b]{0.27\textwidth}
         \centering
         \includegraphics[width=\textwidth]{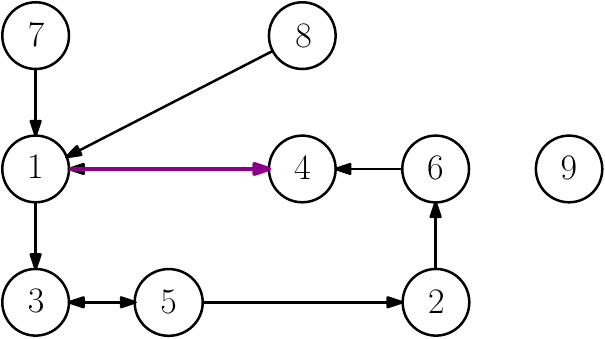}
     \end{subfigure}
     \hfill
     \vspace{7mm}

     \begin{subfigure}[b]{0.27\textwidth}
         \centering
         \includegraphics[width=\textwidth]{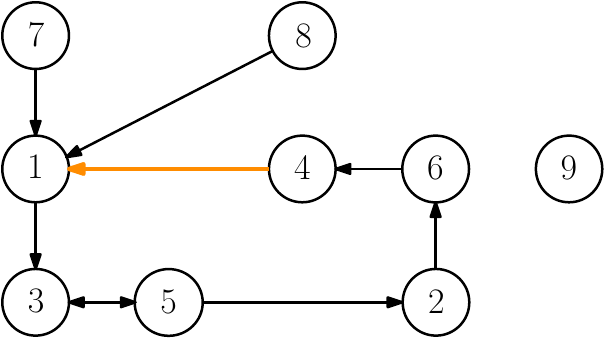}
     \end{subfigure}
     \hfill
     \begin{subfigure}[b]{0.27\textwidth}
         \centering
         \includegraphics[width=\textwidth]{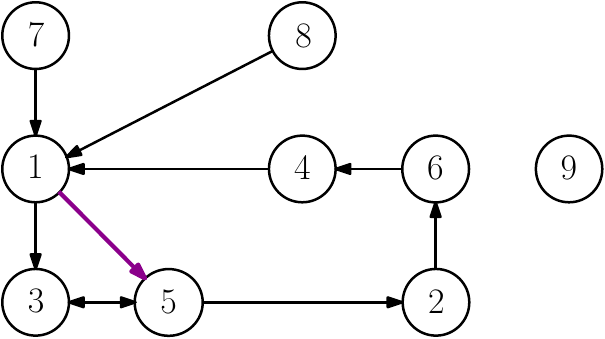}
     \end{subfigure}
     \hfill
     \begin{subfigure}[b]{0.27\textwidth}
         \centering
         \includegraphics[width=\textwidth]{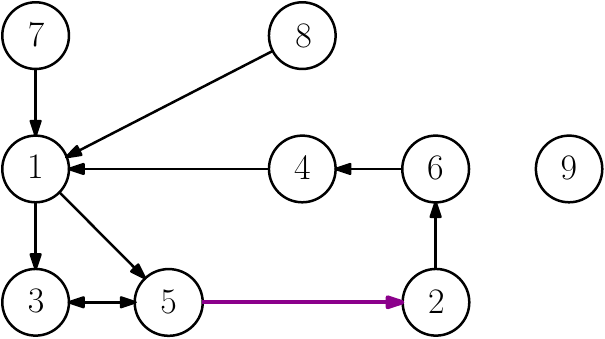}
     \end{subfigure}
     \hfill
     \vspace{7mm}

     \begin{subfigure}[b]{0.27\textwidth}
         \centering
         \includegraphics[width=\textwidth]{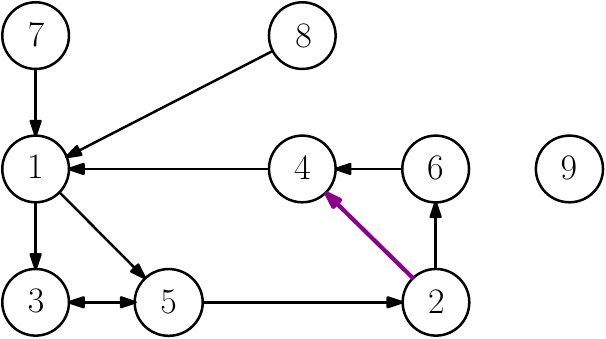}
     \end{subfigure}
     \hfill
     \begin{subfigure}[b]{0.27\textwidth}
         \centering
         \includegraphics[width=\textwidth]{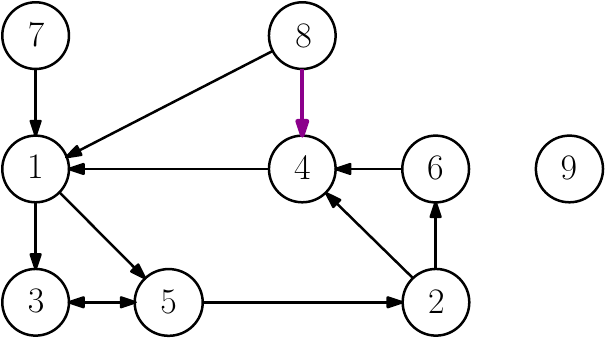}
     \end{subfigure}
     \hfill
     \begin{subfigure}[b]{0.27\textwidth}
         \centering
         \includegraphics[width=\textwidth]{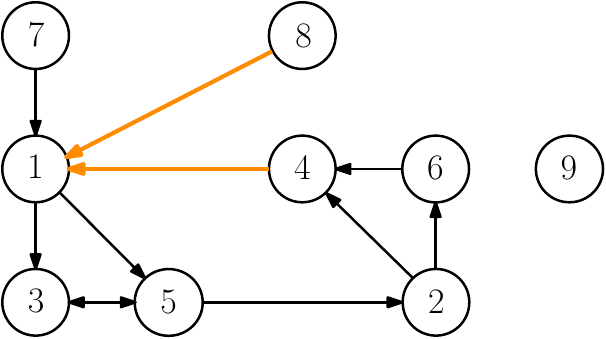}
     \end{subfigure}
     \hfill
     \vspace{7mm}
     \begin{subfigure}[b]{0.27\textwidth}
         \centering
         \includegraphics[width=\textwidth]{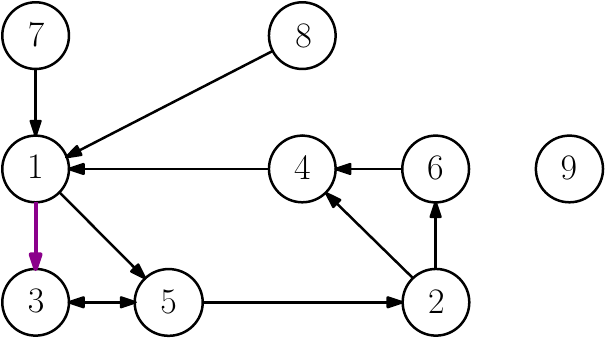}
     \end{subfigure}
     \hfill
     \begin{subfigure}[b]{0.27\textwidth}
         \centering
         \includegraphics[width=\textwidth]{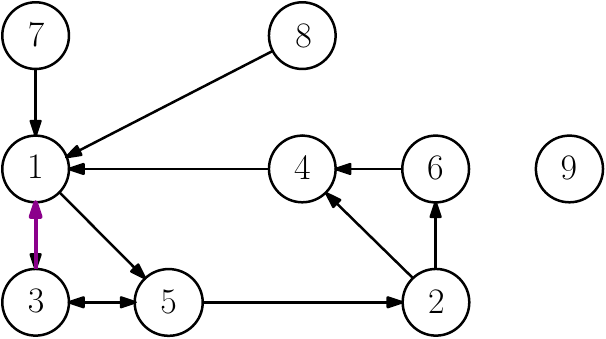}
     \end{subfigure}
     \hfill
     \begin{subfigure}[b]{0.27\textwidth}
         \centering
         \includegraphics[width=\textwidth]{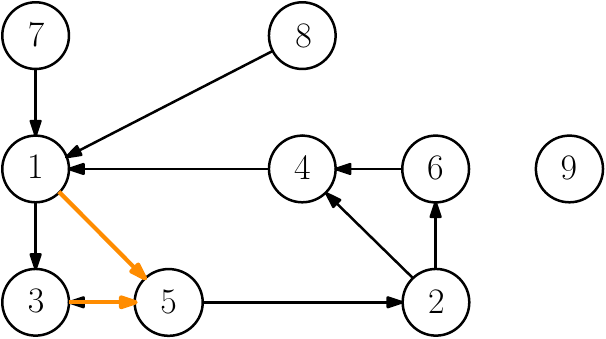}
     \end{subfigure}
     \hfill
     \vspace{7mm}

     \begin{subfigure}[b]{0.27\textwidth}
         \centering
         \includegraphics[width=\textwidth]{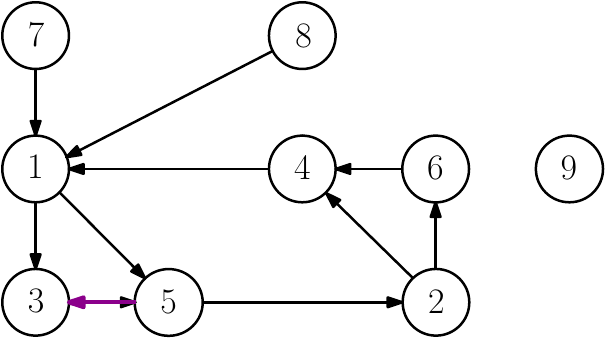}
     \end{subfigure}
     \hfill
     \begin{subfigure}[b]{0.27\textwidth}
         \centering
         \includegraphics[width=\textwidth]{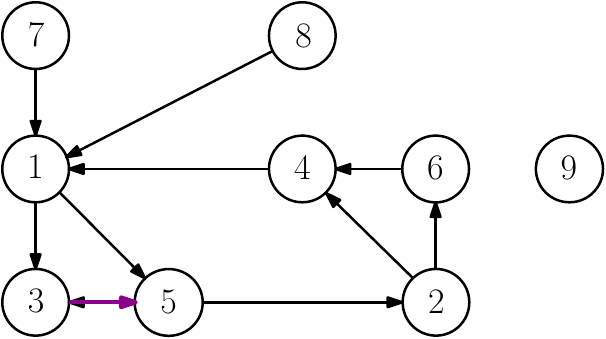}
     \end{subfigure}
     \hfill
     \begin{subfigure}[b]{0.27\textwidth}
         \centering
         \includegraphics[width=\textwidth]{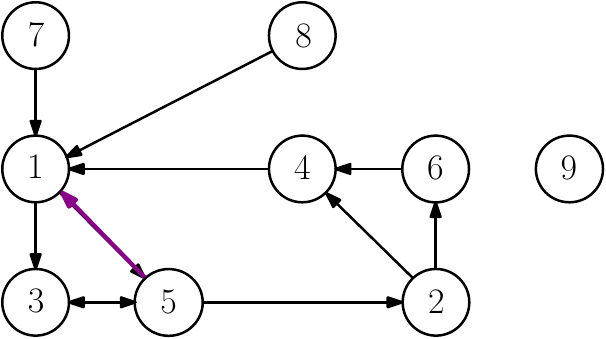}
     \end{subfigure}
     \hfill
     \vspace{7mm}

     \begin{subfigure}[b]{0.27\textwidth}
         \centering
         \includegraphics[width=\textwidth]{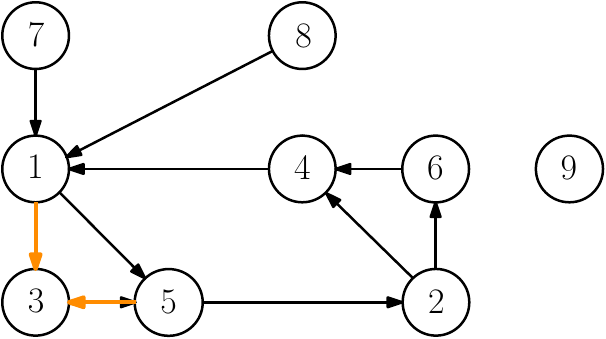}
     \end{subfigure}
     \hfill
     \begin{subfigure}[b]{0.27\textwidth}
         \centering
         \includegraphics[width=\textwidth]{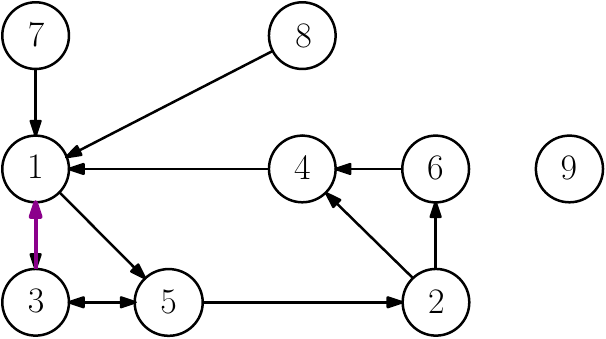}
     \end{subfigure}
     \hfill
     \begin{subfigure}[b]{0.27\textwidth}
         \centering
         \includegraphics[width=\textwidth]{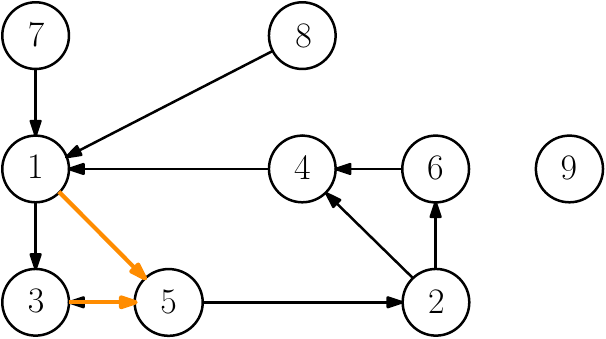}
     \end{subfigure}
     \hfill
\vspace{3mm}
\caption{Attempt 3 in the execution of Algorithm~\ref{alg: SCCR2} (continued).}
\end{figure}
\begin{figure}[H]
\ContinuedFloat
\centering
     \begin{subfigure}[b]{0.27\textwidth}
         \centering
         \includegraphics[width=\textwidth]{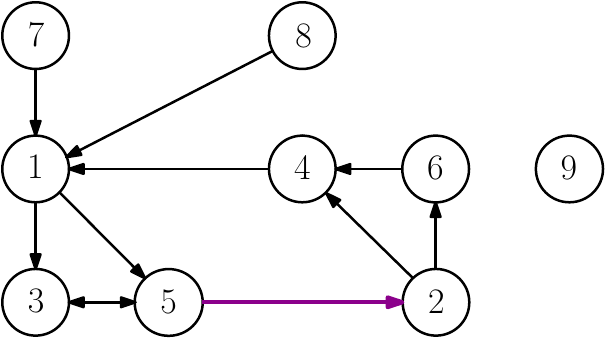}
     \end{subfigure}
     \hfill
     \begin{subfigure}[b]{0.27\textwidth}
         \centering
         \includegraphics[width=\textwidth]{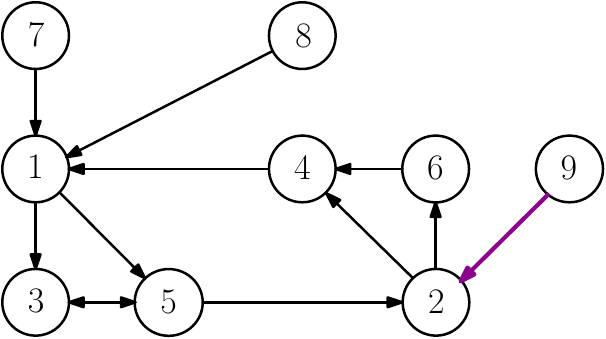}
     \end{subfigure}
     \hfill
     \begin{subfigure}[b]{0.27\textwidth}
         \centering
         \includegraphics[width=\textwidth]{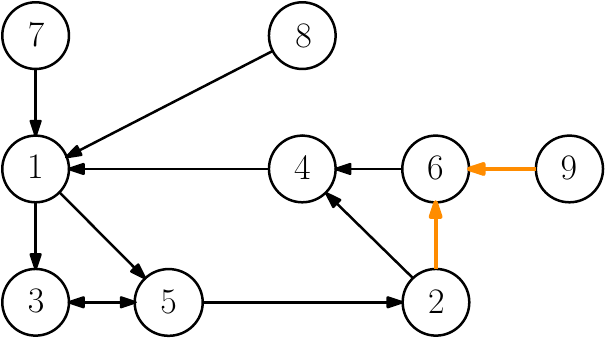}
     \end{subfigure}
     \hfill
     \vspace{7mm}

     \begin{subfigure}[b]{0.27\textwidth}
         \centering
         \includegraphics[width=\textwidth]{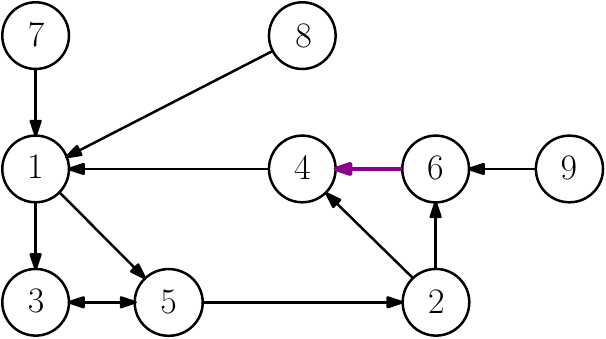}
     \end{subfigure}
     \hfill
     \begin{subfigure}[b]{0.27\textwidth}
         \centering
         \includegraphics[width=\textwidth]{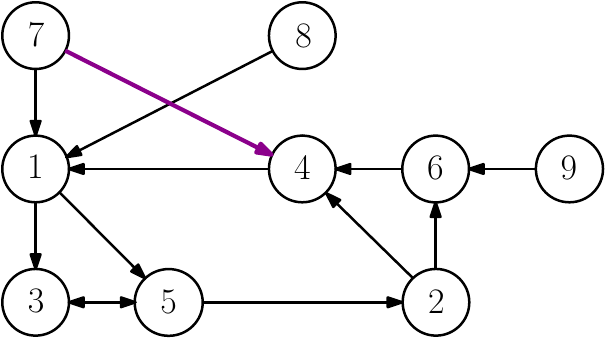}
     \end{subfigure}
     \hfill
     \begin{subfigure}[b]{0.27\textwidth}
         \centering
         \includegraphics[width=\textwidth]{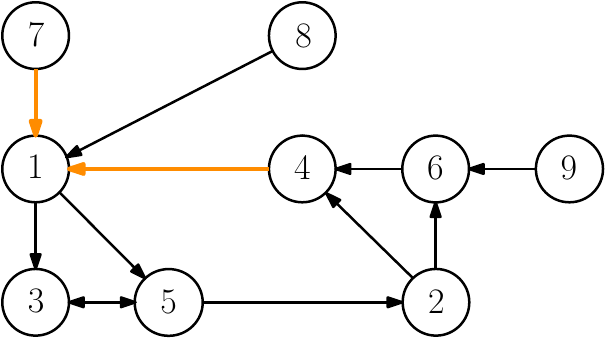}
     \end{subfigure}
     \hfill
     \vspace{3mm}
\caption{Attempt 3 in the execution of Algorithm~\ref{alg: SCCR2} (continued). Each purple edge indicates the edge that is added to the construction at each step. If an unsafe edge is added, the algorithm initiates a correction process. Edges in orange are the ones added or flipped in such process. This attempt ends successfully.}
\end{figure}
\section{Definition of a submodular flow polyhedron} \label{app: submodular flow polyhedron}
In this section, we define a submodular flow polyhedron. Throughout this section, assume $G=(V,E)$ is a directed graph.
 \begin{definition}
    Two subsets $S,T\subseteq V$ are called \textit{crossing} if all the sets $S\cap T$, $S\setminus T$, $T\setminus S$ and $V\setminus (S\cup T)$ are nonempty.

    A family $\mathcal{C}$ of subsets of $V$ is called \textit{crossing} if for any crossing pair $S,T\in \mathcal{C}$, we have $S\cap T\in \mathcal{C}$ and $S\cup T\in \mathcal{C}$.
\end{definition}
\begin{definition}
    Let $\mathcal{C}$ be a crossing family of subsets of $V$. A function $f:\mathcal{C}\to \mathbb{R}$ is called \textit{crossing submodular} if for every crossing pair $S,T\in \mathcal{C}$, we have
    \begin{align*}
        f(S)+f(T) \geq f(S\cap T) + f(S\cup T).
    \end{align*}
\end{definition}
\begin{definition}
    A submodular flow polyhedron corresponding to $G$, a crossing subfamily $\mathcal{C}$ of the subsets of $V$ and a crossing submodular function $f:\mathcal{C}\to \mathbb{R}$ is defined as follows:
    \begin{align*}
        \Set{x\in \mathbb{R}^E | x\left( \delta^-_E(U)\right) - x\left( \delta^+_E(U)\right)\leq f(U) \text{ for all } U\in \mathcal{C}.}.
    \end{align*}
    See Definition~\ref{def: delta-+} for the definitions of $x\left( \delta^-_E(U)\right)$ and $x\left( \delta^+_E(U)\right)$.
\end{definition}
\section{Proofs for Subsection~\ref{subsection: SCCR1}}\label{app: proofs of submodular}
\begin{proof}[Proof of Lemma~\ref{lemma: submodular flow}]
($\Longrightarrow$) Let $C$ be a strongly connected component in the graph $\left([n],\eyj\right)$. Define $x\in \{-1,0\}^{\tilde{A}_C\setminus Y}$ as follows:
\begin{align*}
    \text{For all } e\in \tilde{A}_C\setminus Y,\quad x_e\coloneqq \begin{cases} -1 & \text{if $e\in J$.} \\ 0 & \text{if $e\not\in J$.}\end{cases} 
\end{align*}
We prove $x\in \mathcal{F}(Y)$. Let $U\in 2^C\setminus \{\emptyset,C\}$. Then
\begin{align*}
\left| \gtp(U)\right| - x\left(\gtm(U)\right) + x\left(\gtp(U)\right)
&=\left| \gtp(U)\right| + \left| \gtm(U)\cap J\right| - \left|\gtp(U)\cap J\right| \\
&=\left| \gtjp(U)\right| \geq 1,
\end{align*}
where the last inequality results from the fact that $C$ is a strongly connected component in the graph $\left([n],\eyj\right)$. So, 
\begin{align*}
x\left(\gtm(U)\right) - x\left(\gtp(U)\right) \leq \left| \gtp(U) \right| - 1. 
\end{align*}
Hence, $x\in \mathcal{F}(Y)$.\\

($\Longleftarrow$) Now let $x\in  \mathcal{F}(Y)\cap \{-1,0\}^{\tilde{A}_C\setminus Y}$. Define $J\coloneqq \Set{e\in \tilde{A}_C\setminus Y| x_e=-1}$. We prove $C$ is a strongly connected component in the graph $\left([n],\eyj\right)$. For all $U\in 2^C\setminus \{\emptyset, C\}$, 
\begin{align*}
    x\left(\gtm(U)\right) - x\left(\gtp(U)\right) \leq \left| \gtp(U)\right| -1. 
\end{align*}
Therefore, 
\begin{align*}
    \left|\gtjp(U)\right| &= \left|\gtp(U)\right| + \left|\gtm(U) \cap J\right| - \left| \gtp(U)\cap J\right|\\
    &= \left| \gtp(U)\right| -x\left(\gtm(U)\right) + x\left(\gtp(U)\right) 
    \geq 1. 
\end{align*}
This implies that $C$ is a strongly connected component in $\left([n],\eyj\right)$.
\end{proof}
\begin{table}
        \centering
        \begin{tabular}{m{26mm}|m{26mm}|m{26mm}|m{26mm}|m{26mm}}
           \diagbox[innerwidth=26mm]{$\tilde{E}_C$}{$\eyjx$} 
           & \includegraphics[scale=0.35]{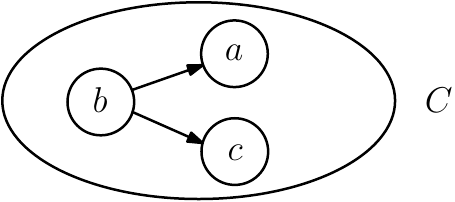}  
           &  \includegraphics[scale=0.35]{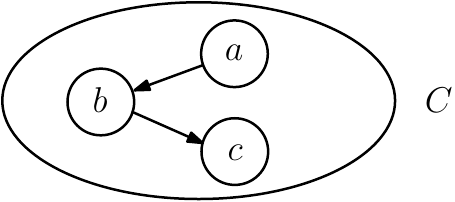} 
           & \includegraphics[scale=0.35]{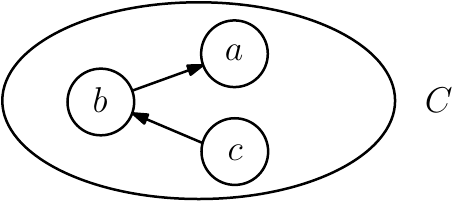} 
           & \includegraphics[scale=0.35]{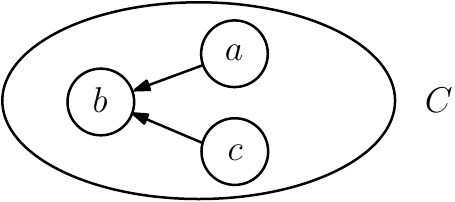} 
           \\
           \hline
            \vspace{2mm}
           \includegraphics[scale=0.35]{figures/fig1.pdf}  & $$0$$ & $$1$$ & $$1$$ & $$2$$ \\
           \hline
            \vspace{2mm}
           \includegraphics[scale=0.35]{figures/fig2.pdf}  & $$-1$$ & $$0$$ & $$0$$ & $$1$$ \\
           \hline
           \vspace{2mm}
           \includegraphics[scale=0.35]{figures/fig4.pdf}  & $$-1$$ & $$0$$ & $$0$$ & $$1$$ \\
           \hline
           \vspace{2mm}
           \includegraphics[scale=0.35]{figures/fig3.pdf}  & $$-2$$ & $$-1$$ & $$-1$$ & $$0$$ 
             \end{tabular}
        \caption{For a removable set $Y\subseteq \tilde{E}_C$, a vertex $x$ of $\mathcal{F}(Y)\cap [-1,0]^{\tilde{A}_C\setminus Y}$ and $(a,b),(c,b) \in A_C$ with $Y\cap \{(a,b),(b,a),(c,b),(b,c)\}=\emptyset$, each cell shows $w_{Y,(a,b,c)}^Tx$ if the orientations of $(a,b),(c,b)$ in the graphs $\left([n],\tilde{E}_C\right)$ and $\left([n],\eyjx\right)$ are as depicted in the corresponding row and column respectively.}
        \label{tab: table1}
    \end{table}
       \begin{table}
        \centering
        \begin{tabular}{m{26mm}|m{26mm}|m{26mm}}
           \diagbox[innerwidth=26mm,height=15mm]{$\tilde{E}_C$}{$\eyjx$} 
           & \includegraphics[scale=0.35]{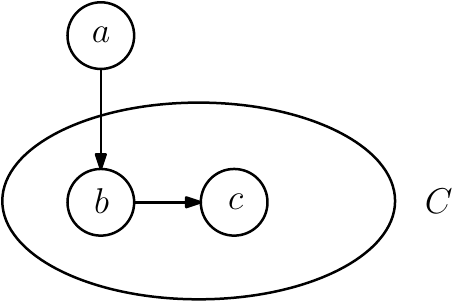}  
           &  \includegraphics[scale=0.35]{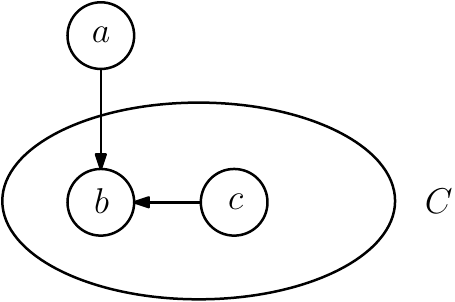} 
           \\
           \hline
            \vspace{2mm}
           \includegraphics[scale=0.35]{figures/fig6.pdf}  & $$0$$ & $$1$$ \\
           \hline
            \vspace{2mm}
           \includegraphics[scale=0.35]{figures/fig5.pdf}  & $$-1$$ & $$0$$ 
             \end{tabular}
        \caption{For a removable set $Y\subseteq \tilde{E}_C$, a vertex $x$ of $\mathcal{F}(Y)\cap [-1,0]^{\tilde{A}_C\setminus Y}$ and $(a,b)\in B_C,(c,b) \in A_C$ with $Y\cap \{(a,b),(b,a),(c,b),(b,c)\}=\emptyset$, each cell shows $w_{Y,(a,b,c)}^Tx$ if the orientation of $(c,b)$ in the graphs $\left([n],\tilde{E}_C\right)$ and $\left([n],\eyjx\right)$ is as depicted in the corresponding row and column respectively.}
        \label{tab: table2}
    \end{table}    
\begin{proof}[Proof of Theorem~\ref{thm: SCCR1}]
($\Longleftarrow$) Let $Y\subseteq \tilde{E}_C$ be a removable set and $x$ be a vertex of the polyhedron $\mathcal{F}(Y)\cap [-1,0]^{\tilde{A}_C\setminus Y}$ satisfying the properties listed above. We prove $\eyjx$ satisfies the properties of Definition~\ref{def: strongly connected component recovery algorithm}. By \cite{EDMONDS1977185}, any submodular flow polyhedron is box-integer if its crossing submodular function is integer-valued. So, every vertex of $\mathcal{F}(Y)\cap [-1,0]^{\tilde{A}_C\setminus Y}$ is an integer-valued vector, and therefore, by Lemma~\ref{lemma: submodular flow}, $C$ is a strongly connected component in the graph $\left([n],\eyjx\right)$. Moreover, $Y$ being removable secures a common child in $C$ for every pair in $ComCh_C$ and prevents every pair in $NoComCh_C$ from having a common child in $C$.

For any $(a,b,c)\in [n]^3$ with $c\in C$, $(a,b),(c,b) \in A_C\cup B_C$ and $Y\cap \{(a,b),(b,a),(c,b),(b,c)\}=\emptyset$, each cell in Tables~\ref{tab: table1} and~\ref{tab: table2} shows $w_{Y,(a,b,c)}^Tx$ if the orientations of $(a,b),(c,b)$ in the graphs $\left([n],\tilde{E}_C\right)$ and $\left([n],\eyjx\right)$ are as illustrated in the corresponding row and column respectively. These tables confirm that $(a,b),(c,b)\in \eyjx$ if and only if 
$w_{Y,(a,b,c)}^Tx = \max_{y\in \{-1,0\}^{\tilde{A}_C\setminus Y}} w_{Y,(a,b,c)}^T y$. This observation along with properties~\ref{thm: SCCR1: property 1} and~\ref{thm: SCCR1: property 2} imply that for any $a,c\in [n]$, $a$ and $c$ are $p$-adjacent in $\left([n],\eyjx\right)$ if and only if $(a,c)\in A_C\cup B_C$ or $(c,a)\in A_C\cup B_C$.

($\Longrightarrow$) Let $E_C\subseteq [n]^2$ be a set with the properties mentioned in Definition~\ref{def: strongly connected component recovery algorithm} and containing at most one direction of each edge. Then $E_C$ has been obtained from $\tilde{E}_C$ by removing the edges in some subset $Y\subseteq \tilde{E}_C$ and flipping the edges in some subset $J\subseteq \tilde{A}_C\setminus Y$. Since in $\left([n],E_C\right)$, no pair in $NoComCh_C$ has a common child in $C$, every pair in $ComCh_C$ does have a common child in $C$, and every $p$-adjacency represented by an omitted edge $e\in Y$ is preserved, $Y$ is a removable set. Choose $x\in \{-1,0\}^{\tilde{A}_C\setminus Y}$ such that $J=J(x)$. Then $E_C=\eyjx$, and by Lemma~\ref{lemma: submodular flow}, $x\in \mathcal{F}(Y)$ given that $C$ is a strongly connected component in $\left([n],E_C\right)$. So, since $x$ is an integer-valued vector, $x$ has to be a vertex of the polyhedron $\mathcal{F}(Y)\cap [-1,0]^{\tilde{A}_C\setminus Y}$. Combining the observation previously made from Tables~\ref{tab: table1} and~\ref{tab: table2} with the fact that for any $a,c\in [n]$, $a$ and $c$ are $p$-adjacent in $\left([n],E_C \right)$ if and only if $(a,c)\in A_C\cup B_C$ or $(c,a)\in A_C\cup B_C$ guarantees properties~\ref{thm: SCCR1: property 1} and~\ref{thm: SCCR1: property 2}. 
\end{proof}

\section*{Acknowledgements}
PS and ER were supported by NSERC Discovery Grant DGECR-2020-00338. Moreover, PS was supported by the Vanier Canada Graduate Scholarship, and ER is a Canada CIFAR AI Chair.

\bibliographystyle{vancouver}
\bibliography{biblio}

\begin{thebibliography}{10}

\bibitem{Friedman2000}
Friedman N, Linial M, Nachman I, Pe'er D.
\newblock Using {B}ayesian networks to analyze expression data.
\newblock Journal of Computational Biology. 2000;7(3-4):601-20.
\newblock PMID: 11108481.
\newblock Available from: \url{https://doi.org/10.1089/106652700750050961}.

\bibitem{pearlcausalitymodels}
Pearl J.
\newblock Causality: Models, Reasoning and Inference.
\newblock 2nd ed. USA: Cambridge University Press; 2009.

\bibitem{Robins2000}
Robins JM, Hernán MA, Brumback B.
\newblock Marginal structural models and causal inference in epidemiology.
\newblock Epidemiology. 2000;11(5):550-60.
\newblock Available from: \url{https://doi.org/10.1097/00001648-200009000-00011}.

\bibitem{SpirtesGlymourClark1993}
Spirtes P, Glymour C, Scheines R.
\newblock Causation, Prediction, and Search. vol.~81.
\newblock Springer New York; 1993.

\bibitem{Mason1953}
Mason SJ.
\newblock Feedback theory-some properties of signal flow graphs.
\newblock Proceedings of the IRE. 1953;41(9):1144-56.
\newblock Available from: \url{https://doi.org/10.1109/JRPROC.1953.274449}.

\bibitem{Mason1956}
Mason SJ.
\newblock Feedback theory-further properties of signal flow graphs.
\newblock Proceedings of the IRE. 1956;44(7):920-6.
\newblock Available from: \url{https://doi.org/10.1109/JRPROC.1956.275147}.

\bibitem{Haavelmo1943}
Haavelmo T.
\newblock The statistical implications of a system of simultaneous equations.
\newblock Econometrica. 1943;11(1):1-12.
\newblock Available from: \url{http://www.jstor.org/stable/1905714}.

\bibitem{Goldberger1972}
Goldberger AS.
\newblock Structural equation methods in the social sciences.
\newblock Econometrica. 1972;40(6):979-1001.
\newblock Available from: \url{http://www.jstor.org/stable/1913851}.

\bibitem{Chickering2002OptimalSI}
Chickering DM.
\newblock Optimal structure identification with greedy search.
\newblock J Mach Learn Res. 2002;3:507-54.
\newblock Available from: \url{https://jmlr.org/papers/v3/chickering02b.html}.

\bibitem{meek1997graphical}
Meek C.
\newblock Graphical models: Selecting causal and statistical models [{PhD} thesis].
\newblock Carnegie Mellon University; 1997.

\bibitem{Verma1990EquivalenceAS}
Verma T, Pearl J.
\newblock Equivalence and synthesis of causal models.
\newblock In: Proceedings of the Sixth Annual Conference on Uncertainty in Artificial Intelligence. UAI '90. USA: Elsevier Science Inc.; 1990. p. 255–270.

\bibitem{Teyssier2005OrderingBasedSA}
Teyssier M, Koller D.
\newblock Ordering-based search: A simple and effective algorithm for learning Bayesian networks.
\newblock In: Proceedings of the Twenty-First Conference on Uncertainty in Artificial Intelligence. UAI'05. Arlington, Virginia, USA: AUAI Press; 2005. p. 584–590.

\bibitem{Raskutti2018}
Raskutti G, Uhler C.
\newblock Learning directed acyclic graph models based on sparsest permutations.
\newblock Stat. 2018;7(1):e183.
\newblock E183 sta4.183.
\newblock Available from: \url{https://onlinelibrary.wiley.com/doi/abs/10.1002/sta4.183}.

\bibitem{Solus2021}
Solus L, Wang Y, Uhler C.
\newblock Consistency guarantees for greedy permutation-based causal inference algorithms.
\newblock Biometrika. 2021 01;108(4):795-814.
\newblock Available from: \url{https://doi.org/10.1093/biomet/asaa104}.

\bibitem{Wang2017}
Wang Y, Solus L, Yang KD, Uhler C.
\newblock Permutation-based causal inference algorithms with interventions.
\newblock In: Proceedings of the 31st International Conference on Neural Information Processing Systems. NIPS'17. Red Hook, NY, USA: Curran Associates Inc.; 2017. p. 5824–5833.

\bibitem{mohammadi2017generalized}
Mohammadi F, Uhler C, Wang C, Yu J.
\newblock Generalized permutohedra from probabilistic graphical models.
\newblock SIAM Journal on Discrete Mathematics. 2018;32(1):64-93.
\newblock Available from: \url{https://doi.org/10.1137/16M107894X}.

\bibitem{Bernstein2020}
Bernstein D, Saeed B, Squires C, Uhler C.
\newblock Ordering-based causal structure learning in the presence of latent variables.
\newblock In: Proceedings of the Twenty Third International Conference on Artificial Intelligence and Statistics. vol. 108. PMLR; 2020. p. 4098-108.

\bibitem{richardson2013discovery}
Richardson T.
\newblock A discovery algorithm for directed cyclic graphs.
\newblock In: Proceedings of the Twelfth International Conference on Uncertainty in Artificial Intelligence. UAI'96. San Francisco, CA, USA: Morgan Kaufmann Publishers Inc.; 1996. p. 454–461.

\bibitem{SAT2013}
Hyttinen A, Hoyer PO, Eberhardt F, J\"{a}rvisalo M.
\newblock Discovering cyclic causal models with latent variables: a general {SAT}-based procedure.
\newblock In: Proceedings of the Twenty-Ninth Conference on Uncertainty in Artificial Intelligence. UAI'13. Arlington, Virginia, USA: AUAI Press; 2013. p. 301–310.

\bibitem{claassen23a}
Claassen T, Mooij JM.
\newblock Establishing {M}arkov equivalence in cyclic directed graphs.
\newblock In: Proceedings of the Thirty-Ninth Conference on Uncertainty in Artificial Intelligence. vol. 216 of Proceedings of Machine Learning Research. PMLR; 2023. p. 433-42.
\newblock Available from: \url{https://proceedings.mlr.press/v216/claassen23a.html}.

\bibitem{Spirtes1995}
Spirtes P.
\newblock Directed cyclic graphical representations of feedback models.
\newblock In: Proceedings of the Eleventh Conference on Uncertainty in Artificial Intelligence. UAI'95. San Francisco, CA, USA: Morgan Kaufmann Publishers Inc.; 1995. p. 491–498.

\bibitem{Bongers_2021}
Bongers S, Forr{\'e} P, Peters J, Mooij JM.
\newblock {Foundations of structural causal models with cycles and latent variables}.
\newblock The Annals of Statistics. 2021;49(5):2885  2915.
\newblock Available from: \url{https://doi.org/10.1214/21-AOS2064}.

\bibitem{forré2017markov}
Forr{\'e} P, Mooij JM.
\newblock {M}arkov properties for graphical models with cycles and latent variables.
\newblock arXiv:171008775. 2017.
\newblock Available from: \url{https://arxiv.org/abs/1710.08775v1}.

\bibitem{richardson1997char}
Richardson T.
\newblock A characterization of {M}arkov equivalence for directed cyclic graphs.
\newblock International Journal of Approximate Reasoning. 1997;17(2):107-62.
\newblock Uncertainty in AI (UAI'96) Conference.
\newblock Available from: \url{https://www.sciencedirect.com/science/article/pii/S0888613X97000200}.

\bibitem{Frydenberg1990}
Frydenberg M.
\newblock The chain graph {M}arkov property.
\newblock Scandinavian Journal of Statistics. 1990;17(4):333-53.
\newblock Available from: \url{http://www.jstor.org/stable/4616181}.

\bibitem{Avis1992APA}
Avis D, Fukuda K.
\newblock A pivoting algorithm for convex hulls and vertex enumeration of arrangements and polyhedra.
\newblock Discrete \& Computational Geometry. 1992;8:295-313.
\newblock Available from: \url{https://doi.org/10.1007/BF02293050}.

\bibitem{githubrepo}
Semnani P, Robeva E. {Code for the paper "Causal structure learning in directed, possibly cyclic, graphical models"}. GitHub;.
\newblock Available from: \url{https://github.com/pardis-semnani/causal-discovery-cyclic-graphical-models}.

\bibitem{EDMONDS1977185}
Edmonds J, Giles R.
\newblock A min-max relation for submodular functions on graphs.
\newblock In: Hammer PL, Johnson EL, Korte BH, Nemhauser GL, editors. Studies in Integer Programming. vol.~1 of Annals of Discrete Mathematics. Elsevier; 1977. p. 185-204.
\newblock Available from: \url{https://www.sciencedirect.com/science/article/pii/S0167506008707349}.

\end{thebibliography}

\end{document}